\documentclass[10pt]{amsart}
\usepackage{amssymb,amscd,graphics,mathrsfs,latexsym,amsmath,amsthm,eucal}
\usepackage{a4wide}
\usepackage[matrix,arrow]{xy}

\usepackage{color}

\usepackage[dvips]{graphicx}
\usepackage{epsfig,enumerate}

\numberwithin{equation}{section}
\newtheorem{theorem}{Theorem}[section]
\newtheorem{proposition}[theorem]{Proposition}
\newtheorem{lemma}[theorem]{Lemma}
\newtheorem{corollary}[theorem]{Corollary}
\theoremstyle{definition}
\newtheorem{definition}[theorem]{Definition}
\theoremstyle{remark}

\newtheorem*{remark*}{Remark}

\DeclareMathOperator{\Tr}{Tr}

\newcommand{\lam}{\lambda}

\newcommand{\Hh}{\mathcal{H}}

\newcommand{\rz}{\mathbb R}
\newcommand{\zz}{\mathbb Z}

\newcommand{\cl}{\rm{cl}}

\def\tr{\text{tr}}
 \DeclareMathOperator{\RE}{Re}

\def\tr{{\rm tr}}
\newcommand{\mc}{\mathcal}
\newcommand{\rr}{\mathbb{R}}
\newcommand{\nn}{\mathbb{N}}
\newcommand{\cc}{\mathbb{C}}
\newcommand{\hh}{\mathbb{H}}

\newcommand{\eps}{\epsilon}

\newcommand{\pl}{\partial}
\newcommand{\x}{\times}

\newcommand{\til}{\widetilde}
\newcommand{\bbar}{\overline}

\newcommand{\cjd}{\rangle}
\newcommand{\cjg}{\langle}

\newcommand{\demi}{\frac{1}{2}}

\newcommand{\ddemi}{\frac{d}{2}}
\newcommand{\tra}{\textrm{Tr}}

\newcommand{\lb}{\textrm{lb}}
\newcommand{\rb}{\textrm{rb}}
\newcommand{\ff}{\textrm{ff}}

\newcommand{\indic}{\operatorname{1\negthinspace l}}

\newcommand{\stra}{\mathop{\hbox{\rm s-Tr}}\nolimits}

\def\sR{\hbox{I\kern-.1667em\hbox{R}}}

\def\tr{\hbox{Tr}}

\def\so{\mathrm{SO}}
\def\ad{\mathrm{ad}}

\def\spin{\mathrm{Spin}}

\def\tr{\mathrm{tr}}

\def\Tr{\mathrm{Tr}}

\begin{document}

\title[Eta invariant and Selberg zeta function of odd type]
{Eta invariant and Selberg Zeta function of odd type over convex co-compact
hyperbolic manifolds}
\author{Colin Guillarmou}
\address{Laboratoire J.A. Dieudonn\'e\\
Universit\'e de Nice Sophia-Antipolis\\
Parc Valrose, 06108 Nice\\France}
\email{cguillar@math.unice.fr}
\author{Sergiu Moroianu}
\address{Institutul de Matematic\u{a} al Academiei Rom\^{a}ne\\
P.O. Box 1-764\\
RO-014700 Bucharest, Romania}
\email{moroianu@alum.mit.edu}
\author{Jinsung Park}
\address{School of Mathematics\\ Korea Institute for Advanced Study\\
 Hoegiro 87\\ Dong\-daemun-gu\\ Seoul 130-722\\
Korea } \email{jinsung@kias.re.kr}


\date{\today}

\begin{abstract}
We show meromorphic extension and give a complete description of the
divisors of a Selberg zeta function of odd type
$Z_{\Gamma,\Sigma}^{\rm o}(\lambda)$ associated to the spinor
bundle $\Sigma$ on an odd dimensional convex co-compact hyperbolic
manifold $\Gamma\backslash\hh^{2n+1}$. As a byproduct we do a full
analysis of the spectral and scattering theory of the Dirac
operator on asymptotically hyperbolic manifolds. We show that
there is a natural eta invariant $\eta(D)$ associated to the Dirac
operator $D$ over a convex co-compact hyperbolic manifold
$\Gamma\backslash \hh^{2n+1}$ and that $\eta(D)=\frac{1}{\pi
i}\log Z_{\Gamma,\Sigma}^{\rm o}(0)$, thus extending Millson's
formula to this setting. We also define an eta invariant for the
odd signature operator, and we show that for Schottky
$3$-dimensional hyperbolic manifolds it gives the argument of a
holomorphic function which appears in the Zograf factorization
formula relating two natural K\"ahler potentials for
Weil-Petersson metric on Schottky space under some assumption on
the exponent of convergence of Poincar\'e series for the group
$\Gamma$.
\end{abstract}
\maketitle

\section{Introduction}

The eta invariant is a measure of the asymmetry of the spectrum
of self-adjoint elliptic operators which has been
introduced by Atiyah, Patodi and Singer
as the boundary term in the index formula for compact manifolds
with boundary \cite{APS}. For an elliptic self-adjoint pseudodifferential
operator $D$ of positive order acting on a bundle over a closed manifold,
it is defined as the value at $s=0$ of the meromorphic function
\[\eta(D,s):=\Tr\left(D (D^2)^{-\frac{s+1}2}\right)
=\frac{1}{\Gamma((s+1)/2)}\int_{0}^\infty t^{s-\demi}
\Tr\left(De^{tD^2}\right)dt,\]
which admits a meromorphic continuation from $\Re(s)\gg 0$ to $s\in\cc$
and is regular at $s=0$. In heuristic terms, $\eta(D):=\eta(D,0)$ computes
the asymmetry $\Tr(D|D|^{-1})$.

By applying Selberg's trace formula, Millson \cite{M} proved that
for any $(4m-1)$-dimensional closed hyperbolic manifold
$X_\Gamma:=\Gamma\backslash\hh^{4m-1}$, the eta invariant
$\eta(A)$ of the odd signature operator $A$ on odd forms
$\Lambda^{\mathrm{odd}}=\oplus_{p=0}^{2m}\Lambda^{2p-1}$ can be
expressed in terms of the geodesic flow on $X_\Gamma$. Millson
defined a Selberg zeta function of odd type by
\begin{equation}\label{defoddzeta}
Z^o_{\Gamma,\Lambda}(\lam):=\exp\left(-\sum_{\gamma\in
\mc{P}}\sum_{k=1}^\infty
\frac{\chi_+(R(\gamma)^k)-\chi_-(R(\gamma)^k)}{|\mathrm{det}({\rm
Id}-P(\gamma)^k)|^\demi}\frac{e^{-\lam k\ell(\gamma)}}{k}\right)
\end{equation}
where $\mc{P}$ denotes the set of primitive closed geodesics in
$X_\Gamma$, $R(\gamma)\in{\rm SO}(4m-2)$ is the holonomy along a
geodesic $\gamma$, $\chi_\pm$ denotes the character associated to
the two irreducible representations of ${\rm SO}(4m-2)$
corresponding to the $\pm i$ eigenspace of $\star$ acting on
$\Lambda^{2m-1}$, $P(\gamma)$ is the linear Poincar\'e map along
$\gamma$, and $\ell(\gamma)$ is the length of the closed geodesic
$\gamma$. Then he showed that
$Z^o_{\Gamma,\Lambda}(\lambda)$ extends meromorphically to
$\lam\in \cc$, its only zeros and poles occur on the line
${\Re(\lambda)=0}$ with order given in terms of the multiplicity
of the eigenvalues of $A$, and the following remarkable identity
holds:
\[e^{\pi i\eta(A)}=Z^o_{\Gamma,\Lambda}(0).\]
It is somehow believed that central value of Ruelle or Selberg type dynamical functions
have some kind topological meaning and this identity, as well as Fried's identity \cite{Fr},
provide striking examples.\\

It is a natural question to try to extend this
identity and to study the meromorphic extension and the zeros and
poles of the zeta function $Z^o_{\Gamma,\Lambda}(\lam)$ on
non-compact hyperbolic manifolds. The first step in this direction
has been done in \cite{P1} for cofinite hyperbolic quotients,
where Millson's identity holds with extra contributions from the
cusps, in the guise of the determinant of the scattering matrix.
In the present work, we carry out this program for \emph{convex
co-compact} manifolds, i.e., geometrically finite hyperbolic
manifolds with infinite volume and no cusps. In dimension
$3$ for the particular case of Schottky groups, our results have
interesting connections to Teichm\"uller theory, as we explain below.
\vspace{0.3cm}

The proof of the meromorphic extension of any reasonable
dynamical zeta function on co-compact hyperbolic manifolds
is contained in the work of Fried \cite{F} using transfer operator techniques
in dynamics, but as explained by
Patterson-Perry \cite{PP}, it extends to the convex co-compact setting
in a natural way. However there is no general description of
the zeros and poles, while we know in the co-compact and cofinite
cases that these are related to spectral and topological data since
the work of Selberg \cite{Sel}. There are now some rather recent
works of Patterson-Perry \cite{PP} and Bunke-Olbrich \cite{BO2}
which give a complete description of the zeros and poles of the
original Selberg zeta function on convex co-compact hyperbolic
(real) manifolds. The case of homogeneous vector bundles is not
yet completely described.

In this paper we study mainly the Dirac operator acting on the
spinor bundle $\Sigma$ over a convex co-compact hyperbolic
manifold $X_\Gamma:=\Gamma\backslash \hh^{2n+1}$. The basic
quantity associated with $\Gamma$ for this case is its exponent
$\delta_\Gamma$ defined to be the smallest number such that
\begin{equation}\label{e:delta}
\sum_{\gamma\in \Gamma} \exp({-(\lambda+n)r_\gamma}) < \infty
\end{equation}
for all $\lambda>\delta_\Gamma$. Here $r_\gamma$ denotes the
hyperbolic distance $d_{\hh^{2n+1}}(m,\gamma m)$ for a fixed point
$m\in \hh^{d+1}$. Note that our definition of $\delta_\Gamma$ is
shifted from the usual one by $-n$. For $\lambda>\delta_\Gamma$,
we define the Selberg zeta function of odd type
$Z^o_{\Gamma,\Sigma}(\lam)$ associated to the spinor bundle
$\Sigma$ exactly like in \eqref{defoddzeta} except that
$R(\gamma)$ denotes now the holonomy in the spinor bundle $\Sigma$
along $\gamma$, and $\chi_\pm$ denote the character of the two
irreducible representations of ${\rm Spin}(2n)$ corresponding to
the $\pm i$ eigenspaces of the Clifford multiplication ${\rm
cl}(T_\gamma)$ with the tangent vector field $T_\gamma$ to
$\gamma$. Like for the hyperbolic space $\hh^{2n+1}$, the Dirac
operator $D$ acting on the spinor bundle $\Sigma$ on a convex
co-compact hyperbolic manifold $X_\Gamma$ has continuous spectrum
the real line $\rr$, and one can define its resolvent for
$\Re(\lam)>0$ in two ways
\[ R_+(\lam):=(D+ i\lam)^{-1} , \textrm{ and }R_-(\lam):=(D-i\lam)^{-1}\]
as an analytic family of bounded operators acting on
$L^2(X_\Gamma,\Sigma)$. We then first show

\begin{theorem}\label{th1}
The Selberg zeta function of odd type $Z^o_{\Gamma,\Sigma}(\lam)$
associated to the spinor bundle $\Sigma$ on a convex co-compact
odd dimensional spin hyperbolic manifold
$X_\Gamma=\Gamma\backslash\hh^{2n+1}$ has a meromorphic extension
to $\cc$ and it is analytic in $\{\Re(\lam)\geq 0\}$. The
resolvents $R_\pm(\lam)$ of the Dirac operator have meromorphic
continuation to $\lam\in\cc$ when considered as operators mapping
$C_0^\infty(X_\Gamma,\Sigma)$ to its dual
$C^{-\infty}(X_\Gamma,\Sigma^*)$, and the poles have finite rank
polar part. A point $\lam_0\in \{\Re(\lam)<0\}$ is a zero or pole
of $Z^o_{\Gamma,\Sigma}(\lam)$ if and only if the meromorphic
extension of $R_+(\lam)$ or of $R_-(\lam)$ have a pole at
$\lam_0$, in which case the order of $\lam_0$ as a zero or pole of
$Z_{\Gamma,\Sigma}^o(\lam)$ (with the positive sign convention for
zeros) is given by
\[{\rm Rank }\, {\rm Res}_{\lam_0}R_-(\lam)-{\rm Rank }\,{\rm
Res}_{\lam_0}R_+(\lam).\]
\end{theorem}

In so far as analysis is concerned, we deal with a much more
general geometric setting in arbitrary dimensions and we prove
various results which were previously known for the Laplacian on
functions. We consider \emph{asymptotically hyperbolic manifolds}
(AH in short). These are complete Riemannian manifolds $(X,g)$
which compactify smoothly to a compact manifold with boundary
$\bar{X}$, whose metric near the boundary is of the form
$g=\bar{g}/x^2$ where $\bar{g}$ is a smooth metric on $\bar{X}$
and $x$ is any boundary defining function of $\pl\bar{X}$ in
$\bar{X}$, and finally such that $|dx|_{\bar{g}}=1$ at
$\pl\bar{X}$, a condition which is equivalent to assuming that the
curvature tends to $-1$ at the boundary. Convex co-compact
hyperbolic manifolds are special cases of AH manifolds. We show
that the spectrum of $D$ on AH manifolds is absolutely continuous
and given by $\rr$, that the resolvents $R_\pm(\lam)$ defined
above have meromorphic extensions to $\lam\in\cc$, and we define
the \emph{scattering operator} $S(\lam):
C^\infty(\pl\bar{X},\Sigma)\to C^{\infty}(\pl\bar{X},\Sigma)$ by
considering asymptotic profiles of generalized eigenfunctions on
the continuous spectrum. The family $S(\lam)$ extends to a
meromorphic family of elliptic pseudo-differential operators
acting on the boundary with the same principal symbol as
$D_{h_0}|D_{h_0}|^{2\lam-1}$ (up to a multiplicative constant),
where $D_{h_0}$ is the Dirac operator induced by the metric
$\bar{g}|_{T\pl\bar{X}}$. The scattering operator is a fundamental
object in the analysis of Selberg zeta function for convex
co-compact manifolds, and we study it thoroughly in this work. We
also show in a related note \cite{GMP} that the construction and
properties of the scattering operator have some nice applications,
for instance the invertibility of $S(\lam)$ except at discrete
$\lam$'s implies that the index of $D_{h_0}^+$ vanishes (the
so-called cobordism invariance of the index), and the operator
$\demi({\rm Id}-S(0))$ is the orthogonal Calder\'on projector of
the Dirac operator $\bar{D}$ on $\bar{g}$, providing a natural way
of constructing the Calder\'on projector without
extending or doubling the manifold $\bar{X}$.
\vspace{0.3cm}

In a second time, we prove that Millson's formula holds on
convex co-compact manifolds for the Dirac operator, and also for the
signature operator under some condition on $\delta_\Gamma$.

\begin{theorem}\label{th2}
Let $X_\Gamma=\Gamma\backslash\hh^{2n+1}$ be an odd
dimensional convex co-compact hyperbolic manifold. Then the
function ${\rm Tr}_{\Sigma}(De^{-tD^2})\in C^\infty(X_\Gamma)$ is
in $L^1(X_\Gamma)$ where ${\rm Tr}_{\Sigma}$ denotes the local
trace on the spinor bundle. The eta invariant $\eta(D)$ can
be defined as a convergent integral by
\begin{equation}\label{e:def-eta}
\eta(D):=\frac{1}{\sqrt{\pi}}\int_0^\infty
t^{-\demi}\Big(\int_{X_\Gamma}{\rm
Tr}_{\Sigma}(De^{-tD^2})(m)\,{\rm dv}(m)\Big)\, dt,
\end{equation}
and the following equality holds
\begin{equation}\label{e:Millson-two}
e^{\pi i\eta(D)}=Z^o_{\Gamma,\Sigma}(0).
\end{equation}
When $2n+1=4m-1$, the eta invariant $\eta(A)$ can also be defined
replacing $D$ by $A$ and $\Sigma$ by
$\Lambda^{\mathrm{odd}}=\oplus_{p=0}^{2m}\Lambda^{2p-1}$ in
\eqref{e:def-eta}, and if the exponent of convergence of Poincar\'e series $\delta_{\Gamma}$ is negative,
we moreover  have $e^{\pi i\eta(A)}=Z^o_{\Gamma,\Lambda}(0)$.
\end{theorem}

The assumption about $\delta_\Gamma$ for the equality $e^{\pi
i\eta(A)}=Z^o_{\Gamma,\Lambda}(0)$ is rather a technical condition
than a serious problem. Most of the analysis we do here for Dirac
operator $D$ goes through without significant difficulties to the
signature operator $A$, but it appears to be slightly more
involved essentially due to the fact that the continuous spectrum
of $A$ has two layers corresponding to closed and co-closed forms.
The complete analysis for forms in all dimensions will be included elsewhere.
\vspace{0.3cm}

To conclude this introduction and to motivate the eta invariant
$\eta(A)$ of the odd signature operator $A$, we describe the particular case of
Schottky $3$-dimensional manifolds with $\delta_\Gamma<0$, where
the eta invariant $\eta(A)$ can be considered as a function on the
space of deformations. The \emph{Schottky space} $\mathfrak{S}_g$
is the space of marked normalized Schottky groups with $g$
generators. It is a complex manifold of dimension $3g-3$, covering
the Riemann moduli space $\mathfrak{M}_g$ and with universal cover
the Teichm\"uller space $\mathfrak{T}_g$. It describes  the
deformation space of the $3$-dimensional hyperbolic Schottky
manifolds $X_\Gamma=\Gamma\backslash\hh^3$ as well as the one of
boundary Riemann surfaces. Like $\mathfrak{T}_g$, the Schottky
space $\mathfrak{S}_g$ has a natural K\"ahler metric, the
Weil-Petersson metric. In \cite{ZT,ZT2}, Takhtajan-Zograf
constructed two K\"ahler potentials of the Weil-Petersson metric
on $\mathfrak{S}_g$, that is,
\[
\partial \overline{\partial} S= \partial\overline{\partial}\Big(-12\pi \log
\frac{\mathrm{Det}\Delta}{\mathrm{det}\, \mathrm{Im}\, \tau}\Big)=2i\, \omega_{WP}
\]
where $\partial$ and $\overline{\partial}$ are the $(1,0)$ and
$(0,1)$ components of the de Rham differential $d$ on
$\mathfrak{S}_g$ respectively, and $\omega_{WP}$ is the symplectic
form of the Weil-Petersson metric; here
$S$ is the so-called \emph{classical Liouville action}, $\mathrm{Det}\Delta$ is the $\zeta$-regularized
determinant of the Laplacian $\Delta$ on the Riemann surface
determined by the hyperbolic metric (i.e. a point of $\mathfrak{S}_g$) and $\tau$ denotes a period matrix.
We show that
\begin{theorem}\label{th3}
The function ${F}$ defined on $\mathfrak{S}_g^0:=\{\Gamma\in
\mathfrak{S}_{g}; \delta_{\Gamma}<0\}$ by
\[F:= \frac{\mathrm{Det}\Delta}{\mathrm{det}\, \mathrm{Im}\, \tau}\,\exp\left(\frac{S}{12\pi}-i\pi \eta(A)\right)\]
is holomorphic. In particular, the eta invariant $\eta(A)$ is a pluriharmonic
function on $\mathfrak{S}_g^0$.
\end{theorem}

The condition $\delta_\Gamma<0$ in Theorem \ref{th3} simplifies
the proof at several stages. But, one can expect that a similar
result still holds over the whole Schottky space $\mathfrak{S}_g$.
This extension problem will be discussed elsewhere.
\vspace{0.3cm}

\textbf{Acknowledgements} The authors thank the KIAS Seoul
(South Korea) where part of this work was started and completed. We
are also grateful to the CIRM in Luminy for funding us through 
a `Research in Pairs' in January 2008. 
C.G.\ was supported by grants ANR05-JCJC-0107091 and NSF0500788, he thanks the MSRI at Berkeley 
where part of this work was done. 
S.M.\  was supported by grant PN-II-ID-PCE 1188 265/2009.

\section{The Dirac operator on real hyperbolic space}

\subsection{Dirac operators over hyperbolic spaces} The
$(d+1)$-dimensional real hyperbolic space is the manifold
$$
\mathbb{H}^{d+1}=\big{\{}\, x\in \mathbb{R}^{d+2}\, | \,
x_0^2+x_1^2+\ldots+x_d^2-x_{d+1}^2=-1,\ x_{d+1}>0 \, \big{\}}
$$
equipped with the metric of curvature $-1$. The orientation
preserving isometries of $\hh^{d+1}$ form the group $\so(d+1,1)$. The
isotropy subgroup of the base point $(0,\ldots,0,1)$ is isomorphic
to $\so(d+1)$. Hence the real hyperbolic space $\hh^{d+1}$ can be
identified with the symmetric space $\so(d+1,1)/\so(d+1)$. Since
$G=\spin(d+1,1)$, $K=\spin(d+1)$ are double coverings of $\so(d+1,1)$,
$\so(d+1)$ respectively, we see that $\so(d+1,1)/\so(d+1)=G/K$ and we use the
identification $\hh^{d+1}\cong G/K$ for our purpose. We denote
the Lie algebras of $G$, $K$ by
$\mathfrak{g}=\mathfrak{spin}(d+1,1)$, $\mathfrak{k}=
\mathfrak{spin}(d+1)$ respectively. The Cartan involution $\theta$
on $\mathfrak{g}$ gives us the decomposition
$\mathfrak{g}=\mathfrak{k}\oplus\mathfrak{p}$ where
$\mathfrak{k},\mathfrak{p}$ are the $1,-1$ eigenspaces of
$\theta$, respectively. The subspace $\mathfrak{p}$ can be
identified with the tangent space $T_o(G/K)\cong
\mathfrak{g}/\mathfrak{k}$ at $o=eK\in G/K$. The invariant metric
of curvature $-1$ over $\hh^{d+1}$ is given by the normalized
Cartan-Killing form
\begin{equation}\label{e:CK}
\langle X, Y\rangle:=-\frac{1}{2d} C(X, \theta Y)
\end{equation}
where the Killing form is defined by $C(X,Y)=\Tr (\ad\, X \circ
\ad\, Y)$ {for} $X,Y\in \mathfrak{g}$.

Let $\mathfrak{a}$ be a fixed maximal abelian subspace of
$\mathfrak{p}$. Then the dimension of $\mathfrak{a}$ is $1$. Let
$M= \spin(d)$ be the centralizer of $A=\exp(\mathfrak{a})$ in
$K$ with Lie algebra $\mathfrak{m}$. We put $\beta$ to be {the}
positive restricted root of $(\mathfrak{g},\mathfrak{a})$. Let
$\rho$ denote the half sum of the positive roots of
$(\mathfrak{g},\mathfrak{a})$, that is,
$\rho=\frac{d}{2}\beta$. From now on, we use the
identification
\begin{equation}\label{e:ident}
\mathfrak{a}_{\mathbb{C}}^* \ {\cong}\ \mathbb{C}
\qquad\text{by}\quad \lambda\, \beta  \longrightarrow  \lambda.
\end{equation}
Let $\mathfrak{n}$ be the positive root space of $\beta$ and
$N=\exp(\mathfrak{n})\subset G$. The Iwasawa decomposition is
given by $G=KAN$. Throughout this paper we use the following Haar
measure on $G$,
\begin{equation}\label{e:measure}
dg= a^{2\rho} dk\, da\, dn = a^{-2\rho} dn\, da\, dk
\end{equation}
where $g=kan$ is the Iwasawa decomposition and
$a^{2\rho}=\exp(2\rho(\log a))$. Here $dk$ is the Haar measure
over $K$ with $\int_K dk=1$, $da$ is the Euclidean Lebesgue
measure on $A$ given by the identification $A \cong \mathbb{R}$
via $a_r=\exp(rH)$ with $H\in\mathfrak{a}$, $\beta(H)=1$, and $dn$
is the Euclidean Lebesgue measure on $N$ induced by the normalized
Cartan-Killing form $\langle \cdot, \cdot\rangle$ given in
\eqref{e:CK}.

The spinor bundle $\Sigma(\hh^{d+1})$ can be identified with  the
associated homogeneous vector bundle over $\hh^{d+1}= G/K$ with the
spin representation $\tau_{d}$ of $K\cong \spin(d+1)$
acting on $V_{\tau_d}=\mathbb{C}^{2^{[d+1/2]}}$, that is,
\begin{equation}\label{e:homog}
\Sigma (\hh^{d+1})= G\times_{\tau_d} V_{\tau_d} \ \longrightarrow \
\hh^{d+1}= G/K.
\end{equation}
Here points of $G\times_{\tau_d} V_{\tau_d}$ are given by
equivalence classes $[g,v]$ of pairs $(g,v)$ under $(gk,v)\sim
(g,\tau_d(k)v)$. Hence the sections of $G\times_{\tau_d}
V_{\tau_d}$ from $G/K$ consist of functions $f:G\to V_{\tau_d}$
with the $K$-equivariant condition,
\[
f(gk)=\tau_d(k)^{-1} f(g)
\]
for $g\in G$, $k\in K$.
Recall that $\tau_d$ is irreducible if $d+1$ is odd, while it splits into $2$ irreducible representations
if $d+1$ is even.

Let us denote by
\[
\nabla: C^\infty(\hh^{d+1},\Sigma(\hh^{d+1})) \ \longrightarrow \
C^\infty(\hh^{d+1},T^*(\hh^{d+1})\otimes \Sigma(\hh^{d+1}))
\]
the covariant derivative induced by the lift of the Levi-Civita
connection to the spinor bundle $\Sigma(\hh^{d+1})$, and by ${{\rm
cl}}:T_m(\hh^{d+1})\to \mathrm{End}\, \Sigma_m(\hh^{d+1})$ the
Clifford multiplication. Then the Dirac operator
${D}_{\hh^{d+1}}$ acting on
$C^\infty_0(\hh^{d+1},\Sigma(\hh^{d+1}))$ is defined by
\[
{D}_{\hh^{d+1}} f(m)= \sum_{j=1}^d {\rm cl}(e_j) \nabla _{e_j}
f(m)\hspace{1cm} \text{for} \quad f\in
C^\infty_0(\hh^{d+1},\Sigma(\hh^{d+1}))
\]
where $(e_j)_{j=1}^{d+1}$ denotes an orthonormal frame of
$T_m(\hh^{d+1})$. The Dirac operator ${D}_{\hh^{d+1}}$ is an
essentially self-adjoint, elliptic and $G$-invariant differential
operator of first order, and we use the same notation for its self
adjoint extension to $L^2(\hh^{d+1}, \Sigma(\hh^{d+1}))$. It is
well known that the spectrum of ${D}_{\hh^{d+1}}$ on
$L^2(\hh^{d+1}, \Sigma(\hh^{d+1}))$ consists only of the
absolutely continuous spectrum $\mathbb{R}$ (for instance, by
Corollary 4.11 in \cite{CP}).

\subsection{The resolvent on $\hh^{d+1}$}

Let us define the resolvent in the half plane $\{\Re(\lam)>0\}$ by
\[R_{\hh^{d+1}}(\lam):=({D}_{\hh^{d+1}}^2+\lam^2)^{-1}\]
which maps $L^2(\hh^{d+1},\Sigma(\hh^{d+1}))$ to itself. Recall
the hypergeometric function $F(a,b,c,z)$ defined by
\[
F(a,b,c,z)=\sum_{k=0}^\infty \frac{\Gamma(a+k)\Gamma(b+k)\Gamma(c)}{\Gamma(a)\Gamma(b)\Gamma(c+k)} z^k
\qquad \text{for} \quad |z|<1.
\]
Then we have from the work of Camporesi \cite[Th. 6.2 and 6.3]{CP}

\begin{proposition}{\bf [Camporesi]}\label{p:CP}
For $\Re(\lam)>0$, the respective Schwartz kernels of $R_{\hh^{d+1}}(\lam)$
and ${D}_{\hh^{d+1}}R_{\hh^{d+1}}(\lam)$ are given by
\begin{equation}\label{resolvforHn}\begin{split}
R_{\hh^{d+1}}(\lam;m,m')=& 2^{-(d+1)}\pi^{-\frac{d+1}{2}}
\frac{\Gamma(\frac{d+1}{2}+\lambda)\Gamma(\lambda)}{\Gamma(2\lambda+1)}(\cosh
(r/2))^{-2(\frac{d}{2}+\lambda)}\\
&\times F(\frac{d+1}{2}+\lambda,\lambda,2\lambda+1,\cosh^{-2}(r/2))\,
U(m,m'),
\end{split}\end{equation}
\begin{equation}\label{e:CP} \begin{split}
{D}_{\hh^{d+1}}R_{\hh^{d+1}}(\lam;m,m')=& 
-2^{-(d+1)}\pi^{-\frac{d+1}{2}}\frac{\Gamma(\frac{d+1}{2}+\lambda)
\Gamma(\lambda+1)}{\Gamma(2\lambda+1)} 
(\cosh (r/2))^{-(d+1)-2\lambda}\sinh(r/2)\\
&\times \, F\big( \frac{d+1}{2}+\lambda, \lambda+1, 2\lambda+1,
\cosh^{-2}(r/2)){\rm cl}(v_{m,m'}) U(m,m')
\end{split}\end{equation}
where $r=d_{\hh^{d+1}}(m,m')$ for $m,m'\in \hh^{d+1}$,
$v_{m,m'}$ is the unit tangent vector at $m$ to the
geodesic from $m'$ to $m$ and $U(m,m')$ is the parallel transport
from $m'$ to $m$ along the geodesic between them. Moreover
$R_{\hh^{d+1}}(\lam)$ has an analytic continuation in
$\cc\setminus\{0\}$ with a simple pole at $\lam=0$ and
${D}_{\hh^{d+1}}R_{\hh^{d+1}}(\lam)$ admits an analytic
continuation to $\lam\in\cc$, as distributions on
$\hh^{d+1}\x\hh^{d+1}$.
\end{proposition}

\begin{remark*} If one denotes
$R_{\hh^{d+1}}(\lam;m;m')=J_\lam(r)U(m,m')$ where
$r=d_{\hh^{n+1}}(m,m')$, the function $J_\lam(r)$ satisfies that
$J_{\lam}(r)-J_{-\lam}(r)$ is smooth in $r$ near $r=0$. This can
be checked using functional equations of hypergeometric functions
but actually follows directly from elliptic regularity since
$J_\lam(r)-J_{-\lam}(r)$ (since the difference of
resolvents too) solves an elliptic ODE. The kernel $D\Pi(\lam;m,m')$ of
$D(R_{\hh^{d+1}}(\lam)-R_{\hh^{d+1}}(-\lam))$ is then also
smooth near the diagonal $m=m'$ and following the proof of
\cite[Th 6.3]{CP}, we see that it can be written under the form
\[D\Pi(\lam;m,m')=-\demi\sinh(r)
\frac{\pl_{\cosh^{-2}(r/2)}\Big((\cosh^{-2}
(\frac{r}{2}))^{-\frac{d}{2}}H_\lam(\cosh^{-2}(\frac{r}{2}))\Big)}{\cosh(\frac{r}{2})^{d+4}}{\rm cl}(v_{m,m'})U(m,m')\]
where $H_{\lam}(\cosh^{-2}(r/2)):=J_\lam(r)-J_{-\lam}(r)$ with $H_\lam(u)$ smooth near $u=1$. Then we deduce
that on the diagonal $D\Pi(\lam;m,m)=0$.
\end{remark*}

\subsection{Dirac operators over convex co-compact hyperbolic
manifolds}\label{ss-hyp} Let $\Gamma$ denote a convex co-compact
torsion-free discrete subgroup of $G= \spin(d+1,1)$ such that its
co-volume $\mathrm{Vol}(\Gamma\backslash G)=\infty$. Hence
\[X_\Gamma:=\Gamma\backslash G/K\]
is a $(d+1)$-dimensional convex co-compact hyperbolic manifold of
infinite volume, which is a spin manifold by construction.  The
boundary $\partial \hh^{d+1}$, which can be identified with $K/M$,
admits a $\Gamma$-invariant decomposition into
$\Omega(\Gamma)\cup\Lambda(\Gamma)$ where $\Omega(\Gamma)\neq
\varnothing$ is open and $\Gamma$ acts freely and co-compactly on
$\hh^{d+1}\cup \Omega(\Gamma)$. Hence $X_\Gamma$ can be
compactified by adjoining the geodesic boundary $\Gamma\backslash
\Omega(\Gamma)$.

By the identification \eqref{e:homog} of the spinor bundle
$\Sigma(\hh^{d+1})$ with the homogeneous vector bundle
$G\times_{\tau_d} V_{\tau_d}$, we can also identify the spinor
bundle $\Sigma(X_\Gamma)$ over $X_\Gamma$ with the locally
homogeneous vector bundle $\Gamma\backslash \big( G\times_{\tau_d}
V_{\tau_d}\big)$. Here $\Gamma$ acts on
$G\times_{\tau_d}V_{\tau_d}$ by $\gamma [g,v]=[\gamma g, v]$ for
$\gamma\in\Gamma$. We can also push down the Dirac operator
${D}_{\hh^{d+1}}$ to $X_\Gamma$, which we denote by $D$. We also use
the same notation for its unbounded self-adjoint extension in
$L^2(X_\Gamma,\Sigma(X_\Gamma))$, that is,
\begin{equation*}
D: \ L^2(X_\Gamma,\Sigma(X_\Gamma)) \, \longrightarrow \,
L^2((X_\Gamma,\Sigma(X_\Gamma)).
\end{equation*}
By Corollary \ref{spectrum} (Corollary 7.9 and Theorem 11.2 in
\cite{BO3}), the Dirac operator $D$ over
$L^2(X_\Gamma,\Sigma(X_\Gamma))$ has no discrete spectrum and only
absolutely continuous spectrum $\mathbb{R}$.

\section{Resolvent of Dirac operator on Asymptotically Hyperbolic manifolds}

In this section, we analyze the resolvent $R(\lam)$ of $D^2$ on an
asymptotically hyperbolic manifold (AH in short) of dimension
$(d+1)$. An asymptotically hyperbolic manifold is a complete
non-compact Riemannian manifold $(X,g)$ which compactifies in a
smooth manifold with boundary $\bar{X}$ and there is a
diffeomorphism $\psi$ (called \emph{product decomposition}) from a
collar neighbourhood $[0,\eps)_x\x\pl\bar{X}$ of the boundary to a
neighbourhood of $\pl\bar{X}$ in $\bar{X}$ so that
\begin{equation}\label{metricg}
\psi^*g=\frac{dx^2+h_x}{x^2}
\end{equation}
for some one-parameter family of metrics $h_x$ on the boundary
$\pl\bar{X}$ depending smoothly on $x\in[0,\eps)$. By abuse of
notations, we will write $x$ for $\psi_*x$, and $x$ is then a
boundary defining function in $\bar{X}$ near $\pl\bar{X}$,
satisfying $|dx|_{x^2g}=1$. A boundary defining functions
satisfying $|dx|_{x^2g}=1$ near the boundary is called
\emph{geodesic boundary defining function}, and it yields a
diffeomorphism $\psi$ like in \eqref{metricg} by taking the flow
of the gradient $\nabla^{x^2g}x$ starting at the boundary.
Following the terminology of \cite{GuiDMJ}, we shall say that
\begin{equation}\label{evenness}
\textrm{ the metric is } even \textrm{ if the Taylor expansion of }h_x\textrm{ at }x=0\textrm{
contains only even powers of }x.
\end{equation}
This property does not depend on the choice of the diffeomorphism $\psi$ but only on $g$, see
\cite[Lemma 2.1]{GuiDMJ}. It is well known that convex co-compact
quotients $X_\Gamma=\Gamma\backslash\hh^{d+1}$ are even AH manifolds
(see \cite{MM}). Note that the metric $h_0$ is not canonical since
it depends on the choice of $\psi$, but its conformal class $[h_0]$ is
canonical with respect to $g$.

\subsection{0-structures, spinor bundle and Dirac operator}

Following the ideas of Mazzeo-Melrose \cite{MM} (and refer to this
paper for more details), there is a natural structure associated
to AH manifolds, this is encoded in the Lie algebra
$\mc{V}_0(\bar{X})$ of smooth vector fields vanishing at the
boundary, whose local basis over $C^{\infty}(\bar{X})$ is given
near the boundary $\pl\bar{X}$ by the vector fields $(x\pl_x,
x\pl_{y_1},\dots, x\pl_{y_{d}})$ if $(x,y_1,\dots,y_{d})$ is a
local chart near a point $p\in \pl\bar{X}$ and $x$ is a smooth
boundary defining function in $\bar{X}$. The algebra is also the
space of smooth section of a bundle $^0T\bar{X}$ with local basis
near $p$  given by $(x\pl_x, x\pl_{y_1},\dots, x\pl_{y_{d}})$
and its dual space is denoted $^0T^*\bar{X}$, with local basis
$(dx/x,dy_1/x,\dots,dy_{d}/x)$. The metric $g$ is a smooth
section of the bundle
of positive definite symmetric form $S^2_+(^0T^*\bar{X})$ of $^0T^*\bar{X}$.
\vspace{0.3cm}

Let us define $\bar{g}:=x^2g$ where $x$ is a boundary defining
function appearing in $\eqref{metricg}$. If $(\bar{X},\bar{g})$ is
orientable, there exists an ${\rm SO}(d+1)$-bundle
${_o}F(\bar{X})\to \bar{X}$ over $\bar{X}$, but also an ${\rm
SO}(d+1)$-bundle ${^0_o}F(\bar{X})\to \bar{X}$ defined using the
$0$-tangent bundle $^0T\bar{X}$ and the metric $g$ smooth on it.
If $(\bar{X},\bar{g})$ admits a spin structure, then there exists
a $0$-spin structure on $(X,g)$ in the sense that there is a ${\rm
Spin}(d+1)$-bundle ${^0_s}F(\bar{X})\to \bar{X}$ which double
covers ${^0_o}F(\bar{X})$ and is compatible with it in the usual
sense. This corresponds to a rescaling of the spin structure
related to $(\bar{X},\bar{g})$. The $0$-spinor bundle
$^0\Sigma(\bar{X})$ can then be defined as a bundle associated to
the ${\rm Spin}(d+1)$ principal bundle ${^0_s}F(\bar{X})$, with
fiber at $p\in\bar{X}$
\[^0\Sigma_p(\bar{X})=({^0_s}F_p\x V_{\tau_d})/\tau_d. \]
The vector field $x\pl_x:=x\nabla^{\bar{g}}(x)$ in the collar neighbourhood is unit normal to all hypersurfaces $\{x={\rm constant}\}$.
The $0$-spinor bundle on $\bar{X}$ splits near the boundary under the form
\[^0\Sigma={^0\Sigma}_+\oplus {^0\Sigma}_- , \quad
\textrm{where } {^0\Sigma}_\pm:= \ker({\rm cl}(x\pl_x)\mp i),\]
note that this
splitting is dependent of the choice of the geodesic boundary
defining function $x$ except at the boundary $\pl\bar{X}$ where it
yields an independent splitting of the spinor (since the first jet
of $x\nabla^{x^2g}x$ is independent of $x$ at $\pl\bar{X}$). To
avoid confusions later (and emphasize the fact that it is only
depending on the conformal class $(\pl\bar{X},[h_0])$), we shall
define ${\rm cl}(\nu)$ the linear map on
$^0{\Sigma}|_{\pl\bar{X}}$
\[{\rm cl}(\nu) \psi:= {\rm cl}(x\pl_x)\psi.\]
At the boundary, ${^0\Sigma}|_{\pl\bar{X}}$ is diffeomorphic to
the spinor bundle $\Sigma(\pl\bar{X})$ on $(\pl\bar{X},h_0)$, this
is not canonical since it depends on $h_0$ and thus on the choice
of $x$, however the splitting above is. Notice
also that in even dimension, the splitting ${^0\Sigma}_+\oplus {^0\Sigma}_-$ near the boundary
is not usual splitting of the spinor bundle.
The Dirac operator near the boundary has the form
\begin{equation}\label{diracAH}
D=x^{\ddemi}\left( {\rm cl}(x\pl_x)x\pl_x + P
\right)x^{-\ddemi}, \quad P\in {\rm Diff}_0^1(\bar{X};\mc{E})
\end{equation}
where $P$ is a first order differential operator in tangential derivatives, which anticommutes with ${\rm cl}(x\pl_x)$ 
and such that $P=xD_{h_0}+O(x^2)$ where $D_{h_0}$ is the Dirac operator on $\pl\bar{X}$ equipped with the metric $h_0$.
If the metric $g$ is even, it is easy to see that locally near any point $y'$ of the boundary, if 
$(x\pl_x,xY_1,\dots,xY_{d})$ is an orthonormal frame near $y'$ and $(x,y)$ are coordinates on $[0,\eps)\x\pl\bar{X}$ there, 
then $P$ is of the form
\[P=\sum_{i=1}^d P_i(x^2,y;\nabla^{\bar{g}}_{xY_i})\]
for some differential operators $P_i$ of order $1$ and with smooth coefficients in
$(x^2,y)$. This can be checked for instance by using that $D=x^{\ddemi+1}\bar{D}x^{-\ddemi}$ where $\bar{D}$ is the Dirac
operator for the metric $\bar{g}:=x^2g$ which is smooth in the coordinates $({\bf x}=x^2,y)$ down to ${\bf x}=0$.
From these properties, it is straightforward to check that if $g$ is even, then for $x$ geodesic boundary defining function 
fixed, $D$ preserves the space
$\mc{A}_\pm\subset C^{\infty}(\bar{X},^0\Sigma)$ of smooth spinors which have expansion at the boundary of the form
\begin{equation}\label{defA}
\sigma\sim_{x\to 0} \sum_{j=0}^\infty x^j\psi_j, \textrm{ with }\psi_{2j}\in \Sigma_\pm(\pl\bar{X}) \textrm{ and }
\psi_{2j+1}\in \Sigma_\mp(\pl\bar{X}).\end{equation}

\subsection{The stretched product}\label{stretched}
Following Mazzeo-Melrose \cite{MM}, we define the  stretched
product $\bar{X}\x_0\bar{X}$ as the blow-up
$[\bar{X}\x\bar{X},\Delta_\pl]$ of $\bar{X}\x\bar{X}$ around the
diagonal in the boundary $\Delta_\pl:=\{(y,y)\in
\pl\bar{X}\x\pl\bar{X}\}$. The blow-up is a smooth manifold with
codimension $2$ corners, and $3$ boundary hypersurfaces, the left
boundary denoted $\lb$, the right boundary denoted $\rb$ and the
new face, called `front face' and denoted $\ff$, obtained from the
blow-up. The blow-down map is denoted $\beta:\bar{X}\x_0\bar{X}\to
\bar{X}\x\bar{X}$ and maps ${\rm int}(\lb)$ to $\pl\bar{X}\x X$,
${\rm int}(\rb)$ to $X\x\pl\bar{X}$ and $\ff$ to $\Delta_{\pl}$.
The face $\ff$ is a bundle over $\Delta_\pl\simeq \pl\bar{X}$ with
fibers a quarter of $d$-dimensional sphere. Let us use the
boundary defining function $x$ in \eqref{metricg}, which induces
$x:=\pi_L^*x$ and $x':=\pi_R^*x$ as boundary defining functions of
$\bar{X}\x\bar{X}$ where $\pi_R,\pi_L$ are the right and left
projection $\bar{X}\x \bar{X}\to\bar{X}$. The fibre $\ff_p$ of the
front face $\ff$ (with $p=(y',y')\in \pl\bar{X}\x \pl\bar{X}$) is,
by definition of blow-up, given by the quotient
\begin{equation}\label{fp}
\ff_p=
\Big(\Big(N_p(\Delta_{\pl},\pl\bar{X}\x\pl\bar{X})\x(\rr^+\pl_x)
\x (\rr^+\pl_{x'})\Big) \setminus\{0\}\Big) \Big{/}\{(w,t,u)\sim
s(w',t',u'), s>0\}\end{equation} where in general $N(M,Y)$ denotes
the normal bundle of a submanifold $M$ in a manifold $Y$. Since
$T_{y'}\pl\bar{X}$ is canonically isomorphic to
$N_p(\Delta_\pl,\pl\bar{X}\x \pl\bar{X})$ by $z\in
T_{y'}\pl\bar{X}\to (z,-z)\in T_p(\pl\bar{X}\x \pl\bar{X})$,
$h_0(y')$ induces a metric on $N_p(\Delta_\pl,\pl\bar{X}\x
\pl\bar{X})$. Then $\ff_p$ is clearly identified with the quarter
of sphere
\[\ff_p\simeq \{w+t\pl_x+u\pl_{x'}\in 
N_p(\Delta_\pl,\pl\bar{X}\x \pl\bar{X})\x(\rr^+\pl_x) \x 
(\rr^+\pl_{x'}), t^2+u^2+|w|_{h_0(y')}^2=1\}.\]
In projective coordinates $(s:=t/u,z:=w/u)\in(0,\infty)\x \rr^{d}$,
the interior of the front face fiber $\ff_p$ is a half-space diffeomorphic to $\hh^{d+1}$.
In the same way we define the blow-up $\bar{X}\x_0\pl\bar{X}$ of $\bar{X}\x\pl\bar{X}$ around
$\Delta_\pl$ and the blow-up $\pl\bar{X}\x_0\pl\bar{X}$ of $\pl\bar{X}\x\pl\bar{X}$ around $\Delta_\pl$.
The first one is canonically diffeomorphic to the face ${\rm rb}$ of $\bar{X}\x_0\bar{X}$ while the second
one is canonically diffeomorphic to ${\rm lb}\cap{\rm rb}$.

The manifold $\bar{X}\x\bar{X}$ carries the bundle
\[\mc{E}={^0\Sigma}(\bar{X})\boxtimes {^0\Sigma}^*(\bar{X})\]
which on the diagonal is isomorphic to ${\rm End}({^0\Sigma})$.
This bundle lifts under $\beta$ to a bundle over $\bar{X}\x_0\bar{X}$, still denoted by $\mc{E}$, whose fiber at the front face $\ff$ is given
by ${^0\Sigma}_{y'}(\bar{X})\boxtimes {^0\Sigma}_{y'}^*(\bar{X})$ everywhere on the fiber $\ff_p$ (here $p=(y',y')\in\Delta_\pl$)
if ${^0\Sigma}_{y'}(\bar{X})$ is the fiber of ${^0\Sigma}(\bar{X})$ at the point $y'\in \pl\bar{X}$.

On a manifold with corners $M$ with a smooth bundle $E\to M$, let
us denote by $\dot{C}^\infty(M,E)$ the space of smooth section of
$E$ which vanish to all order at the (topological) boundary and
let $C^{-\infty}(M,E^*)$ be its dual, the elements of which are
called extendible distributions. Then $\beta^*$ is an isomorphism
between $\dot{C}^{\infty}(\bar{X}\x\bar{X};\mc{E})$
and $\dot{C}^\infty(\bar{X}\x_0\bar{X};\mc{E})$ and also between their duals, meaning that
distributions on $\bar{X}\x\bar{X}$ can be as well considered on
the stretched product. In what follows, we consider the Schwartz
kernel $K_A\in C^{-\infty}(\bar{X}\x\bar{X},\mc{E})$ of an operator $A:
\dot{C}^\infty(\bar{X},{^0\Sigma})\to
C^{-\infty}(\bar{X},{^0\Sigma})$ defined by
\[ \cjg A\psi,\phi\cjd = \cjg K_A,\phi\boxtimes\psi\cjd\]
where $\cjg.,.\cjd$ is the duality pairing using the volume density of the metric.
By abuse of notations we will write $A(m,m')$ for $K_A(m,m')$ and the bundle
$\mc{E}$ at the diagonal will be identified to ${\rm End}({^0\Sigma})$.

\subsection{Pseudo-differential operators}

We define the space $\Psi_0^{m,\alpha,\beta}(\bar{X};\mc{E})$ for $m\in\rr$, $\alpha,\beta\in\cc$ as in \cite{MaJDG,MM}, and refer the reader
to these references for more details. An operator $A$ is in $\Psi_0^{m,\alpha,\beta}(\bar{X};\mc{E})$ if its Schwartz kernel $K_A$ lifts to $\bar{X}\x_0\bar{X}$ to a distribution $\beta^*(K_A)$ which can be decomposed as a sum $K^1_A+K^2_A$ with $K^1_A\in \rho_{\rm lb}^{\alpha}\rho_{\rm lb}^\beta C^{\infty}(\bar{X}\x_0\bar{X};\mc{E})$ and $K^2_A\in I_{\rm cl}^{m}(\bar{X}\x_0\bar{X},\Delta; \mc{E})$ where $I_{\rm cl}^{m}(\bar{X}\x_0\bar{X},\Delta; \mc{E})$ denotes the space of distribution on $\bar{X}\x_0\bar{X}$ classically conormal to the lifted diagonal $\Delta:=\bbar{\beta^*(\{(m,m)\in X\x X\})}$
of order $m$ and vanishing to infinite order at the left and right boundaries ${\rm lb}\cup{ \rm rb}$.

\subsection{Microlocal structure of the resolvent on $\hh^{d+1}$}
We want to describe the resolvent kernel as a conormal
distribution on a compactification of $\hh^{d+1}\x\hh^{d+1}$, in
order to show later that a similar result holds for convex
co-compact quotients and more generally asymptotically hyperbolic
manifolds. Here we let $\bbar{\hh^{d+1}}$ be the natural
compactification of $\hh^{d+1}$, i.e. the unit ball in
$\rr^{d+1}$.

\begin{lemma}\label{resolvhd}
The analytically extended resolvent $R_\pm^{\hh^{d+1}}(\lam):= (D_{\hh^{d+1}}\pm
i\lambda)^{-1}$ of the Dirac operator $D_{\hh^{d+1}}$ on
$\hh^{d+1}$ is in the space $\Psi_0^{-1,\lam+\ddemi,\lam+\ddemi}(\bbar{\hh^{d+1}};\mc{E})$.
\end{lemma}
\begin{proof} This is very similar to the case of the Laplacian on functions dealt with in \cite{MM}, since we have
$\cosh(r/2)^{-2}\in
\rho_{\lb}\rho_{\rb}C^{\infty}(\bbar{\hh^{d+1}}\x_0\bbar{\hh^{d+1}}\setminus
\Delta)$ and is a smooth function of the distance, and since, by
the remark after Lemma \ref{ptp}, the lift of $U(m,m)$ by
$\beta^*$ is smooth on $\bbar{\hh^{d+1}}\x_0\bbar{\hh^{d+1}}$,
combining the formulae \eqref{resolvforHn}, \eqref{e:CP}
 proves our claim, except for the singularity
at the diagonal $\Delta$. The conormal diagonal singularity can
easily seen by applying the first step of the parametrix of
Mazzeo-Melrose (we refer to the proof of Proposition
\ref{extresolv} below) near the diagonal, indeed the construction
shows that there exists
$Q^0_{\pm}(\lam)\in\Psi_0^{-1,\infty,\infty}(\bbar{\hh^{d+1}};\mc{E})$
such that $(D_{\hh^{d+1}}\pm
i\lam)(R^{\hh^{d+1}}_\pm(\lam)-Q^0_\pm(\lam))$ and
$(R^{\hh^{d+1}}_\pm(\lam)-Q^0_\pm(\lam))(D_{\hh^{d+1}}\pm i\lam)$
have a smooth kernel in a neighbourhood of $\Delta$ down to the
front face $\ff$, and so by $0$-elliptic regularity
$R_\pm^{\hh^{d+1}}(\lam)-Q^0_\pm(\lam)$ is smooth near $\Delta$
down to $\ff$.
\end{proof}

Note that the resolvent can be also considered as a convolution
kernel on $\hh^{d+1}$ with a conormal singularity at the center
$0\in \hh^{d+1}$.

\subsection{The parametrix construction of Mazzeo-Melrose}
We can construct the resolvent $R_\pm(\lam):=(D\pm i\lam)^{-1}$
through a pseudo-differential parametrix, following Mazzeo-Melrose
\cite{MM} or Mazzeo \cite{MaJDG}. We will not give the full
details since this is a straightforward application of the paper
\cite{MM} and the analysis of the resolvent
${{R_\pm^{\hh^{d+1}}}}(\lam)$ on the model space $\hh^{d+1}$. This
will be done in $3$ steps. If $E,F$ are smooth bundles over
$\bar{X}$, we will say that a family of operator
$A(\lambda):\dot{C}^\infty(\bar{X},E)\to C^{-\infty}(\bar{X},F)$
depending meromorphically on a parameter $\lambda\in\cc$ is
\emph{finite meromorphic} if the polar part of $A(\lambda)$ at any
pole is a finite rank operator.

\begin{proposition}\label{extresolv}
Let $(X,g)$ be a spin asymptotically hyperbolic manifold and $D$
be a Dirac operator over $X$. Then the resolvent
$R_\pm(\lam)=(D\pm i\lam)^{-1}$ extends from $\{\Re(\lam)>0\}$ to
$\cc\setminus (-\nn/2)$ as a finite meromorphic family of
operators in
$\Psi_0^{-1,\lam+\ddemi,\lam+\ddemi}(\bar{X};\mc{E})$. Moreover
$R_\pm(\lambda)$ maps $\dot{C}^\infty(\bar{X},{^0\Sigma})$ to
$x^{\lam+\ddemi}C^{\infty}(\bar{X},{^0\Sigma})$ and for all
$\sigma\in\dot{C}^\infty(\bar{X},{^0\Sigma})$, we have
$[x^{-\ddemi-\lam}R_\pm(\lam)\sigma]|_{x=0}\in
C^{\infty}(\pl\bar{X},{^0\Sigma_\mp})$.
\end{proposition}
\begin{proof}
First, we construct an operator $Q^0_\pm(\lam)\in
\Psi^{-1,\infty,\infty}_0(\bar{X},\mc{E})$
supported near the interior diagonal such that  $(D\pm
i\lam)Q^0_\pm(\lam)=\mathrm{Id}+K^0_\pm(\lam)$ with $K^0_\pm(\lam)\in
\Psi^{-\infty,\infty,\infty}_0(\bar{X},\mc{E})$,
thus a smooth kernel on $\bar{X}\x_0\bar{X}$ and whose support
actually does not intersect the right and left boundary. Note that
this can be done thanks to the ellipticity of $D$ and it can be
chosen analytic in $\lam$, moreover notice also that
$Q^0_\pm(\lam)(D\pm
i\lam)-\mathrm{Id}\in\Psi^{-\infty,\infty,\infty}_0(\bar{X};\mc{E})$  by standards arguments of
pseudo-differential calculus. The error $K^0_\pm(\lam)$ is a
priori not compact on any weighted space
$x^{s}L^2(X,^0\Sigma)$ so this parametrix is not sufficient
for our purpose. To be compact on such a space, it would be enough
to have vanishing of the error on the front face $K^0_\pm(\lam)|_{\ff}=0$.

Next we need to solve away the term at the front face $\ff$, i.e.,
$K^0_\pm(\lam)|_{\ff}$. We can use the normal operator of $D$: the
normal operator $N_{y'}(D)$ of $D$ at $y'\in\pl\bar{X}$ is an operator
acting on the space $X_{y'}:=\{z\in {^0T}_{y'}\bar{X},
\frac{dx}{x}(z)>0\}$ obtained by freezing coefficients of $D$ at
$y'\in\pl\bar{X}$, when considered as polynomial in the $0$-vector
fields $x\pl_x,x\pl_y$. Here the spinor bundle over $X_{y'}$ is
trivial, i.e., it is given by $X_{y'}\x {^0\Sigma}_{y'}(\bar{X})$ where
${^0\Sigma}_{y'}(\bar{X})$ is the fiber of ${^0\Sigma}$  at the
boundary point ${y'}\in\pl\bar{X}$. The half space $X_{y'}$ equipped
with the metric $g$ frozen at ${y'}$ (the metric here is considered
as a symmetric tensor on the 0-cotangent space ${^0T^*\bar{X}}$)
is isometric to $\hh^{d+1}$ and the operator $N_{y'}(D)$ corresponds to
the Dirac operator on $\hh^{d+1}$ using this isometry. Moreover this
half space is also canonically identified to the interior of the
front face fibre $\ff_p$ with basis point $p=(y',y')\in\Delta_\pl$.
One has from \cite{MM} that the composition $(D\pm i\lam)G_{{\pm}}$
for $G_\pm\in\Psi^{-\infty,\alpha,\beta}_0(\bar{X};\mc{E})$ is in the calculus
$\Psi_0(\bar{X})$ and the restriction at the front face fiber
$\ff_p$ is given by
\[((D\pm i\lam)G_{\pm})|_{\ff_p}=N_{y'}(D\pm i\lam).G_\pm|_{\ff_p}\]
which is understood as the action of the differential operator
$N_{y'}(D\pm i\lam)$ on the conormal distribution $G_\pm|_{\ff_p}$ on
$\ff_p\simeq X_{y'}$. Thus to solve away the error term at $\ff$, it
suffices to find an operator $Q^1_\pm(\lam)$ in the calculus such that
\[N_{y'}(D\pm i\lam).Q^1_\pm(\lam)|_{\ff_p}=-K^0_\pm(\lam)|_{\ff_p}\]
for all ${y'}\in\pl\bar{X}$. This can be done smoothly in $y'$ by
taking
${Q^1_\pm(\lam)}|_{\ff_p}:=-R^{X_{y'}}_{\pm}(\lam)(K^0_\pm(\lam)|_{\ff_p})$
where $R^{X_{y'}}_{\pm}(\lam)$ is the analytically extended
resolvent of $N_{y'}(D\pm i\lam)\simeq (D_{\hh^{d+1}}\pm i\lam)$
on $X_{y'}\x {^0\Sigma}_{y'}(\bar{X})\simeq \hh^{d+1}\x
{^0\Sigma}(\hh^{d+1})$, and then defining ${Q^1_\pm(\lam)}$ to be
a distribution on $\bar{X}\x_0\bar{X}$ whose restriction to each
fiber $\ff_p$ is
$-R^{X_{y'}}_{\pm}(\lam)(K^0_\pm(\lam)|_{\ff_p})$. As we studied
above, the resolvent $R^{\hh^{d+1}}_{\pm}(\lam)$ is analytic in
$\lam$ and it maps $\dot{C}^\infty(\hh^{d+1},{^0\Sigma})$ to
$\rho^{\lam+\ddemi}C^{\infty}(\bbar{\hh^{d+1}},{^0\Sigma})$ if
$\rho$ is a boundary defining function of the compactification of
$\hh^{d+1}$, moreover the leading asymptotic term is of the form
$\rho^{\lam+\ddemi}\psi_\mp$ for some $\psi_\mp\in
C^\infty(\pl\bbar{\hh^{d+1}},{^0\Sigma_\mp})$. Thus, the
composition $R^{X_{y'}}_{\pm}(\lam)(K^0_\pm(\lam)|_{\ff_p})$ is a
conormal distribution in the class
$(\rho_{\lb}\rho_{\rb})^{\lam+\ddemi}C^{\infty}(\ff_p,{\rm
End}({^0\Sigma}_{y'}))$ and it is then possible to find
$Q^1_\pm(\lam)\in
\Psi_0^{-\infty,\lam+\ddemi,\lam+\ddemi}(\bar{X};\mc{E})$ with the
correct restriction at $\ff$. Let $P_\pm$ denotes the canonical
projection $P_\pm:{^0\Sigma}(\pl\bar{X})\to
{^0\Sigma_\pm}(\pl\bar{X})$. The restriction of a conormal kernel
in $\Psi_0^{-\infty,0,\beta}(\bar{X};\mc{E})$ at $\lb$ can be
considered as a section $C^{-\infty}(\pl\bar{X}\x \bar{X};\mc{E})$
conormal to all boundary faces. From the mapping property of
$R^{\hh^{d+1}}_{\pm}(\lam)$ just discussed, it is possible to
choose $Q^1_\pm(\lam)$ such that
$P_\pm[\rho_{\lb}^{-\lam-\ddemi}Q^1_\pm(\lam)]|_{\lb}=0$, which
will be important for the next step. Then we get $(D\pm
i\lam)(Q^0_\pm(\lam)+Q^1_\pm(\lam))=\mathrm{Id}+K^1_\pm(\lam)$
where $K^1_\pm(\lam)\in
\rho_{\ff}\Psi_0^{-\infty,\lam+\ddemi,\lam+\ddemi}
(\bar{X},\mc{E})$.

The final terms in the parametrix are those at the left boundary, solved away through the indicial equation for $z\in\cc$
\begin{align}\label{indicialeq}
(D\pm i\lam)&x^{\ddemi+z}(\psi_++\psi_-)\\
=&i(z\pm\lam)x^{\ddemi+z}\psi_+ +i(-z\pm\lam)x^{\ddemi+z}\psi_-+
O(x^{\ddemi+z+1}), \quad \forall \psi_\pm\in
C^\infty(\pl\bar{X},{^0\Sigma_\pm}).\notag
\end{align}
which is an easy consequence of \eqref{diracAH}. Lifting $D$ as
acting on the left variable on the space $\bar{X}\x_0\bar{X}$, it
satisfies the same type of indicial equation: if $G\in
\Psi_0^{-\infty,\alpha+\ddemi,\beta+\ddemi}(\bar{X};\mc{E})$ for
some $\alpha,\beta\in\cc$, then $(D\pm
i\lam)G\in\Psi_0^{-\infty,\alpha+\ddemi,\beta+\ddemi}(\bar{X};\mc{E})$
and the leading term at $\lb$ is
\begin{equation}\label{leadingterm}
[i(\alpha\pm\lam)\rho_{\lb}^{\alpha+\ddemi}P_+[(\rho_{\lb}^{-\alpha-\ddemi}G)|_{\lb}]
+i(-\alpha\pm\lam)\rho_{\lb}^{\alpha+\ddemi}P_-[(\rho_{\lb}^{-\alpha-\ddemi}G)|_{\lb}]
\end{equation}
where  the restriction at $\lb$ is considered as a section
$C^{-\infty}(\pl\bar{X}\x \bar{X};\mc{E})$ (conormal
to all boundary faces). Then since for $\alpha=\lam$, the term \eqref{leadingterm}
vanishes if
${{P_{\pm}(\rho_{\lb}^{-\lambda-\ddemi}G)|_{\lb}=0}}$, one clearly has
$K^1_\pm(\lam)\in\rho_{\ff}\Psi_0^{-\infty,\lam+\ddemi+1,\lam+\ddemi}(\bar{X},\mc{E})$ thanks to the choice of
$Q^1_\pm(\lam)$ and now, since $\alpha=\lam+j$ for $j\in\nn$ is
not solution of the indicial equation above when ${\lam\notin-\nn/2}$, it is possible
by induction and using Borel Lemma to construct a term
$Q^2_\pm(\lam)\in \rho_{\ff}\Psi_0^{-\infty,\lam+\ddemi+1,\lam+\ddemi}
(\bar{X},\mc{E})$, holomorphic in
${\cc\setminus(-\nn/2)}$  such that $(D\pm
i\lam)Q^2_\pm(\lam)=-K^1_\pm(\lam)+K^2_\pm(\lam)$ for some
operator  $K^2_\pm(\lam)\in
\rho_{\ff}\Psi_0^{-\infty,\infty,\lam+\ddemi}
(\bar{X},\mc{E})$.

By \cite{MazzeoCPDE}, the error term $K^2_\pm(\lam)$ is now
compact on $\rho^zL^2$ for all $z\in[0,\infty)$ such that
$\Re(\lam)+z>0$. Fix $\lam_0$ such that $\Re(\lam_0)>0$ where
$R_\pm(\lam_0)$ is bounded on $L^2(X)$. Then we use a standard
argument: we can add a finite rank term $Q^3_\pm(\lam)=Q^3_\pm
(\lam_0)\in \Psi^{-\infty,\infty,\infty}(\bar{X};\mc{E})$ to
$Q^2_\pm(\lam)$ in case $(\mathrm{Id}+K^2_\pm(\lam_0))$ has kernel
in $\rho^zL^2$, so that that $\mathrm{Id}+K^3_\pm(\lam_0)$ is
invertible if $K^3_\pm(\lam):=(D\pm
i\lam)(\sum_{i=0}^3Q^i_\pm(\lam))-{\rm Id}$  (see the proof of Th.
7.1 in \cite{MM} for more details). The operator
$Q_\pm(\lam):=\sum_{i=0}^3Q^i_\pm(\lam)$ is bounded from
$\rho^{z}L^2$ to $\rho^{-z}L^2$ if $\Re(\lam)+z>0$ by
\cite{MazzeoCPDE}, and so Fredholm theorem proves that
$R_\pm(\lam)=Q_\pm(\lam)(\mathrm{Id}+K^3_\pm(\lam))^{-1}$ on the
weighted space $\rho^{z}L^2$ for $z\in[0,\infty)$ such that
$z+\Re(\lam)>0$. Finally, writing
$(\mathrm{Id}+K^3_\pm(\lam))^{-1}=\mathrm{Id}+T_\pm(\lam)$, we see
that
$T_\pm(\lam)=-K^3_\pm(\lam)+K^3_\pm(\lam)(\mathrm{Id}+K^3_\pm(\lam))^{-1}K^3_\pm(\lam)$
and the same arguments as \cite{MM} show that $R_\pm(\lam)\in
{\Psi_0^{-{1},\lam+\ddemi,\lam+\ddemi}}(\bar{X};\mc{E})$.  The
mapping property of $R_\pm(\lam)$ acting on
$\dot{C}^\infty(\bar{X},{^0\Sigma})$ follows from Mazzeo
\cite{MazzeoCPDE}, and the fact that for all $\sigma\in
\dot{C}^\infty(\bar{X},{^0\Sigma})$ we have
$R_\pm(\lam)\sigma=x^{\lam+\ddemi}\psi_\mp+O(x^{\lam+\ddemi+1})$
for some ${\psi_\mp}\in C^{\infty}(\pl\bar{X},{^0\Sigma_\mp})$ is a
straightforward consequence of the indicial equation
\eqref{indicialeq}.
\end{proof}

Let us now discuss the nature of the spectrum of $D$. We start by an application of Green formula, usually
called \emph{boundary pairing property} (compare to \cite[Prop. 3.2]{GRZ})
\begin{lemma}\label{boundarypairing}
Let $\Re(\lam)=0$ and $\sigma_i=x^{\ddemi-\lam}\sigma_i^-+x^{\ddemi+\lam}\sigma_i^+$ for $i=1,2$ and
$\sigma_i^{\pm}\in C^{\infty}(\bar{X},{^0\Sigma})$,
then if $(D^2+\lam^2)\sigma_i=r_i\in \dot{C}^{\infty}(\bar{X},{^0\Sigma})$, one has
\[\int_X (\cjg\sigma_1,r_2\cjd-\cjg r_1,\sigma_2\cjd) {\rm dv}_g=2\lam \int_{\pl\bar{X}}
\cjg\sigma_1^-|_{\pl\bar{X}},{\sigma_2^-|_{\pl\bar{X}}}\cjd-\cjg
\sigma_1^+|_{\pl\bar{X}},\sigma_2^+|_{\pl\bar{X}}\cjd {\rm
dv}_{h_0}\] where {$\cjg \ , \, \cjd$} denotes the scalar product
with respect to $g$ on $X$ and to $h_0$ on $\pl\bar{X}$.
\end{lemma}
\begin{proof}
The proof of this Lemma is straightforward by using integration by parts in $\{x\geq \eps\}$ and letting $\eps\to 0$.
\end{proof}

As a corollary of the resolvent extension and this Lemma, we obtain the following
\begin{corollary}\label{spectrum}
On a spin asymptotically hyperbolic manifold, the extended
resolvent $R_\pm(\lam)=(D\pm i\lam)^{-1}$ is holomorphic on the
imaginary line $i\rr$, consequently the spectrum of $D$ is $\rr$
and absolutely continuous.
\end{corollary}
\begin{proof} In view of the meromorphy of $R_\pm(\lam)$, it clearly
suffices from Stone's formula to prove that $R_\pm(\lam)$ has no
pole on the imaginary line. Assume $\lam_0$ is such a pole with
order $p$, then the most singular coefficient of the Laurent
expansion is a finite rank operator whose range is made of
generalized eigenspinors $\sigma$ solving $(D\pm i\lam_0)\sigma=0$
and $\sigma\in x^{\lam_0+\ddemi}C^\infty(\bar{X},{^0\Sigma})$. In
particular it satisfies ${(D^2+\lam_0^2)\sigma=0}$ and by applying
Lemma \ref{boundarypairing} with $\sigma_1=\sigma_2=\sigma$ we see
that $(x^{-\lam_0-\ddemi}\sigma)|_{x=0}=0$ and so $\sigma\in
x^{\lam_0+\ddemi+1}C^{\infty}(\bar{X},{^0\Sigma})$. Now from the
indicial equation \eqref{indicialeq}, this implies $\sigma\in
\dot{C}^{\infty}(\bar{X},{^0\Sigma})$ if $\lam_0\not=0$. Then
Mazzeo's unique continuation theorem \cite{MaAJM} says that for a
class of operator including $D^2$, there is no eigenfunction
vanishing to infinite order at the boundary except $\sigma\equiv 0$,
we deduce that $\sigma=0$ and thus by induction this shows that
the polar part of Laurent expansion of $R(\lam)$ at $\lam_0$ is
$0$. Now there remains the case $\lam_0=0$. First from
self-adjointness of $D$, we easily get
\begin{equation}\label{selfadj}
\lam ||\sigma||_{L^2}\leq ||(D\pm i\lam)\sigma||_{L^2}
\end{equation}
for all $\lam>0$ and $\sigma$ in the $L^2$-Sobolev space
$H^1(X,{^0\Sigma})$ of order $1$, and this implies that
$R_{\pm}(\lam)$ has a pole of order at most $1$ at $\lam=0$, i.e.
one has $R_\pm(\lam)=A_\pm\lam^{-1}+B_\pm(\lam)$ for some
$B_\pm(\lam)$ holomorphic (these can be considered as operators
from $\dot{C}^\infty(\bar{X},{^0\Sigma})$ to its dual). By
\eqref{selfadj}, we also see by taking $\lam\to 0$ that
$||A_\pm\sigma||_{L^2}\leq ||\sigma||_{L^2}$ for all $\sigma\in
\dot{C}^\infty(\bar{X},{^0\Sigma})$ and so $A_\pm$ is bounded on
$L^2$ and also in $\ker(D)$ by the functional equation $(D\pm
i\lam)R_\pm(\lam)={\rm Id}$. Now in view of the structure of the
kernel of $R_\pm(\lam)$, it is not hard to check (e.g. see
\cite{GuiMRL}) that the elements in the range of $A_\pm$ are
harmonic spinors of the form $\sigma\in
x^{\ddemi}C^{\infty}(\bar{X},{^0\Sigma})$, which can only be $L^2$
if the leading asymptotic $(x^{-\ddemi}\sigma)|_{x=0}=0$, i.e. if
$\sigma\in x^{\ddemi+1}C^{\infty}(\bar{X},{^0\Sigma})$. Using
again the indicial equation, \eqref{indicialeq}, we deduce that
$\sigma\in \dot{C}^\infty(\bar{X},{^0\Sigma})$ and thus $\sigma=0$
by Mazzeo's unique continuation theorem, so $A_\pm=0$.
\end{proof}

Another corollary of Proposition \ref{extresolv} is
\begin{corollary}\label{resolvforD}
The resolvent $R(\lam):=(D^2+\lam^2)^{-1}$ extends meromorphically
to $\lam\in \cc\setminus (-\nn/2)$ with poles of finite
multiplicity, except at $\lam=0$ where it has a simple infinite
rank pole with residue $(2i)^{-1}(R_-(0)-R_+(0))$. Moreover
$R(\lam)$ is an operator in
$\Psi_0^{-2,\lam+\ddemi,\lam+\ddemi}(\bar{X},\mc{E})$ and for any $\sigma\in
\dot{C}^{\infty}(\bar{X},{^0\Sigma})$, one has $R(\lam)\sigma \in
x^{\ddemi+\lam}C^{\infty}(\bar{X},{^0\Sigma})$.
\end{corollary}
\begin{proof}: The extension and the structure of $R(\lam)$ are a consequence of Proposition \ref{extresolv}
since $R(\lam)=(2i\lam )^{-1}(R_-(\lam)-R_+(\lam))$. As for the
mapping property, this is a consequence of mapping properties of
operators in $\Psi_0^{*,*,*}(\bar{X})$ in \cite{MazzeoCPDE}. The
question of the simple pole at $\lam=0$ is also clear since
$R_\pm(\lam)$ are holomorphic at $\lam=0$. It remains to show that
the residue $\Pi_0$ is infinite rank. One way to prove it is to
consider the asymptotic of $\Pi_0\phi$ if $\phi\in
C_0^\infty(\bar{X};{^0\Sigma})$. First, both $R_-(0)$ and $R_+(0)$
have infinite rank since $DR_\pm(0)={\rm Id}$ on
$C_0^\infty(\bar{X};{^0\Sigma})$, but moreover if $\phi$ is smooth
compactly supported, $R_\pm(0)\phi$ has an asymptotic of the form
$x^\ddemi\psi_\pm$ at the boundary where $\psi_\pm\in
C^\infty(\pl\bar{X};{^0\Sigma_\mp})$ according to Proposition
\ref{extresolv}. If $\psi_{{\pm}}=0$, then
$R_\pm(0)\phi=O(x^{\ddemi+1})$, and by the indicial equation
\eqref{indicialeq} it must vanish to infinite order at
$\pl\bar{X}$, which by Mazzeo's unique continuation theorem
implies that $R_\pm(0)\phi=0$, a contradiction. This then shows
that the range of $R_+(0)$ on $C_0^\infty(\bar{X};{^0\Sigma})$
does not intersect the range of $R_-(0)$ acting on the same space,
concluding the proof.
\end{proof}
\begin{remark*} By self-adjointness of $D^2$,
one deduces easily that $R(\lam)^*=R(\bar{\lam})$, or in terms of kernels
\[R(\lam;m,m')^*=R(\bar{\lam};m',m) \quad \forall m, m'\in X, \quad m\not=m',\]
here $A^*\in\Sigma_{m'}\boxtimes\Sigma_{m}$ means the adjoint of
$A\in\Sigma_m\boxtimes \Sigma_{m'}$ if we identify the dual
$\Sigma_m^*$ with $\Sigma_m$ via the Hermitian product induced
by the metric $g$.
\end{remark*}

\subsection{Another parametrix construction when the curvature is constant near $\infty$}\label{caseconvexccompact}
When $(X,g)$ is asymptotically hyperbolic with constant curvature
outside a compact set (which is the case of a convex co-compact
quotient $X_\Gamma=\Gamma\backslash \hh^{d+1}$), one may use a
simplified construction similar to that of Guillop\'e-Zworski
\cite{GZ} for the Laplacian on functions.

Indeed, there exists a covering of a neighbourhood of $\pl\bar{X}$ by open sets $(U_j)$ with isometries
\[\begin{gathered}
\iota_j: (U_j,g) \to (B,g_{\hh^{d+1}}),\\
\textrm{ where }B :=\{(x_0,y_0)\in (0,\infty)\x\rr^{d},
x_0^2+|y_0|^2<1\} \quad \textrm{ and }\quad
g_{\hh^{d+1}}=\frac{dx_0^2+|dy_0|^2}{x_0^2}.
\end{gathered}\]
We denote by $\bar{B}:=\{(x_0,y_0)\in [0,\infty)\x\rr^{d},x_0^2+|y_0|^2<1 \}$ a half-ball in $\hh^{d+1}$
and $\pl\bar{B}:=\bar{B}\cap\{x_0=0\}$.
We shall also use the notation $\iota_j$ for the restriction 
$\iota_j|_{U_j\cap\pl\bar{X}}$.

Note that the function $x_0$ in $B$ is not pulled-back to a boundary defining
function putting the metric $g$ under the form $g=(dx^2+h_x)/x^2$, but we have $(x/\iota_j^*x_0)|_{\pl\bar{X}}=\iota_j^*\eta_j$
for some functions $\eta_j\in C^\infty(\pl\bar{B})$.
Through $\iota_j$, the spinor bundle on $\hh^{d+1}$ pulls-back to the spinor bundle ${^0\Sigma}(\bar{X})|_{U_j}$ but the splitting
induced by ${\rm cl}(x_0\pl_{x_0})$ does not correspond to the splitting ${^0\Sigma}_+\oplus{^0\Sigma}_-$, except at
the boundary $x_0=0$, since the eigenspaces of ${\rm cl}(\iota_j^*(x_0\pl_{x_0}))$  are the eigenspaces of ${\rm cl}(x\pl_x)$
when restricted to the boundary.

The Dirac operator $D_{\hh^{d+1}}$ pulls back to $D$ in $U_j$,  consequently one may choose a partition
of unity $\iota^*_j\chi^1_j$  near $\pl\bar{X}$ associated to $(U_j)_j$, that is $\sum_j \iota^*_j\chi^1_j=\chi$
where $\chi \in C^\infty(\bar{X})$ is equal to $1$ near
$\pl\bar{X}$. Take now $\chi^2_j\in C_0^{\infty}(\bar{B})$ some functions which are equal to $1$ on the support of $\chi^1_j$.
Let us define $\phi^i_j$ on $\pl\bar{B}$ by $\chi^i_j(0,y_0)=\phi^i_j(y_0)$ so that $\sum_{j}\iota_j^*\phi^1_j=1$ on $\pl\bar{X}$
and $\phi^2_j\phi^1_j=\phi^1_j$.

The first parametrix we can use for
$(D^2+\lam^2)$ is
\begin{equation}\label{defE0}
R_0(\lam)=\sum_{j}\iota_j^*\chi^2_j R_{\hh^{d+1}}(\lam)\chi^1_j{\iota_j}_* + Q_0(\lam)
\end{equation}
where $Q_0(\lam)\in \Psi_0^{-2;\infty,\infty}(\bar{X};\mc{E})$ is holomorphic, compactly supported and solves
$(D^2+\lam^2)Q_0(\lam)=1-\chi+K_0(\lam)$ for some $K_0(\lam)\in C^\infty_0(X\x X;\mc{E})$. Here
$\iota_j^*$ denotes the pull back on sections of the spinor bundle and ${\iota_j}_*:=(\iota_j^{-1})^*$.
We obtain
\[(D^2+\lam^2)R_0(\lam)=\mathrm{Id} + \sum_{j}\iota_j^*[D^2_{\hh^{d+1}},\chi^2_j]R_{\hh^{d+1}}(\lam)\chi^1_j{\iota_j}_*+
K_0(\lam).\] The last term  $K_0(\lam)$ is clearly
compact on all weighted spaces $x^NL^2(X,{^0\Sigma})$ while the
first one is not. Since on $\hh^{d+1}$ one has
\[D_{\hh^{d+1}}^2=x_0^\ddemi\left( -(x_0\pl_{x_0})^2{\rm Id}+x_0^2D_{\rr^{d}}+ix_0AD_{\rr^{d}}\right)x_0^{-\ddemi} , \quad A:=\left[
\begin{matrix}
1 & 0\\
0 & -1
\end{matrix}\right]=-i{\rm cl}(x_0\pl_{x_0}),\]
in the splitting induced by ${\rm cl}(x_0\pl_{x_0})$. The operator
$[D_{\hh^{d+1}}^2,\chi^2_j]$ can be written as follows,
\begin{align*}
[D_{\hh^{d+1}}^2,\chi^2_j]=&d x_0\pl_{x_0}
\chi^2_j(x_0,y_0)){\rm Id}\\
& -x_0^2 (\pl_{x_0}^2 \chi^2_j(x_0,y_0)){\rm Id}
+x_0^2[D_{\rr^{d}},\chi^{2}_j(x_0,y_0)]+ix_0[AD_{\rr^{d}},\chi^2_j(x_0,y_0)].
\end{align*}
Using the fact that $(\nabla
\chi^2_j)\chi^1_j=0$ and the expression of $R_{\hh^{d+1}}(\lam)$, we
deduce that
\[ [D_{\hh^{d+1}}^2,\chi^2_j]R_{\hh^{d+1}}(\lam)\chi_j^1 \in x_0^{\lam+\ddemi+1}{x_0'}^{\lam+\ddemi}C^{\infty}(
\bar{B}\x\bar{B},\mc{E})\]
where $(x_0,y_0,x_0',y_0')$ are the natural coordinates on $B\x B$.  This error term can be solved away using the indicial equation explained
above for the general AH case and one can thus construct, for all $N\in \nn$, an operator $R_N(\lam)\in x^{\lam+\ddemi+1}{x'}^{\lam}C^{\infty}(\bar{X}\x\bar{X}, \mc{E})$ such that
\[(D^2+\lam^2)(R_0(\lam)+R_N(\lam))=\mathrm{Id}+K_N(\lam), \quad K_N(\lam)\in x^{\lam+\ddemi+N}{x'}^{\lam+\ddemi}C^{\infty}(\bar{X}\x\bar{X}, \mc{E})\]
and $K_N(\lam)$ is compact on $x^{N'}L^2(X,{^0\Sigma})$ if
$0<N'<N$ and $\Re(\lam)>-N+N'$ and $\Re(\lam)>-N'$. All these
terms are holomorphic in $\lam$ except possibly at $-\nn/2$ where first order
poles come from the indicial equation and at $\lam=0$ where $R_{\hh^{d+1}}(\lam)$ has an infinite rank pole.
As above for the general
case, we can take an asymptotic series using Borel Lemma, which
gives an operator $R_\infty(\lam)\in
x^{\lam+\ddemi+1}{x'}^{\lam+\ddemi}C^{\infty}(\bar{X}\x\bar{X};
\mc{E})$, holomorphic in $\lam\notin -\nn_0/2$ so that
$(D^2+\lam^2)(R_0(\lam)+R_{\infty}(\lam))={\rm
Id}+K_{\infty}(\lam)$ for some $K_\infty(\lam)\in
x^{\infty}{x'}^{\lam+\ddemi}C^{\infty}(\bar{X}\x\bar{X})$. And
again , as in the proof of Proposition \ref{extresolv}, up to the
addition of a residual finite rank term for $R_\infty(\lam)$, we
can assume that there is $\lam_0$ with $\Re(\lam_0)>0$ such that
 ${\rm Id}+K_{\infty}(\lam_0)$ is invertible on $x^NL^2(X)$ for all $N>0$.
The extended resolvent of $D^2+\lam^2$ is thus given by
\[R(\lam)=(R_0(\lam)+R_\infty(\lam))({\rm Id}+K_{\infty}(\lam))^{-1},\]
it is finite meromorphic in $\cc\setminus(-\nn_0/2)$. Moreover
standard arguments show that
$(\mathrm{Id}+K_{\infty}(\lam))^{-1}=\mathrm{Id}+S_\infty(\lam)$
for some $S_\infty(\lam)\in
x^{\infty}{x'}^{\lam+\ddemi}C^{\infty}(\bar{X}\x\bar{X},\mc{E})$
and so
\[
R(\lam)-R_0(\lam)({\rm Id}+S_\infty(\lam))\in x^{\lam+\ddemi+1}{x'}^{\lam+\ddemi}C^{\infty}(\bar{X}\x\bar{X}, \mc{E}).
\]
Using the composition results of Mazzeo
\cite[Th.\,3.15]{MazzeoCPDE}, we get
$R_0(\lam)S_\infty(\lam)\in
(xx')^{\lam+\ddemi}C^{\infty}(\bar{X}\x\bar{X};\mc{E})$ so
\begin{equation}\label{RmoinsR0}
R(\lam)-R_0(\lam)\in x^{\lam+\ddemi}{x'}^{\lam+\ddemi}C^{\infty}(\bar{X}\x\bar{X}, \mc{E}).
\end{equation}
Similarly, using the remark following Corollary \ref{resolvforD},
we deduce that the kernel
\begin{equation}\label{Rtransp1}
(m,m')\to R(\lam;m,m')-R_0(\bar{\lam};m',m)^* \in
{x}^{\lam+\ddemi}{x'}^{\lam+\ddemi}C^{\infty}(\bar{X}\x\bar{X},
\mc{E})
\end{equation}
and $R_0(\bar{\lam};m',m)^*$  is given, by symmetry of
$R_{\hh^{d+1}}(\lam)$, by
\begin{equation}\label{Rtransp2}
R_0(\bar{\lam};m',m)^*= \sum_{j}(\iota_j^*\chi^1_j R_{\hh^{d+1}}(\lam)
\chi^2_j{\iota_j}_*)(m,m') + Q_0(\bar{\lam})^*.
\end{equation}
The expressions \eqref{RmoinsR0}, \eqref{Rtransp1} and \eqref{Rtransp2}
will be very useful in what follows for obtaining an explicit formula of the scattering operator
modulo a smoothing term.

\section{Scattering and Eisenstein series}

\subsection{Definitions and properties}
Similarly to the Laplacian on functions, we can define Eisenstein
series and scattering operator for Dirac operator. The Eisenstein
series $E(\lam)$ is an operator mapping
$C^{\infty}(\pl\bar{X},^0\Sigma)\to C^\infty(X,{^0\Sigma})$ and
for all $\psi$, $E(\lam)\psi$ is a non $L^2$-solution of
$(D^2+\lam^2)\sigma=0$; more precisely it is defined using the
following
\begin{lemma}\label{generalizedeig}
Let $\psi\in C^\infty(\pl\bar{X},^0\Sigma)$, and
$\lam\in\cc\setminus(-\nn/2)$ not a pole of $R(\lam)$, then there
exists $\sigma\in C^{\infty}(X,{^0\Sigma})$ solution of
$(D^2+\lam^2)\sigma=0$, unique when $\Re(\lam)\geq 0$, and such
that there exist $\sigma^\pm \in C^{\infty}(\bar{X},{^0\Sigma})$
with $\sigma^-|_{\pl\bar{X}}=\psi$ and
$\sigma=x^{\ddemi-\lam}\sigma^-+x^{\ddemi+\lam}\sigma^+$. Moreover
$\sigma^{\pm}$ depend meromorphically on
$\lam\in\cc\setminus(-\nn/2)$.
\end{lemma}
\begin{proof}
This is essentially the same construction as for the Laplacian on
functions in \cite{GRZ}: using the indicial equation
\eqref{indicialeq} and Borel lemma, it is possible to construct a
spinor $\sigma_\infty\in
x^{-\lam+\ddemi}C^{\infty}(\bar{X},{^0\Sigma})$, holomorphic in
$\cc\setminus(\zz/2)$ such that
$(D^2+\lam^2)\sigma_\infty=O(x^\infty)$ and
$(x^{\lam-\ddemi}\sigma_\infty)|_{\pl\bar{X}}=\psi$. Note that this
spinor is meromorphic in $\lam\in\cc$ with only simple poles at
$\nn/2$ coming from the roots of the indicial equation. Then we
can set
\begin{equation*}\label{e:def-Eisen}
\sigma:=\sigma_{\infty}-R(\lam)(D^2+\lam^2)\sigma_\infty.\end{equation*}
If $\lam$ is not a pole of $R(\lam)$, this solves the problem and
defines $\sigma^\pm$ by using the mapping property of $R(\lam)$
stated in Corollary \ref{resolvforD}. The meromorphy of
$\sigma^\pm$ is also a consequence of the construction and of the
meromorphy of $R(\lam)$. The uniqueness of the solution is due to
the fact that for two solutions $\sigma_1,\sigma_2$ of the
problem, the indicial equation implies that $\sigma_1-\sigma_2\in
x^{\lam+\ddemi}C^\infty(\bar{X},{^0\Sigma})$, and then for
$\Re(\lam)>0$ this would be $L^2$ and $\lam$ would be a pole of
the resolvent in the physical plane.  For $\Re(\lam)=0$, this can
be proved using an application of Green formula like in
\cite{GRZ}: if $\til{\sigma}_1,\til{\sigma}_2$ are two solutions
of the problem, then the difference
$\til{\sigma}_1-\til{\sigma}_2$ is in
$x^{\ddemi+\lam}C^{\infty}(\bar{X},{^0\Sigma})$ by the indicial
equation and it is also in the kernel of $D^2+\lam^2$, so we may
apply the Lemma \ref{boundarypairing} with
$\sigma_1=\til{\sigma}_1-\til{\sigma}_2$ and
$\sigma_2:=R(-\lam)\varphi$ where $\varphi\in
\dot{C}^{\infty}(\bar{X},{^0\Sigma})$ is chosen arbitrarily. This
clearly implies that $\int_X
\cjg\til{\sigma}_1-\til{\sigma}_2,\varphi\cjd{\rm dv}_g=0$ and
thus $\til{\sigma}_1=\til{\sigma}_2$.
\end{proof}

\begin{remark*} By uniqueness of the solution, $\sigma$ and
$\sigma^\pm|_{\pl\bar{X}}$ depend linearly on $\psi$.
\end{remark*}

\begin{definition}
The Eisenstein series is the operator $E(\lam):C^{\infty}(\pl\bar{X},{^0\Sigma})\to
C^{\infty}(X,{^0\Sigma})$ defined by $E(\lam)\psi:=\sigma$ where
$\sigma$ is the smooth spinor in Lemma \ref{generalizedeig}.
\end{definition}

\begin{definition}
The scattering operator $S(\lam):C^{\infty}(\pl\bar{X},{^0\Sigma})\to
C^{\infty}(\pl\bar{X},{^0\Sigma})$ is defined by
$S(\lam)\psi:=\sigma^+|_{\pl\bar{X}}$ where $\sigma^+$ is the smooth spinor in
Lemma \ref{generalizedeig}.
\end{definition}

It is rather easy to prove that the scattering operator is
off-diagonal with respect to the splitting
${^0\Sigma}(\pl\bar{X})={^0\Sigma}_+(\pl\bar{X})\oplus{^0\Sigma}_-(\pl\bar{X})$.
To that end, we give an alternative construction of the
Eisenstein series $E(\lam)\psi$ when $\psi\in {^0\Sigma}_+$ or
$\psi\in{^0\Sigma}_-$.
Let us first define a useful meromorphic function on $\cc$
\begin{equation}\label{defclam}
C(\lam):=2^{-2\lam}\frac{\Gamma(\demi-\lam)}{\Gamma(\demi +\lam)}.
\end{equation}
which satisfies $C(\lam)C(-\lam)=1$.

\begin{lemma}\label{alternative}
Let $\psi\in C^\infty(\pl\bar{X},^0{\Sigma}_\pm)$, and
$\lam\in\cc\setminus(-\nn/2)$ be not a pole of $R_{{\pm}}(\lam)$, then there exists a unique
$\sigma\in C^{\infty}(X,{^0\Sigma})$ solution of $(D\pm
i\lam)\sigma=0$ and such that there exists $\sigma^\pm \in
C^{\infty}(\bar{X},{^0\Sigma})$ with $\sigma^-|_{\pl\bar{X}}=\psi$
and $\sigma=x^{\ddemi-\lam}\sigma^-+x^{\ddemi+\lam}\sigma^+$.
Moreover one has $\sigma^+|_{\pl\bar{X}}\in
C^{\infty}(\pl\bar{X},{\Sigma_{\mp}})$ and $\sigma^\pm$ are
meromorphic in $\lam\in\cc\setminus(-\nn/2)$. If in addition the
metric $g$ is even, then $\sigma^\pm/C(\lam)$ are holomorphic in
$\{\Re(\lam)\geq 0\}$ where $C(\lam)$ is the function in
\eqref{defclam}.
\end{lemma}
\begin{proof}
Recall the indicial equation for $(D\pm i\lam)$: let $j\in\nn$ and
$\psi_\pm\in C^{\infty}(\pl\bar{X},{^0\Sigma}_\pm)$ then there exist
some smooth spinor $F_{\lam,j}$ near $\pl\bar{X}$ such that
\begin{equation}\label{indicialequ}
x^{\lam-\ddemi}(D\pm i\lam)x^{-\lam+\ddemi+j}(\psi_++\psi_-)=ix^{j}\Big((j-\lam\pm\lam)\psi_++(\lam-j\pm\lam)\psi_-\Big)+
x^{j+1}F_{\lam,j}. 
\end{equation} 
Using this indicial equation inductively, we can construct for all $\psi\in C^{\infty}(\pl\bar{X},{^0\Sigma}_\pm)$ a formal
Taylor series, and thus a true spinor $\sigma_{\infty,\pm}\in
x^{-\lam-\ddemi}C^{\infty}(\bar{X},{^0\Sigma})$ by Borel lemma,
such that $(D\pm i\lam)\sigma_{\infty,\pm}=O(x^\infty)$ and
$(x^{\lam-\ddemi}\sigma_{\infty,\pm})|_{\pl\bar{X}}=\psi$. This can be done holomorphically in $\lam$ 
as long as $\lam$ is not a root of the indicial equation \eqref{indicialequ}. 
The $\lam$ so that the indicial numbers $j-\lam\pm \lam$ and
$\lam-j\pm\lam$ in \eqref{indicialequ} vanish are at $\nn/2$ and
they vanish only on the $\Sigma_\mp$ part of the bundle, therefore since $C(\lam)$ has first order poles at $1/2+\nn_0$, we see that
$\sigma_\infty/C(\lam)$ can be chosen holomorphically in
$\lam\in\cc\setminus \nn$ and that it has at most poles of order
$1$ at each $\lam=k$ with $k\in\nn$; since it does not involve new
arguments we refer the reader who needs more details to the paper
of Graham-Zworski \cite{GRZ} where it was studied in the case of
the Laplacian on functions. Now consider the case of a metric $g$ even. 
Since $\mc{A}_+,\mc{A}_-$ defined in\eqref{defA} are preserved by $x^{-\ddemi}Dx^{\ddemi}$ if $g$ is even, then 
clearly $x^{\lam-\ddemi}(D\pm i\lam)x^{-\lam+\ddemi}=\lam(-{\rm cl}(x\pl_x)\pm i)+x^{-\ddemi}Dx^\ddemi$ also preserves
both $\mc{A}_+,\mc{A}_-$. In particular, if $\psi_-=0$ in \eqref{indicialequ}, then $x^{2j+1}F_{\lam,2j}\in \mc{A}_+$
and $x^{2j+2}F_{\lam,2j+1}\in\mc{A}_-$, while the converse is true if $\psi_+=0$. 
This implies that  
the spinor $\sigma_{\infty,\pm}$ can be taken so that $x^{\lam-\ddemi}\sigma_{\infty,\pm}\in\mc{A}_\pm$
and the $\lam$ which are actually solution of the indicial equation \eqref{indicialequ} for $D\pm i\lam$ are
only at $1/2+\nn_0$.
The spinor $\sigma_\infty/C(\lam)$ can be taken holomorphic also at $\lam\in\nn$. It remains
to set
\begin{equation}\label{caspm}
\sigma:=\sigma_{\infty,\pm}- R_\pm(\lam)(D\pm i\lam)\sigma_{\infty,\pm}
\end{equation}
which solves our problem, using the mapping property of
$R_\pm(\lam)$ stated in Proposition \ref{extresolv}.
\end{proof}

By uniqueness, the solution in Lemma \ref{alternative} is clearly
the same than the one of Lemma \ref{generalizedeig} when the
initial data $\psi$ is either in $\sigma_+|_{\pl{\bar{X}}}$ or
$\sigma_-|_{\pl{\bar{X}}}$, which implies
\begin{corollary}
The scattering operator $S(\lam)$ maps
$C^{\infty}(\pl\bar{X},{^0\Sigma}_\pm)$ to
$C^{\infty}(\pl\bar{X},\Sigma_\mp)$.
\end{corollary}

Let us define the natural projection and inclusion
\[P_\pm:C^{\infty}(\pl\bar{X};{^0\Sigma})\to C^{\infty}(\pl\bar{X};{^0\Sigma}_\pm); \quad
I_\pm: C^\infty(\pl\bar{X};{^0\Sigma}_\pm)\to
C^{\infty}(\pl\bar{X};{^0\Sigma})\] and also the maps corresponding
to the two off-diagonal components of $S(\lam)$
\[\ \ \ \ S_\pm(\lam):=P_\mp S(\lam)I_\pm : C^\infty(\pl\bar{X};{^0\Sigma}_\pm)\to
C^{\infty}(\pl\bar{X};{^0\Sigma}_\mp)\]
\[E_\pm (\lam):=E(\lam)I_\pm P_{\pm} :C^\infty(\pl\bar{X};{^0\Sigma})\to
C^{\infty}(X;{^0\Sigma}).\]

\subsection{Some relations between resolvent, scattering operator and Eisenstein series}

Like for the Laplacian on functions, the Schwartz kernels of
$R(\lam)$, $E(\lam)$ and $S(\lam)$ are related by the following
\begin{proposition}\label{kernels}
Let $\lam\in\cc$ be such that $\lam\notin -\nn/2$ and $\lam$ not a
pole of $R(\lam)$, then the Schwartz kernel $E(\lam;z,y')$ and
$E_{\pm}(\lam;z,y')$ in $C^{-\infty}(\bar{X}\x \pl\bar{X};\mc{E})$
of respectively $E(\lam)$ and $E_\pm(\lam)$ can be expressed by
\begin{equation}\label{kernelEla}
\begin{gathered}
E(\lam;z,y')=2\lam[{x'}^{-\ddemi-\lam}R(\lam;z,x',y')]|_{x'=0},\\
E_\pm(\lam;z,y')=[{x'}^{-\ddemi-\lam}R_\pm(\lam;z,x',y')]|_{x'=0}{\rm cl}(\nu)
\end{gathered}
\end{equation}
where we use the product decomposition $(x',y')\in
[0,\eps)\x\pl\bar{X}$ near the boundary in the right variable of
$\bar{X}\x\bar{X}$. If in addition $\Re(\lam)<-\ddemi$,
the Schwartz kernel $S(\lam;y,y')$ of $S(\lam)$ is in
$C^0(\pl\bar{X}\x\pl\bar{X};\mc{E})$ and can be expressed by
\begin{equation}\label{kernelSla}
S(\lam;y,y')=[x^{-\ddemi-\lam}E(\lam;x,y,y')]_{x=0}, \quad
S_\pm(\lam;y,y')=[x^{-\ddemi-\lam}E_\pm(\lam;x,y,y')]|_{x=0}.
\end{equation}
\end{proposition}
\begin{proof}
Let $z\in X$ and $\sigma_\infty$ as in Lemma \ref{generalizedeig},
then
$E(\lam)\psi=\sigma_{\infty}-R(\lam)(D^2+\lam^2)\sigma_{\infty}$.
The first statement of the Lemma is simply obtained by integration
by part in $(x',y')$ of
\[\int_{x'\geq \eps}\cjg R(\lam;z,x',y'),(D^2+\lam^2)\sigma_{\infty}(x',y')\cjd{\rm dv}_g(x',y')\]
and letting $\epsilon \to 0$, this  gives the term
$\sigma_\infty(z)$ plus a term
\[\int_{\pl\bar{X}} \big{[}{x'}^{-d}\Big(\cjg R(\lam;z,x',y'),\nabla_{x'\pl_{x'}}\sigma_{\infty}(x',y')\cjd
-\cjg
\nabla_{x'\pl_x'}R(\lam;z,x',y'),\sigma_{\infty}(x',y')\cjd\Big)\big{]}_{x'=0}\,
{\rm dv}_{h_0}.\]
But from the analysis of the
resolvent $R(\lam)$, we have for $z\in X$ and as $x'\to 0$
\[\begin{gathered}
R(\lam;z,x',y')= {x'}^{\ddemi+\lam}(L(\lam;z,y')+
O(x')),\\
\nabla_{x'\pl_{x'}}R(\lam;z,x',y')=
(\ddemi+\lam){x'}^{\ddemi+\lam}(L(\lam;z,y')+
O(x'))\\
\sigma_{\infty}(x',y')={x'}^{\ddemi-\lam}(\psi(y')+O(x')), \quad
\nabla_{x'\pl_{x'}}\sigma_{\infty}(x',y')=(\ddemi-\lam){x'}^{\ddemi-\lam}(\psi(y')+O(x'))
\end{gathered}\]  for some
$L\in C^{\infty}(X\x \pl\bar{X},\mc{E})$ and
where $\psi\in C^{\infty}(\pl\bar{X},\Sigma )$ is arbitrarily
chosen. We can then deduce that
$E(\lam;z,y')=2\lam L(\lam;z,y')$ as distributions in
$C^{\infty}(X\x \pl\bar{X};\mc{E})$. Using the
structure of $R(\lam)$ in Proposition \ref{extresolv}, we observe
that the kernel of $E(\lam)$ is also a distribution in
$C^{-\infty}(\bar{X}\x\pl\bar{X};\mc{E})$ since its
lift to $\bar{X}\x_0\pl\bar{X}$ is a conormal distribution on
$\bar{X}\x_0\pl\bar{X}$, more precisely it is an element in
$\rho_{\lb}^{\lam+\ddemi}\rho_{\ff}^{-\lam-\ddemi}C^{\infty}(\bar{X}\x_0\pl\bar{X};\mc{E})$. This is exactly the same argument for the
$E_\pm(\lam)$ formula in \eqref{kernelEla} by using the
representation \eqref{caspm} and integration by part.

Now for the scattering operator, we take {{$\Re(\lam)<-\ddemi$}}
and use the definition of $S(\lam)\psi$ to deduce that
\[S(\lam)\psi= (x^{-\lam-\ddemi}E(\lam)\psi)|_{x=0}.\]
From the fact that the lift of the kernel $x^{-\lam-\ddemi}E(\lam)$ to
$\bar{X}\x_0\pl\bar{X}$ is in
$\rho_{\ff}^{-2\lam-d}C^{\infty}(\bar{X}\x_0\pl\bar{X};{\rm
End}({^0\Sigma}))$, thus in $C^0(\bar{X}\x\pl\bar{X})$, we see
that
\[\int_{\pl\bar{X}}x^{-\lam-\ddemi}E(\lam;x,y,y')\psi(y'){\rm dv}_{h_0}(y')\]
is $C^0(\bar{X},{^0\Sigma})$ and its restriction to $\pl\bar{X}$
is given by the pairing of $[x^{-\lam-\ddemi}E(\lam)]|_{x=0}$ with
$\psi$, which ends the proof for $S(\lam)$. The argument is the
same for $S_\pm(\lam)$.
\end{proof}

In the same way as $E(\lam)$, we define the operator
$E^\sharp(\lam),E^\sharp_\pm(\lam)$ so that the Schwartz kernel of
$E^\sharp(\lam), E^\sharp_\pm(\lam)$ are given for
$\lam\in\cc\setminus(-\nn/2)$ not a pole of $R(\lam)$ by
\begin{equation}\label{esharp}
E^\sharp(\lam;y,m'):=2\lam[x^{-\lam-\ddemi}R(\lam;x,y,m')]|_{x=0}, \quad
E^\sharp_\pm(\lam;y,m'):=-{\rm cl}(\nu)[x^{-\lam-\ddemi}R_\pm(\lam;x,y,m')]|_{x=0}
\end{equation}
using the product decomposition
$[0,\eps)_x\x\pl\bar{X}$ near $\pl\bar{X}$.
Like for the analysis
of the kernel of $E(\lam)$ above, the structure of the kernel
$E^\sharp(\lam)$ on the blow-up $\pl\bar{X}\x_0\bar{X}$ is clear
from the analysis of $R(\lam)$. Note also that
$E_\pm^\sharp(\lam)$ maps $\dot{C}^{\infty}(\bar{X},{^0\Sigma})$
to $C^{\infty}(\pl\bar{X};{^0\Sigma}_\pm)$ by using Corollary
\ref{resolvforD} and
\[(E_\pm^{\sharp}(\lam)f)(y)=-{\rm cl}(\nu)\lim_{x\to 0}\int_{X}x^{-\lam-\ddemi}R_\pm(\lam;x,y,m')\sigma(m'){\rm dv}_g(m').\]
We see also from the remark following Corollary \ref{resolvforD} that
\begin{equation}\label{esharpestar}
E^\sharp(\bar{\lam};y,m')=E(\lam;m',y)^*, \qquad \qquad
E^\sharp_\pm(\bar{\lam};y,m')=E_\mp(\lam;m',y)^*.
\end{equation}
when these are considered as linear maps from
${^0\Sigma}_{m'}$ to ${^0\Sigma}_y$.

\begin{lemma}\label{functeq1}
Let $m,m'\in X$, then for $\lam\notin\zz/2$ neither a pole of
$R(\lam)$ nor of $R(-\lam)$, we have
\begin{equation}\begin{gathered}\label{Rla-R-la}
R(\lam;m,m')-R(-\lam;m,m')=(2\lam)^{-1}\int_{\pl\bar{X}} E(\lam;m,y)E^\sharp(-\lam;y,m'){\rm dv}_{h_0}(y),\\
R_\pm(\lam;m,m')-R_\mp(-\lam;m,m')=-\int_{\pl\bar{X}} E_\pm(\lam;m,y){\rm cl}(\nu) E_\mp^\sharp(-\lam;y,m') {\rm dv}_{h_0}(y)
\end{gathered}\end{equation}
or in terms of operators
\[R(\lam)-R(-\lam)=(2\lam)^{-1}E(\lam)E^\sharp(-\lam), \quad R_\pm(\lam)-R_\mp(-\lam)=-E_\pm(\lam){\rm cl}(\nu) E_\mp^\sharp(-\lam).\]
\end{lemma}
\begin{proof} This is a straightforward application of Green formula
and does not involve anything more than in the proof given by Guillop\'e
\cite{Guillope} for the Laplacian on functions on a surface. It is
based on the fact that
$(D^2+\lam^2)R(\lam;m,m')=(D^2+\lam^2)R(-\lam;m,m')=\delta(m-m')$
and $(D\pm i\lam)R_\pm(\lam;m,m')=(D\pm
i\lam)R_\mp(-\lam;m,m')=\delta(m-m')$, where $\delta(m-m')$
denotes the Dirac mass on the diagonal.
\end{proof}
A corollary of this is some functional equations relating
$E(\lam),E^\sharp(\lam)$ and $S(\lam)$.
\begin{corollary}\label{SvsE}
The following meromorphic identities hold
\[E(\lam)=E(-\lam)S(\lam),\quad \quad \ E^{\sharp}(\lam)=S(\lam)E^{\sharp}(-\lam)
 ,\]
\[E_\mp(\lam)=E_\pm(-\lam)S_\mp(\lam), \quad  \ E^\sharp_\pm(\lam)
=S_\pm(\lam)E_\mp^\sharp(-\lam). \]
\end{corollary}
\begin{proof}
Let us consider the second identity: assume $\Re(\lam)<-\ddemi$,
it suffices to multiply the first line of \eqref{Rla-R-la} by
$x(m)^{-\lam-\ddemi}$ and take the limit as $x(m)\to 0$ when
$m'\in X$ is fixed, the limit makes sense in view of our analysis
of the Schwartz kernels of $R(\lam), E(\lam), E^{\sharp}(\lam)$
and $S(\lam)$. Then we use Proposition \ref{kernels} and the
definition of $E^\sharp(\lam)$ and this gives the proof of the
second identity of the Corollary, at least for
$\Re(\lam)<-\ddemi$, but this extends meromorphically to
$\lam\in\cc$. The proofs of other identities are similar.
\end{proof}

\subsection{Properties of $S(\lam)$}

\begin{proposition}\label{spseudo}
For $\lam$ such that $\Re(\lam)<0$, $\lam\notin -\nn/2$ and $\lam$
not a pole of $R(\lam)$, the operator $S(\lam)$ is a classical
pseudo-differential operator on $\pl\bar{X}$ of order $2\lam$,
with principal symbol
\begin{equation}\label{defC}
\sigma_{\rm pr}(S(\lam))(\xi)=C(\lam){\rm cl}(\nu) |\xi|^{2\lam-1}_{h_0}{\rm
cl}(\xi), \textrm{ with
}C(\lam):=2^{-2\lam}\frac{\Gamma(-\lam+1/2)}{\Gamma(\lam+1/2)}.
\end{equation}
Moreover $S(\lam)$ can be
meromorphically extended to $\cc\setminus(-\nn/2)$ as a family of
pseudo-differential operators in $\Psi^{2\lam}(\pl\bar{X};\mc{E})$.
\end{proposition}
\begin{proof} Let $\beta:\bar{X}\x_0\bar{X}\to \bar{X}\x\bar{X}$ be the blow-down map,
$\pl\bar{X}\x_0\pl\bar{X}:=[\pl\bar{X},\pl\bar{X},\Delta_\pl]$ be
the blow-up of $\pl\bar{X}\x\pl\bar{X}$ around the diagonal
$\Delta_\pl$ and $\beta_\pl:\pl\bar{X}\x_0\pl\bar{X}$ the
associated blow-down map. Then the expression \eqref{kernelSla}
can also be written for $\Re(\lam)<0$ ($S(\lam)$ and $R(\lam)$
denote also the Schwartz kernel)
\begin{equation}\label{descriptionS}
S(\lam)=2\lam{\beta_\pl}_*\Big(\beta^{*}((xx')^{-\lam-\ddemi}R(\lam))|_{\lb\cap\rb}\Big)
\end{equation}
where $\lb\cap\rb$ is naturally identified with
$\pl\bar{X}\x_0\pl\bar{X}$. For more details, we refer to the
article of Joshi-Sa Barreto \cite{JSB} which deals with the
Laplacian on functions. Now using the fact that $R(\lam)\in
\Psi_0^{-2,\lam+\ddemi,\lam+\ddemi}(\bar{X},\mc{E})$, we
deduce that
\[((xx')^{-\lam-\ddemi}R(\lam))|_{\lb\cap\rb}\in \rho_{\ff,\pl}^{-2\lam-d}C^{\infty}(\pl\bar{X}\x_0\pl\bar{X},\mc{E})\]
where $\rho_{\ff,\pl}:=\rho_{\ff}|_{\lb\cap\rb}$ is a boundary
defining function of the boundary (i.e. the face obtained by
blowing-up) of $\pl\bar{X}\x_0\pl\bar{X}$. This shows that the
kernel $S(\lam)$ is classically (or polyhomogeneous) conormal to
the diagonal and the leading singularity at $y=y'=p$ given in
polar coordinates in the conormal bundle is given by
\[S(\lam;y,y')\sim c(\lam)|y-y'|^{-2\lam-d}U_p(p'), \quad p'=\frac{y-y'}{|y-y'|}\in S^{d-1} \]
for some $c(\lam)\in\cc$ and where $U_p(p')\in{\rm
End}({^0\Sigma}_p(\bar{X}))$ denote the limit of the parallel
transport $U(e_p,m)$ in the fiber $\ff_p$ when $m\to p'\in
\lb\cap\rb\cap \ff_p\simeq S^{d-1}$ (here $e_p$ is the center of
$\ff_p$ defined by the intersection of the interior diagonal with
$\ff_p$ and identified with the center of hyperbolic space). Thus
we obtain $S(\lam)\in\Psi^{2\lam}(\pl\bar{X},\mc{E})$, moreover
the expression \eqref{descriptionS} can be meromorphically
extended to $\cc\setminus(-\nn/2)$ as a distribution classically
conormal to the diagonal, thus as a family $S(\lam)\in
\Psi^{2\lam}(\pl\bar{X};\mc{E})$. As for the principal symbol, we
use the expression of $U_p(p')$ for $\hh^{d+1}$ given in Corollary
\ref{upp'} and Fourier transform to obtain \eqref{defC}. Notice
that there might be first order poles of infinite multiplicity at
$\nn/2$ coming from the meromorphic extension of the distribution
$|y-y'|^{2\lam+j}$ to $\lam\in \cc$. This phenomenon is described
in \cite{GRZ} for the case of $\Delta$ on functions.
\end{proof}

Like for functions, the scattering operator is a unitary operator
on the continuous spectrum and satisfies a functional equation
\begin{lemma}\label{fcteq}
The operator $S(\lam)$ is unitary on $\{\Re(\lam)=0\}$, it
satisfies $S(\lam)S(-\lam)={\rm Id}$ for $\lam$ such that $S(\pm
\lam)$ is defined, and it is conformally covariant in the sense
that for another choice $\hat{x}=e^{\omega}x$ of geodesic boundary
defining function, the corresponding scattering operator is
$\hat{S}(\lam)=e^{-(\lam+\ddemi)\omega_0}S(\lam)e^{(\ddemi-\lam)\omega_0}$,
where $\omega_0=\omega|_{\pl\bar{X}}$.
\end{lemma}
\begin{proof} The functional equation is a straightforward consequence of the uniqueness in Lemma
\ref{generalizedeig} or the first equality of Corollary
\ref{SvsE}.  The unitarity follows easily from Lemma
\ref{boundarypairing} by taking the solutions $\sigma_1,\sigma_2$
of Lemma \ref{generalizedeig} for two initial data $\psi_1,\psi_2
\in C^{\infty}(\pl\bar{X},{^0\Sigma})$. The conformal covariance
of $S(\lam)$ is straightforward by using the uniqueness of the
solution in Lemma \ref{generalizedeig}.
\end{proof}

\begin{corollary}\label{corol1}
If the metric $g$ is even in the sense of \eqref{evenness}, the operator
$\til{S}(\lam):=S(\lam)/C(\lam)$ is finite meromorphic in $\cc$, and it is holomorphic in
$\{\Re(\lam)\geq 0\}$.
\end{corollary}
\begin{proof}
The analyticity in the right half-plane is a consequence of the last statement in Lemma \ref{alternative} and the fact
that $S(\lam)\psi=\sigma^+|_{\pl\bar{X}}$ with the notation of this Lemma.
We already know the meromorphic extension outside $-\nn/2$ so we can write, using Proposition \ref{spseudo},
\[S(\lam)/C(\lam)={\rm cl}(\nu)
(D_{h_0}+i)(|D_{h_0}|+1)^{2\lam-1}(\mathrm{Id}+K(\lam))\] for some
$K(\lam)$ compact on $L^2(\pl\bar{X},{^0\Sigma})$ and analytic in
$\{\Re(\lam)\geq 0\}$. We know from Lemma \ref{fcteq} that
$\mathrm{Id}+K(\lam)$ is invertible for almost all
$\lambda\in\cc$, so we may use Fredholm analytic theorem to show
that $(S(\lam)/C(\lam))^{-1}$ is a meromorphic family of operators
with poles of finite multiplicity at most in $\Re(\lam)>0$, so by
the functional equation in Lemma \ref{fcteq}, we deduce that
$S(\lam)/C(\lam)$ is meromorphic in $\Re(\lam)<0$ with poles of
finite multiplicity.
\end{proof}

We give another corollary of the properties of $\til{S}(\lam)$.
\begin{corollary}\label{corol2}
For an AH manifold with a metric $g$ even in the sense of \eqref{evenness},
the resolvent $R_\pm(\lam)$ is finite meromorphic.
\end{corollary}
\begin{proof}
According to Proposition \ref{extresolv}, the only problem for the meromorphy of $R_\pm (\lam)$ can be at $-\nn/2$, so consider the half plane $\{\Re(\lam)<0\}$.
Since $R(-\lam),E(-\lam),E^{\sharp}(-\lam)$ are holomorphic in $\{\Re(\lam)<0\}$, the result
is a straightforward consequence of Corollary \ref{corol1} together with the formula
\[R(\lam)=R(-\lam) +(2\lam)^{-1}C(\lam)E(-\lam)\frac{S(\lam)}{C(\lam)}E^{\sharp}(-\lam),\]
itself a consequence of Lemma \ref{Rla-R-la} and Corollary \ref{SvsE}.
\end{proof}

\subsection{Representation of $E(\lam)$ and $S(\lam)$ in the case $X_\Gamma$}

Using Proposition \ref{kernels} and \eqref{RmoinsR0}, we obtain
directly an explicit representation modulo a smoothing term. We
use the functions $\eta_j, \phi^i_j,\chi^i_j$ of Section
\ref{caseconvexccompact}. We denote $S_{\hh^{d+1}}(\lam)$ and
$E_{\hh^{d+1}}(\lam)$ the scattering operator and Eisenstein
series on $\hh^{d+1}$, defined using a defining function of
$\pl\bbar{\hh^{d+1}}$ which is equal to $x_0$ on the half ball
$\bar{B}$ in the model
$\hh^{d+1}=\{(x_0,y_0)\in(0,\infty)\x\rr^{d-1}\}$, i.e. in terms
of distribution kernel on $\bar{B}\x \pl\bar{B}$ and on
$\pl\bar{X}\x\pl\bar{B}$
\begin{equation*}
\begin{gathered}
E_{\hh^{d+1}}(\lam;x_0,y_0,y'_0):=2\lam[{x'_0}^{-\lam-\ddemi}R_{\hh^{d+1}}(\lam;x_0,y_0,x'_0,y'_0)]_{x'_0=0},
\\
S_{\hh^{d+1}}(\lam;y_0,y'_0):=[x_0^{-\lam-\ddemi}E_{\hh^{d+1}}(\lam;x_0,y_0,y'_0)]|_{x_0=0}.\end{gathered}
\end{equation*}
\begin{lemma}\label{representations}
If $\lam\in \cc\setminus(-\nn_0/2)$ is not a pole of $R(\lam)$, then the Eisenstein
series $E(\lam)$ for $D^2$ on a convex co-compact quotient
$X:=\Gamma\backslash \hh^{d+1}$ has {the} kernel
$E(\lam)=E_0(\lam)+E_\infty(\lam)$ where
\[\begin{gathered}
E_0(\lam):=\sum_{j}\iota_j^*\chi_j^2 E_{\hh^{d+1}}(\lam)\phi_j^1\eta_j^{-\lam-\ddemi}{\iota_j}_*,\\
 E_\infty(\lam)=2\lam[{x'}^{-\lam-\ddemi}(R(\lam)-R_0(\lam))]|_{x'=0}\in x^{\lam+\ddemi}C^{\infty}(\bar{X}\x\pl\bar{X};\mc{E}).
\end{gathered}\]
Similarly, the scattering operator $S(\lam)$ for $D^2$ on $X$ has
the kernel $S(\lam)=S_0(\lam)+S_\infty(\lam)$ where
\[\begin{gathered}
S_0(\lam):=\sum_{j}\iota_j^*\eta_j^{-\lam-\ddemi}\phi_j^2 S_{\hh^{d+1}}(\lam)\phi_j^1\eta_j^{-\lam-\ddemi}{\iota_j}_*,\\
 S_\infty(\lam)=2\lam[(xx')^{-\lam-\ddemi}(R(\lam)-R_0(\lam))]|_{x=x'=0}\in C^{\infty}(\pl\bar{X}\x\pl\bar{X};\mc{E}).
\end{gathered}\]
\end{lemma}

\section{Selberg zeta function of odd type and Eta invariant}\label{geometricside}

In this section, we will assume that the dimension $d+1=2n+1$ is odd
except in Lemma \ref{l:r-dis}.

\subsection{Odd heat kernel of Dirac operator on $\hh^{d+1}$}
\label{ss-har} By the identification \eqref{e:homog},  the kernel
of the odd heat operator ${D}_{\hh^{d+1}} e^{-t{D}_{\hh^{d+1}}^2}$
 on $L^2(\hh^{d+1}, \Sigma(\hh^{d+1}))$ can be considered as
a $\tau_d$-radial function $P_t$ over $G$. Hence there exists a
function $P_t$ from $G$ to $\mathrm{End}(V_{\tau_d})$ satisfying
the $K$-equivariance condition
\begin{align}\label{e:equiv}
P_t(k_1 g k_2)=\tau_d(k_2)^{-1} P_t(g) \tau_d(k_1)^{-1} \qquad
\text{for} \quad g\in G, k_1,k_2\in K
\end{align} such that
\begin{align}\label{e:Pt}
{D}_{\hh^{d+1}}e^{-t{D}_{\hh^{d+1}}^2}(gK,hK)=P_t(h^{-1}g)
\qquad\text{for}\quad g,h\in G.
\end{align}
Let us remark that $P_t(h^{-1}g)$ and $P_t(k_1^{-1}h^{-1}gk_2)$
for $k_1,k_2\in K$  give the same map by the condition
\eqref{e:equiv}, so that the right hand side of \eqref{e:Pt} does
not depend on the choice of the representatives of the $K$-cosets.

Recalling the Cartan decomposition $G=KA^+K$ with $A^+:=\{
a_r=\exp(rH)\,|\, r>0\}$,  any element $g\in G$ can be written as
$g=ha_r k$ where $a_r=\exp(rH)$ and $r$ is the same as the
hyperbolic distance $d_{\hh^{d+1}}(eK, gK)$ between two points
$eK$ and $gK$ in $\hh^{d+1}\cong G/K$. Now let us recall that the
spin representation $\tau_d$ decomposes into two half spin
representations $\sigma_+$,$\sigma_-$ when restricting to $M={\rm
Spin}(d)$,
\[
\tau_d|_{M} = \sigma_+ \oplus \sigma_-,
\]
hence the representation space $V_{\tau_d}$ also decomposes into
$V_{\sigma^+}\oplus V_{\sigma^-}$ as $M$-representation spaces. By
Schur's lemma there exists a function
$p^\pm_t:\mathbb{R}\to\mathbb{C}$ such that
\[
P_t(a_r)|_{V_{\sigma_\pm}}= p^\pm_t(r) \,
\mathrm{Id}_{V_{\sigma_\pm}}
\]
where $a_r\in A^+$. As in the proof of Theorem 8.5 of \cite{CP}
using Theorem 8.3 of \cite{CP}, one can easily derive

\begin{proposition}\label{p:scalar}
The scalar components of $D_{\hh^{d+1}}e^{-tD_{\hh^{d+1}}^2}$ are given by
\begin{equation}\label{e:scalar}
p^\pm_t(r)=\pm\frac{\sinh(r/2)}{ i 2^{3n+3/2}
\Gamma(n+3/2)t^{3/2}} \left(-\frac{d}{d(\cosh r)}\right)^{n} r
\sinh^{-1} (r/2)\, e^{-\frac{r^2}{4t}}.
\end{equation}
\end{proposition}
Let us observe that the equalities \eqref{e:equiv} and
\eqref{e:scalar} determine the odd heat kernel
${D}_{\hh^{d+1}}e^{-t{D}_{\hh^{d+1}}^2}$ by the Cartan decomposition
$G=KA^+K$.

\subsection{Odd heat kernels over convex co-compact hyperbolic manifolds}
By a usual construction, the kernel of the odd heat operator
$De^{-tD^2}$ over $X_\Gamma$  is given (as an automorphic kernel) by
\begin{equation}\label{e:heatkernel}
De^{-tD^2}(g , h )=\sum_{\gamma\in\Gamma}
{D}_{\hh^{d+1}}e^{-t{D}_{\hh^{d+1}}^2}(g, \gamma h )
\end{equation}
where $g $, $h$ denote points in $X_\Gamma=\Gamma\backslash\hh^{d+1}$
which we view as a fundamental domain in $\hh^{d+1}$ with identified
sides through $\Gamma$. By \eqref{e:equiv}, \eqref{e:scalar} and some elementary
calculations,

\begin{equation}\label{e:odd-est}
|| {D}_{\hh^{d+1}}e^{-t{D}^2_{\hh^{d+1}}}(g, \gamma h) || \leq C\,
t^{-\frac32}\,
e^{-\ddemi r_\gamma(g,h)-\frac{r_\gamma(g,h)^2}{4t}}
\sum_{\substack{0\leq j\leq n+1\\ 0\leq k \leq
n}} r_\gamma(g,h)^j t^{-k} .
\end{equation}
Here $C$ is a positive constant independent on $t,r_\gamma(g,h)$
where $r_\gamma(g,h):=d_{\hh^{d+1}}(g,\gamma h)$.

\begin{lemma}\label{l:r-dis} Let $\mc{F}$ be a fundamental domain of $\Gamma$ and $\til{x}$ be a boundary defining function
of $X_\Gamma$ which we view as well as a function on $\mc{F}$. There are positive constants
$C_1,C_2,C_3$ such that for all $\gamma\in \Gamma$ with $r_\gamma:=d_{\hh^{d+1}}(e,\gamma e)>C_3$ and all $g,h\in \mc{F}$
\[ C_1\, e^{-r_\gamma}\til{x}(g)\til{x}(h) \leq e^{-r_\gamma(g,h)} \leq C_2\, e^{-r_\gamma}\til{x}(g)\til{x}(h).\]
\end{lemma}
\begin{proof} By conjugating by an isometry of $\hh^{d+1}$,
we can assume that in the half space model $\hh^{d+1}=\rr_x^+\x\rr_y^d$ the points
$(0,0)$ and $\infty$ are in the limit set $\Lambda(\Gamma)$ of the group $\Gamma$.
Then, since the group is convex co-compact,
a fundamental domain $\mc{F}$ satisfies the following: there exists $C>0$ and $\eps>0$ such that
\[\mc{F}\subset B:=\{z\in[0,\infty)\x\rr^d; |z|\leq C, d_{\rm eucl}(z;\Lambda(\Gamma))\geq \eps\}\]
where $d_{\rm eucl}$ denotes the Euclidean distance in $\rr^{d+1}$.
Let $\gamma\in \Gamma$ an isometry which fixes points $p^1_\gamma,p_\gamma^2\in\rr^d_y\cap B$.
Composing a translation $z\to z-(p^1_\gamma+p^2_\gamma)/2$
with a rotation in the $\rr^d_y$ variable, we define an isometry $q_\gamma$
which maps $p^1_\gamma$ to $p_\gamma:=(0,|p^1_\gamma-p_\gamma^2|/2,0\dots,0)$ and
$p^2_\gamma$ to $-p_\gamma$. Notice that
\[q_\gamma(\mc{F})\subset B':=\{z\in [0,\infty)\x\rr^d;
|z|\leq 2C, d_{\rm eucl}(z;\pm p_\gamma)\geq \eps\}.\]
Now we can define the isometry $s_\gamma$ obtained by composing the dilation $z\to |p_\gamma|^{-1}z$
with the map
\[(z_0,y_2,\dots,y_d)\to (\frac{1+z_0}{z_0-1},y_2,\dots,y_d)\]
where $z_0=(y_1+ix)\in\cc$ denotes complex coordinates for $\rr^+_{x}\x\rr_{y_1}$.
This isometry $s_\gamma$ maps $p_\gamma$ to $\infty$ and $-p_\gamma$ to $0$. Moreover,
it is easy to see that for $z\in B'$ we have
\[\frac{\eps}{4C}\leq |s_\gamma(z)|\leq \frac{4C}{\eps}.\]
We conclude that $t_\gamma:=s_\gamma\circ q_\gamma$ maps $\mc{F}$ to $\{z\in\rr^+\x\rr^d, \varepsilon\leq|z|\leq\varepsilon^{-1}\}$
for some $\varepsilon>0$ which do not depend on $\gamma$. But $t_\gamma\circ \gamma\circ t_\gamma^{-1}$
is an isometry fixing the line $\{y=0\}$ and thus can be written
under the form
\[t_\gamma\circ \gamma\circ t_\gamma^{-1}:(x,y)\to e^{r_\gamma}(x,\mc{O}_\gamma(y))\]
for some $\mc{O}_\gamma(y)\in {\rm SO}(d)$. Moreover, possibly by
composing $t_\gamma$ with an inversion, we can assume than
$r_\gamma>0$ and it is clearly given by $d_{\hh^{d+1}}(e,\gamma
e)$. Then we have for $m=(x,y),m'=(x',y')$ in the half-space model
\[\cosh^2 (d_{\hh^{d+1}}(m,m')/2) = \frac{|y-y'|^2+|x+x'|^2}{4xx'}\]
and by writing $t_\gamma g=(x,y)$ and $t_\gamma h=(x',y')$ in the half-space model,
\begin{equation}\label{e:r-dis}
\cosh^2(r_\gamma(g,h)/2) =  e^{r_\gamma} \frac{| e^{-r_\gamma} y
-\mc{O}_\gamma(y')|^2 +| e^{-r_\gamma} x+x'|^2}{4xx'}.
\end{equation}
But since $t_\gamma g,t_\gamma h\in t_\gamma(\mc{F})$, one has $\varepsilon\leq (x^2+|y|^2)^\demi\leq \varepsilon^{-1}$
and the same for $(x',y')$, which from \eqref{e:r-dis} implies that $\cosh^2(r_\gamma(g,h)/2)e^{-r_\gamma}xx'\in [\varepsilon^{2},\varepsilon^{-2}]$ for $r_\gamma$
large enough (depending only $\varepsilon$).
It remains to use the that $t_\gamma$ is an isometry of $\hh^{d+1}$
and so there exists $C>1$ such that $(x\circ t_\gamma)/\til{x}\in [C^{-1},C]$
on $\mc{F}$, which ends the proof.
\end{proof}

By \eqref{e:delta}, \eqref{e:odd-est} and Lemma \ref{l:r-dis}, we
get that the right hand side of \eqref{e:heatkernel} converges uniformly in $g,h$ in
a fundamental domain. We denote by ${\rm Tr}_{\Sigma}$ the local
trace over ${^0\Sigma}_m(X_\Gamma)\cong V_{\tau_d}$ for $m$ in a fundamental domain of $\Gamma$.

\begin{proposition}\label{p:odd-kernel} For $m$ in a fundamental domain of $\Gamma$, we have
\begin{equation}\label{e:odd-kernel}
\mathrm{Tr}_{\Sigma} ( D e^{-tD^2})(m)=
\sum_{\gamma\in\Gamma\setminus\{ e\}} \mathrm{Tr}_{\Sigma} (
{D}_{\hh^{d+1}} e^{-t{D}_{\hh^{d+1}}^2} )(m,\gamma m).
\end{equation}
\end{proposition}
\begin{proof} It is enough to show that $\mathrm{Tr}_{\Sigma}( {D}_{\hh^{d+1}} e^{-t{D}_{\hh^{d+1}}^2})(m,m)=0$,
which is a consequence of
\[\mathrm{Tr}_{\Sigma}( {D}_{\hh^{d+1}}e^{-t{D}_{\hh^{d+1}}^2})(m,m)=\mathrm{Tr}_{\Sigma}(P_t)(e)=d(\sigma_\pm)(p^+_t(0)+p^-_t(0))=0\]
where $d(\sigma_\pm)$ denotes the dimension of $V_{\sigma_\pm}$ and the last equality
follows from \eqref{e:scalar}.
\end{proof}
By equations \eqref{e:odd-est}, \eqref{e:r-dis} and Proposition
\ref{p:odd-kernel}, we deduce that there is $\eps>0$ such that
\[\big |\mathrm{Tr}_{\Sigma}( D e^{-tD^2})(m)\big| \leq
C_\epsilon(t)\, x(m)^{2(\ddemi+\epsilon)}\,
\sum_{\gamma\in\Gamma\setminus\{e\}} e^{-(\ddemi+ \epsilon) r_\gamma}\]
where $C_\epsilon(t)$ is a constant depending only on $\epsilon,t$ and $x$
a boundary defining function.
Hence the local trace function
$\mathrm{Tr}_{\Sigma}(De^{-tD^2})(m, m)$ is integrable over
$X_\Gamma$. Now we can define
\begin{equation}\label{e:trace}
\mathrm{Tr}(De^{-tD^2}):= \int_{X_\Gamma}
\mathrm{Tr}_{\Sigma}(De^{-t D^2})( m)\, {\rm dv}(m)
\end{equation}
where ${\rm dv}(m)$ denotes the metric over $X_\Gamma$ induced
from the hyperbolic metric ${\rm dv}_{\hh^{d+1}}$.

By our assumption on $\Gamma$,  $\Gamma\setminus\{e\}$ consists of
hyperbolic elements and decomposes into $\Gamma$-conjugacy classes
of hyperbolic elements. We denote by $\Gamma_{\mathrm{hyp}}$ the
set of $\Gamma$-conjugacy classes of hyperbolic elements. Each
element $[\gamma]$ in the set $\Gamma_{\mathrm{hyp}}$ corresponds
to a closed geodesic $C_\gamma$ in $X_\Gamma$. We denote by
$l(C_\gamma)$ the length of $C_\gamma$ and by $j(\gamma)$  the
positive integer such that $\gamma=\gamma_0^{j(\gamma)}$ with a
primitive $\gamma_0$. A primitive hyperbolic element $\gamma$
means that it can not be given by a power of any other elements in
$\Gamma$, so that $\Gamma$-conjugacy class of a primitive $\gamma$
corresponds to a prime geodesic $C_\gamma$ in $X_\Gamma$. The
trace of the monodromy in $\Sigma(X_\Gamma)\cong \Gamma\backslash
\big(G\times_{\tau_d} V_{\tau_d}\big)$ along a closed geodesic
$C_\gamma$ is given by $\chi_{\sigma_+}(m_\gamma) +
\chi_{\sigma_-}(m_\gamma)$ since any hyperbolic element $\gamma$
can be conjugated to $ m_\gamma a_\gamma\in MA^+$. A closed
geodesic $C_\gamma$ corresponds to a fixed point of the geodesic
flow on the unit sphere bundle over $X_\Gamma$. The Poincar\'e map
$P(C_\gamma)$ is the differential of the geodesic flow at
$C_\gamma$, which is given by $P(C_\gamma)=\mathrm{Ad}( m_\gamma
a_\gamma)$ if $\gamma= m_\gamma a_\gamma$. The unit sphere bundle
$S X_\Gamma$ of $X_\Gamma$ is given by $\Gamma\backslash G/M$,
and its tangent bundle $T S X_\Gamma$ is given by
\[TS X_\Gamma = \Gamma\backslash G\times_M
(\overline{\mathfrak{n}}\oplus\mathfrak{a}\oplus \mathfrak{n})
\]
where $\overline{\mathfrak{n}}=\theta(\mathfrak{n})$ and $M$ acts
on $\overline{\mathfrak{n}}\oplus\mathfrak{a}\oplus \mathfrak{n}$
by the adjoint action $\mathrm{Ad}$.  Hence $P(C_\gamma)$
preserves the decomposition
$\overline{\mathfrak{n}}\oplus\mathfrak{a}\oplus \mathfrak{n}$. We
denote by $P(C_\gamma)_{|\mathfrak{n}}$,
$P(C_\gamma)_{|{\overline{\mathfrak{n}}}}$ its restriction to
$\mathfrak{n}$, $\overline{\mathfrak{n}}$ part respectively. Now
we put
\begin{equation}\label{e:Dgamma}
D(\gamma):= | \det(
P(C_\gamma)_{|\mathfrak{n}\oplus\bar{\mathfrak{n}}}
-\mathrm{Id})|^{1/2} =e^{-nl(C_\gamma)}|\det(
P(C_\gamma)_{|\mathfrak{n}}-\mathrm{Id})|= e^{n l(C_\gamma)}\det(
\mathrm{Id}-P(C_\gamma)_{|\bar{\mathfrak{n}}}).
\end{equation}

\begin{proposition} \label{p:odd-heat}
The following identity holds
\begin{align}\label{e:odd-heat}
\mathrm{Tr} ( D e^{-tD^2}) =\frac{2\pi i}{(4\pi t)^{\frac32}}
\sum_{[\gamma]\in \Gamma_{\mathrm{hyp}}} \frac{l(C_\gamma)^2}
{j(\gamma) D(\gamma)} (
\chi_{\sigma_+}(m_\gamma)-\chi_{\sigma_-}(m_\gamma) )
e^{-\frac{l(C_\gamma)^2}{4t}}.
\end{align}
\end{proposition}

\begin{proof}
By equalities \eqref{e:Pt}, \eqref{e:odd-kernel} and
\eqref{e:trace},
\begin{equation}\label{e:tr-hyp}
\mathrm{Tr} ( D e^{-tD^2})=\sum_{\gamma\in\Gamma\setminus\{e\}}
\int_{\Gamma\backslash G} \mathrm{Tr}_{\Sigma}\, P_t(g^{-1}\gamma
g)\, d(\Gamma g).
\end{equation}
By Theorem 2.2 in \cite{M}, the scalar function
$p_t(g):=\mathrm{\Tr}_{\Sigma}\,P_t(g)$ is in the Harish-Chandra
$L^1$-space. (Note that $p_t(g)$ should not be confused with
$p_t^\pm(r)$ in the subsection \ref{ss-har}.) Hence we can follow
the well known path of Selberg on p.\ 63--66 of his famous paper
\cite{Sel} to obtain
\begin{equation}\label{e:hyp-orbit}
\sum_{\gamma\in\Gamma\setminus\{e\}} \int_{\Gamma\backslash G}
p_t(g^{-1}\gamma g)\, d(\Gamma g)=\sum_{[\gamma]\in
\Gamma_{\mathrm{hyp}}} \mathrm{vol}(\Gamma_\gamma\backslash
G_{\gamma}) \int_{G_\gamma\backslash G} p_t(g^{-1}\gamma g)\,
d(G_\gamma g)
\end{equation}
where $\Gamma_\gamma$, $G_\gamma$ denote the centralizer of  $\gamma$
in $\Gamma$ and $G$ respectively. Now we show the following
equality,
\begin{equation}\label{e:hyp5}
\mathrm{vol}(\Gamma_\gamma\backslash G_\gamma)\,\cdot
\int_{G_\gamma\backslash G} p_t(g^{-1}\gamma g)\, d(G_\gamma
g)=\frac{2\pi i}{(4\pi t)^{\frac32}}\frac{l(C_\gamma)^2}{
j(\gamma)D(\gamma)} (
\chi_{\sigma_+}(m_\gamma)-\chi_{\sigma_-}(m_\gamma)
)e^{-\frac{l(C_\gamma)^2}{4t}} .
\end{equation}
We may assume that a hyperbolic element $\gamma\in \Gamma$ has the
form $m_\gamma a_\gamma   \in MA^+$. If $\gamma\in MA^+$,
\begin{equation}\label{e:hyp1}
\int_{G_\gamma\backslash G} p_t(g^{-1}\gamma g)\, d(G_\gamma
g)=\mathrm{vol}(G_\gamma/A)^{-1}\,\int_{G/A} p_t(g \gamma
g^{-1})\, d(gA).
\end{equation}
We also have
\begin{equation}\label{e:hyp2}
\int_{G/A} p_t(g \gamma g^{-1})\, d(gA) = D(\gamma )^{-1} F_{p_t}(
m_\gamma a_\gamma)
\end{equation}
where the Abel transform of $p_t$ is given by
\[
F_{p_t}(  m_\gamma a_\gamma)= a_\gamma^\rho \int_{N}\int_K p_t(k\,
m_\gamma a_\gamma n
 k^{-1})\, dk\, dn
\] with $a_\gamma^\rho=\exp(n l(C_\gamma))$.
By Theorem 6.2 in \cite{Wal}, we have
\begin{equation}\label{e:hyp3}
\Theta_{\sigma,\lambda}(p_t)=\int_M\int^\infty_{-\infty} F_{p_t}(m
\exp(rH))\, \tr\, \sigma(m)\, e^{i\lambda r}\, dr\, dm.
\end{equation}
Here $\Theta_{\sigma,\lambda}(p_t)$ is defined by
\[
\Theta_{\sigma,\lambda}(p_t):=\Tr\,
\pi_{\sigma,\lambda}(p_t)=\Tr\,\int_G
p_t(g)\pi_{\sigma,\lambda}(g)\, dx
\]
and for $(\sigma, H_\sigma)\in \widehat{M}$ (where $\widehat{M}$
denotes the set of equivalence classes of irreducible unitary
representations of $M$)  and $\lambda \in
\mathfrak{a}_{\mathbb{C}}^*$ the principal representation
$\pi_{\sigma,\lambda}:=\mathrm{Ind}^G_{MAN}(\sigma \otimes
e^{i\lambda}\otimes \mathrm{Id})$ of $G$ acts on the space
$$
\Hh_{\sigma,\lambda}:=\{ \, f: G\to H_\sigma \, | \, f(xman)=
a^{-(i\lambda+\rho)} \sigma(m)^{-1} f(x), \, f|_K \in L^2(K)\, \}
$$
by the left translation $\pi_{\sigma,\lambda}(g)f(x)= f(g^{-1}x)$.
Applying the Fourier inversion theorem and the Peter-Weyl theorem
to the equality \eqref{e:hyp3}, we get
\begin{equation}\label{e:hyp4}
F_{p_t}( m_\gamma
a_\gamma)=\sum_{\sigma\in\widehat{M}}\overline{\tr\,
\sigma(m_\gamma)}\, \frac{1}{2\pi}\int^\infty_{-\infty}
\Theta_{\sigma,\lambda}(p_t) e^{-il(C_\gamma) \lambda}\, d\lambda.
\end{equation}
Now let observe that $\Theta_{\sigma,\lambda}(p_t)$ vanishes
unless $\sigma=\sigma_\pm$ since  $\tau_d|_M=\sigma_+\oplus
\sigma_-$. Moreover, we have
\begin{equation}\label{e:Fourier-Tr}
\Theta_{\sigma_\pm,\lambda}(p_t)= \pm \lambda e^{-t\lambda^2}
\end{equation}
as in Proposition 3.1 in \cite{P1} by (4.5) in \cite{MS1}.
Note that the analysis for this does not depend on $\Gamma$, but
is performed over $G$.
Combining \eqref{e:hyp1}, \eqref{e:hyp2}, \eqref{e:hyp4},
\eqref{e:Fourier-Tr} and observing that
$\mathrm{vol}(\Gamma_\gamma\backslash
G_\gamma)/\mathrm{vol}(G_\gamma/A)=l(C_\gamma)
/j(\gamma)$, we conclude
\begin{align*}\label{e:hyp5}
\mathrm{vol}(\Gamma_\gamma\backslash G_\gamma)\,
\int_{G_\gamma\backslash G} p_t(g^{-1}\gamma g)\, d(G_\gamma g)
=&\sum_{\sigma\in \widehat{M}} \frac{l(C_\gamma){\overline{\tr\,
\sigma(m_\gamma)}}}{j(\gamma) D(\gamma)}\,
\frac{1}{2\pi}\int^\infty_{-\infty} \Theta_{\sigma,\lambda}(h)
e^{-il(C_\gamma) \lambda}\, d\lambda\\
=&\ \frac{l(C_\gamma)\big(
\chi_{\sigma_+}(m_\gamma)-\chi_{\sigma_-}(m_\gamma)\big)}{j(\gamma)D(\gamma )}\,
\frac{1}{2\pi} \int^\infty_{-\infty} \lambda e^{-t\lambda^2}
e^{il(C_\gamma) \lambda} \, d\lambda\\
=&\ \frac{2\pi i}{(4\pi t)^{\frac32}} \frac{l(C_\gamma)^2\big(
\chi_{\sigma_+}(m_\gamma)-\chi_{\sigma_-}(m_\gamma)\big)}{ j(\gamma)
D(\gamma )}\,  \,
e^{-\frac{l(C_\gamma)^2}{4t}}.
\end{align*}
Taking the sum over $[\gamma]\in \Gamma_{\mathrm{hyp}}$ of this
equality and by \eqref{e:tr-hyp} and \eqref{e:hyp-orbit}, we
obtain \eqref{e:odd-heat}.
\end{proof}
From Proposition \ref{p:odd-heat}, putting $c:=\mathrm{min}
_{[\gamma]\in\Gamma_{\mathrm{hyp}}} l(C_\gamma)>0$ we obtain the
\begin{corollary}\label{c:asymp}
The following estimates hold
\[\mathrm{Tr}(De^{-tD^2}) =  O( e^{-c^2/t}) \, \text{ as }t\to 0;\quad
\mathrm{Tr}(De^{-tD^2})  = O(t^{-3/2}) \, \text{ as }
t\to\infty.\]
\end{corollary}

\subsection{Selberg zeta function of odd type}
We define the Selberg zeta functions attached to half spinor
representations $\sigma_\pm$ by
\begin{equation}\label{e:def-Sel-halh}
Z_\Gamma({\sigma}_\pm,\lambda):=\exp\Big(-\sum_{[\gamma]\in\Gamma_{\mathrm{hyp}}}
\frac{\chi_{\sigma_\pm}(m_\gamma)}{j(\gamma)D(\gamma)}
e^{-\lambda\, l(C_\gamma)}\Big)
\end{equation}
for $\Re(\lambda)>\delta_\Gamma$. It is
easy to see that $Z_\Gamma(\sigma_\pm,\lambda)$ absolutely
converges for $\Re(\lam)>\delta_\Gamma$ by \eqref{e:delta}.

\begin{proposition}\label{p:zeta-exp} For $\Re(\lambda)>\delta_\Gamma$,
\begin{equation}\label{e:zeta-exp}
Z_\Gamma(\sigma_\pm,\lambda)=\prod_{[\gamma]\in
\mathrm{P}\Gamma_{\mathrm{hyp}}} \prod_{k=0}^\infty \mathrm{det}
\left(\mathrm{Id}-{\sigma}_\pm(m_\gamma)\otimes
S^k(P(C_\gamma)_{|\overline{\mathfrak{n}}})
e^{-(\lambda+n)l(C_\gamma)}\right)
\end{equation}
where $\mathrm{P}\Gamma_{\mathrm{hyp}}$ is the set of
$\Gamma$-conjugacy classes of primitive hyperbolic elements, and
for an endomorphism $L:V\to V$, $S^k(L)$ denotes the action
of $L$ on the symmetric tensor product $V^{\otimes
k}_{\mathrm{sym}}$.
\end{proposition}

\begin{proof} It is easy to see that $\log$ of the right hand side
\eqref{e:zeta-exp} is the same as
\begin{align*}
& \sum_{[\gamma]\in \mathrm{P}\Gamma_{\mathrm{hyp}}}
\sum_{k=0}^\infty\ \mathrm{Tr}\, \log
\left(\mathrm{Id}-\sigma_\pm(m_\gamma)\otimes
S^k(P(C_\gamma)_{|{\overline{\mathfrak{n}}}})
e^{-(\lambda+n)l(C_\gamma)}\right)\\
=& -\sum_{[\gamma]\in \mathrm{P}\Gamma_{\mathrm{hyp}}}
\sum_{k=0}^\infty \ \sum_{j=1}^\infty\ j^{-1} \mathrm{Tr}
\left(\sigma_\pm(m_\gamma)\otimes
S^k(P(C_\gamma)_{|{\overline{\mathfrak{n}}}})
e^{-(\lambda+n)l(C_\gamma)}\right)^j\\
=&-\sum_{[\gamma]\in \Gamma_{\mathrm{hyp}}} \sum_{k=0}^\infty\
j(\gamma)^{-1} \mathrm{Tr} (\sigma_\pm(m_\gamma)) \, \mathrm{Tr}(
S^k(P(C_\gamma)_{|{\overline{\mathfrak{n}}}}))
e^{-(\lambda+n)l(C_\gamma)}\\
=& -\sum_{[\gamma]\in
\Gamma_{\mathrm{hyp}}}j(\gamma)^{-1}\mathrm{det}(
\mathrm{Id}-P(C_\gamma)_{|{\overline{\mathfrak{n}}}})^{-1}
\mathrm{Tr} (\sigma_\pm(m_\gamma)) \, e^{-(\lambda+n)l(C_\gamma)}.
\end{align*}
Now these equalities complete the proof if we use
\eqref{e:Dgamma}.

\end{proof}

The Selberg zeta function of odd type is defined by
\begin{equation}\label{e:def-Sel-odd}
Z^o_{\Gamma,\Sigma}(\lambda)=\frac{Z_\Gamma(\sigma_+,\lambda)}{Z_{\Gamma}(\sigma_-,\lambda)}
\qquad\ \text{for}\quad \RE(\lambda)>\delta_\Gamma.
\end{equation}
Note that the definition in \eqref{e:def-Sel-odd} is shifted by
$-n$ from the one in \cite{M}, \cite{P1}.
 From this definition, the following equality follows easily
\begin{equation}\label{e:zeta-resolvent}
\pl_{\lam}\log Z_{\Gamma,\Sigma}^o(\lambda)=\sum_{[\gamma]\in
\Gamma_{\mathrm{hyp}}} l(C_\gamma) j(\gamma)^{-1} D(\gamma)^{-1} (
\chi_{\sigma_+}(m_\gamma)-\chi_{\sigma_-}(m_\gamma) ) e^{-\lambda
l(C_\gamma)}.
\end{equation}
By Proposition \ref{p:odd-heat},  and from the identity
\[
\int^\infty_0 e^{-t\lambda^2 } (4\pi t)^{-\frac32}
e^{-\frac{r^2}{4t}} \, dt = \frac{e^{-\lambda r}}{4\pi r}
\]
valid for $\RE(\lam^2)>0,$  we have
\begin{corollary}\label{c:resolvent-eq}
For $\RE(\lambda)>\max(\delta_\Gamma,0)$,
\begin{align}\label{e:resolvent-eq}
\int^\infty_0 e^{-t\lambda^2} \mathrm{Tr}(De^{-tD^2})\, dt
=&\frac{ i}{2} \partial_{\lam} \log Z_{\Gamma,\Sigma}^o(\lambda).
\end{align}
\end{corollary}

Let us define the function $\omega_\lam$ on the fundamental
domain $X_\Gamma$ by
\begin{equation*}
\label{defomega} \omega_\lam(m):= {\rm
Tr}_{\Sigma}[DR(\lam;m,m')-DR_{\hh^{d+1}}(\lam;m,m')]_{m=m'}
\end{equation*}
where $R_{\hh^{d+1}}(\lam)$ is the meromorphic extension of the
resolvent of $D_{\hh^{d+1}}^2$. The kernel of
$D(R(\lam)-R_{\hh^{d+1}}(\lam))$ is smooth in $X_\Gamma\x X_\Gamma$
and, since automorphic on $\hh^{d+1}$ with respect to $\Gamma$, it
induces a smooth kernel on $X\x X$. From the analysis of the
resolvent above, we see that it is in the class
$(xx')^{\lam+\ddemi}C^{\infty}(\bar{X}\x\bar{X};\mc{E})$
so its restriction to the diagonal followed by the local trace is
in $x^{2\lam+d}C^{\infty}(\bar{X})$ and thus integrable when
$\Re(\lam)>0$; moreover it is meromorphic in $\cc\setminus
(-\nn/2)$.
By reversing the order of integration and trace in
\eqref{e:resolvent-eq},
we can write for $\Re(\lam)>\max(\delta_{\Gamma},0)$
\begin{equation}\label{extensionZ}
\pl_\lam\log\,
Z_{\Gamma,\Sigma}^o(\lam)=\frac{\pl_{\lam}Z^o_{\Gamma,\Sigma}(\lam)}{Z_{\Gamma,\Sigma}^{o}(\lam)}=-2i\int_{X_\Gamma}
\omega_\lam(m)\, {\rm dv}(m)
\end{equation}
and the integral of $\omega_\lam(m)$ can be decomposed under the
form $\int_{x>\eps_0}$ and $\int_{x<\eps_0}$ for some boundary
defining function $x$, so that it can be decomposed under the sum
of a meromorphic function of $\lam$ and of
\[\lim_{\eps\to 0}\int_{\eps}^{\eps_0}x^{2\lam-1}\omega'_\lam(x,y)\,dx\,{\rm dv}_{h_0}(y)\]
for some $\omega'_\lam$ smooth and meromorphic in $\lam$. As
$\eps\to 0$, this has an expansion of the form $
A(\lam)+\sum_{j=0}^\infty \eps^{2\lam+j}C_j(\lam)$ for some
meromorphic $C_j(\lam)$, $A(\lam)$,  and for
$\lam>\max(\delta_\Gamma,0)$ we have \eqref{extensionZ} which is
equal to $A(\lam)+\int_{x(m)>\eps_0}\omega_\lam(m)\,{\rm dv}(m)$.
This shows
\begin{lemma}\label{formdlogzeta}
The function
$\pl_{\lam}Z^o_{\Gamma,\Sigma}(\lam)/Z_{\Gamma,\Sigma}^{o}(\lam)$
has a meromorphic extension to $\cc\setminus(-\nn/2)$ given by the
value
\[\frac{\pl_{\lam}Z^o_{\Gamma,\Sigma}(\lam)}{Z_{\Gamma,\Sigma}^{o}(\lam)}=-2i\, {\rm FP}_{\eps\to 0}
\int_{x(m)>\eps} \omega_\lam(m)\, {\rm dv}(m)\] where ${\rm
FP}_{\eps \to 0}$ means the finite part as $\eps \to 0$, that is
the constant coefficient in the expansion in powers of $\eps$, and
$\omega_\lam$ is given in \eqref{defomega}.
\end{lemma}

In addition, using that ${\rm
Tr}_\Sigma[DR_{\hh^{d+1}}(\lam;m,m')-DR_{\hh^{d+1}}(-\lam;m,m')]_{m=m'}=0$
for all $\lam\in\cc$, we obtain
\begin{equation}\label{zetavsresolvent}
\frac{\pl_{\lam}Z^o_{\Gamma,\Sigma}(\lam)}{Z_{\Gamma,\Sigma}^{o}(\lam)}-\frac{\pl_\lam
Z_{\Gamma,\Sigma}^o(-\lam)}{Z_{\Gamma,\Sigma}^o(-\lam)}= -2i\,
{\rm FP}_{\eps\to 0} \int_{x(m)>\eps} {\rm
Tr}_{\Sigma}[D\Pi(\lam;m,m')]_{m=m'}\, {\rm dv}(m)
\end{equation}
where $\Pi(\lam;m,m'):=R(\lam;m,m')-R(-\lam;m,m')$.

\subsection{Eta invariant} We define the eta invariant of $D$ by
\begin{equation}\label{e:def}
\eta(D)=\frac{1}{\sqrt{\pi}}\int^\infty_0 t^{-\frac12}
\mathrm{Tr}(De^{-tD^2})\, dt
\end{equation}
where the trace ${\rm Tr}$ means the integral of the local trace
like in \eqref{e:trace}. Note that the integral on the right hand
side of \eqref{e:def} is finite by Corollary \ref{c:asymp}.
The meromorphic extension of $Z^{o}_{\Gamma,\Sigma}(\lam)$ and its analyticity on $[0,\infty)$ will be proved
in next section, it implies directly the following result:
\begin{theorem}\label{t:first main thm}
The eta invariant of the Dirac operator over a convex co-compact hyperbolic
manifold $X_\Gamma$ satisfies
\begin{equation}
\exp({i\pi \eta(D)})= Z_{\Gamma,\Sigma}^o(0).
\end{equation}
\end{theorem}
\begin{proof}
We start by writing
\[t^{-1/2}=\frac{2}{\sqrt{\pi}}\int_0^\infty e^{-\lam^2t}d\lam,\]
then we have by \eqref{e:resolvent-eq}
\[\eta(D)=\frac{2}{\pi}\int_0^\infty\int_0^\infty e^{-\lam^2t}{\rm Tr}(De^{-tD^2})\, d\lam\, dt=\frac{i}{\pi}\int_0^\infty
\frac{\pl_\lam
Z_{\Gamma,\Sigma}^o(\lam)}{Z_{\Gamma,\Sigma}^o(\lam)}d\lam\] and
this concludes the proof by Theorem \ref{mainth}, in particular,
the meromorphic extension of $Z^o_{\Gamma,\Sigma}(\lam)$ over
$\cc$ with $\lam=0$ as a regular value.
\end{proof}

\section{Spectral side of trace formula, Maass-Selberg relation}

With the only exception of Theorem \ref{mainth},
the dimension $d+1$ in this Section can be either odd or even, and
$X$ can be any asymptotically hyperbolic manifold with constant curvature
near infinity.
By convention, if $J(\lam)$ is an operator depending on $\lam$, we
shall use the following notation throughout the paper
\[\til{J}(\lam):=J(\lam)/C(\lam), \quad {\rm with }\quad C(\lam)=2^{-2\lam}\frac{\Gamma(-\lam+1/2)}{\Gamma(\lam+1/2)} \]
where the constant $C(\lam)$, already introduced in \eqref{defC}, satisfies
$C(\lam)C(-\lam)=1$.

\subsection{The Maass-Selberg relation}
We now describe the Maass-Selberg relation in order to study the
singularities of the odd Selberg zeta function
in terms of scattering data.

A corollary of the Lemma \ref{functeq1} is that the kernel of
$\Pi(\lam):=R(\lam)-R(-\lam)$ is smooth on $X\x X$. Actually, in
the Mazzeo-Melrose construction described before, one can choose
the same term $Q_0(\lam)$ for the parametrix of $R(\lam)$ and
$R(-\lam)$, proving directly that $\Pi(\lam)$ is the sum of a term
whose lift under $\beta$ is smooth on $\bar{X}\x_0\bar{X}\setminus
(\lb\cup\rb)$ with a term in
$(xx')^{\lam+\ddemi}C^{\infty}(\bar{X}\x\bar{X};\mc{E})+
(xx')^{-\lam+\ddemi}C^{\infty}(\bar{X}\x\bar{X};\mc{E})$. The
local trace of $\Pi(\lam)$, i.e. the trace of the endomorphism
$\Pi(\lam;m,m)$, is thus in
\[\tra(\Pi(\lam;m,m))\in C^{\infty}(\bar{X})+x^{2\lam+d}C^{\infty}(\bar{X})+x^{-2\lam+d}C^{\infty}(\bar{X}).\]
From the composition properties of $\Psi_0^{*,*,*}(\bar{X};
\mc{E})$ in Mazzeo \cite{MazzeoCPDE}, the operator $D\Pi(\lam)$
has a kernel which has the exact same properties as $\Pi(\lam)$
and thus its local trace satisfies
\[\tra((D\Pi)(\lam;m,m))\in C^{\infty}(\bar{X})+x^{2\lam+d}C^{\infty}(\bar{X})+x^{-2\lam+d}C^{\infty}(\bar{X}).\]

\begin{lemma}\label{maassselberg2}
Let $\eps>0$ and $\lam\in \cc\setminus (\zz/2)$ neither a pole of
$R(\lam)$ nor of $R(-\lam)$, then
\begin{align*}
&\lefteqn{\int_{x(m)>\eps}{\rm Tr}_{\Sigma}(D\Pi)(\lam;m,m){\rm dv}_g(m)}\\
=&-\frac{1}{2} \eps^{-d}\int_{x(m)=\eps}\int_{\pl\bar{X}} {\rm
Tr}_{\Sigma}\Big({\rm
cl}(\nu)\pl_\lam\til{E}(\lam;x,y,y')\til{E}^\sharp(-\lam;y',x,y)\Big)
{\rm dv}_{h_0}(y'){\rm dv}_{h_\eps}(y).
\end{align*}
\end{lemma}
\begin{proof}
First observe that
\[D\til{E}(\lam)=-\lam \til{E}(\lam){\rm cl}(\nu),\quad \til{E}^\sharp(-\lam)D=-\lam {\rm cl}(\nu)\til{E}^\sharp(-\lam)\]
which is a consequence of Lemma \ref{alternative} and the remark
that follows. Then we get for $t\in\cc$ small
\[\frac{1}{2t}\Big(D\til{E}(\lam+t)\til{E}^\sharp(-\lam)-\til{E}(\lam+t)(\til{E}^\sharp(-\lam)D)\Big)=-\demi
E(\lam+t){\rm cl}(\nu)\til{E}^\sharp(-\lam).\] From Lemma
\ref{functeq1}, the limit as $t\to 0$ on the right hand side is
$D\Pi(\lam)=-\demi \til{E}(\lam){\rm cl}(\nu)
\til{E}^\sharp(-\lam)$, which by taking the local trace and using
the fact that $\tra(AB)=\tra(BA)$ gives
\[{\rm Tr}_{\Sigma}(D\Pi(\lam;m,m))=\demi\lim_{t\to 0}\frac{1}{t}\tra_\Sigma\Big(\til{E}^\sharp(-\lam)D\til{E}(\lam+t)(m,m)
-(\til{E}^\sharp(-\lam)D)\til{E}(\lam+t)(m,m)\Big).\] In
particular, remark that the local trace on the right hand side has
to vanish at $t=0$. We use Green formula on $\{x>\eps\}$
\begin{align*}
&\int_{x(m)>\eps}{\rm Tr}_{\Sigma}(D\Pi(\lam;m,m)){\rm dv}_g(m) \\
=&-\frac{\eps^{-d}}{2} \lim_{t\to 0}\frac{1}{t}\int_{x(m)=\eps}\int_{\pl\bar{X}}
\tra_{\Sigma}(\til{E}^\sharp(-\lam;y',m)){\rm cl}(\nu)\til{E}(\lam+t;m,y'))\,{\rm dv}_{h_0}(y')\,{\rm dv}_{h_\eps}(y)\\
=&-\frac{\eps^{-d}}{2} \int_{x(m)=\eps}\int_{\pl\bar{X}}
\tra_{\Sigma}(\til{E}^\sharp(-\lam;y',m)){\rm
cl}(\nu)\pl_{\lam}\til{E}(\lam;m,y'))\,{\rm dv}_{h_0}(y')\, {\rm
dv}_{h_\eps}(y).
\end{align*}
The Lemma is proved using the cyclicity of the local trace once again.
\end{proof}

In a second time, we want to consider the expression
in Lemma \ref{maassselberg2} as $\eps\to 0$. For that we
introduce the representation of $\til{E}(\lam)$ given in Lemma
\ref{representations} and a similar one for
$\til{E}^\sharp(\lam)$:
$\til{E}^\sharp(\lam)=\til{E}^{\sharp}_0(\lam)+\til{E}^\sharp_{\infty}(\lam)$
obtained by restricting \eqref{Rtransp2} times $x(m)^{-\lam-\ddemi}$ at
$m\in \pl\bar{X}$ and using \eqref{Rtransp1}, that is
\begin{equation}\label{esharptil}
\begin{gathered}
\til{E}_0^\sharp(\lam):=\sum_j\iota_j^*\phi_j^1\eta_j^{-\lam-\ddemi}\til{E}^\sharp_{\hh^{d+1}}(\lam)\chi^2_j{\iota_j}_*,\\
\til{E}_\infty^\sharp(\lam):=2\lam[x^{-\lam-\ddemi}(\til{R}(\lam)-\til{R}_0^\sharp(\lam))]_{x=0}\in
{x'}^{\lam+\ddemi}C^{\infty}(\pl\bar{X}\x\bar{X};\mc{E})
\end{gathered}
\end{equation}
where again $\til{E}^\sharp_{\hh^{d+1}}(\lam)$ is the
corresponding operator on $\hh^{d+1}$ like in Lemma
\ref{representations} and $R^\sharp_0(\lam):=R_0(\bar{\lam})^*$.
Notice that, using the same arguments as in \eqref{esharpestar}, we obtain
$\til{E}_0^\sharp(\lam;y',m)=\til{E}_0(\bar{\lam};m,y')^*$.
Similarly we have
$\til{S}(\lam)=\til{S}_0^{\sharp}(\lam)+\til{S}^{\sharp}_{\infty}(\lam)$
with
\begin{equation}\label{ssharptil}
\begin{gathered}
\til{S}_0^\sharp(\lam):=\sum_j\iota_j^*\phi_j^1\eta_j^{-\lam-\ddemi}\til{S}_{\hh^{d+1}}(\lam)\eta_j^{-\lam-\ddemi}\phi^2_j{\iota_j}_*,\\
\til{S}_\infty^\sharp(\lam):=2\lam[{xx'}^{-\lam-\ddemi}(\til{R}(\lam)-\til{R}_0^\sharp(\lam))]_{x=x'=0}\in
C^{\infty}(\pl\bar{X}\x\pl\bar{X};\mc{E}).
\end{gathered}\end{equation}
and $\til{S}_{\hh^{d+1}}(\lam)$ is the operator on $\hh^{d+1}$ like in
Lemma \ref{representations}. Then we can prove

\begin{proposition}\label{firstred}
The meromorphic identity holds in $\lam\in \cc\setminus (\zz/2)$,
\begin{align}\label{sommedesE}
&{\int_{x(m)>\eps}{\rm Tr}_\Sigma(D\Pi(\lam;m,m)){\rm
dv}_{g}(m)}\notag\\ =&-\frac{1}{2}
\eps^{-d}\int_{x(m)=\eps}\int_{\pl\bar{X}}\Big[{\rm
Tr}_\Sigma\Big({\rm cl}(x\pl_x)\pl_\lam\til{E}(\lam;x,y,y')\til{E}_\infty^\sharp(-\lam;y',x,y)\Big)\\
&\hspace{3cm} +{\rm Tr}_\Sigma\Big({\rm cl}(x\pl_x)\pl_\lam
\til{E}_\infty(\lam;x,y,y')\til{E}^\sharp_0(-\lam;y',x,y)\Big)
\Big]{\rm dv}_{h_0}(y'){\rm dv}_{h_\eps}(y)\notag \end{align}
where ${\rm Tr}_\Sigma$ means the local trace on ${\rm
End}({^0\Sigma})$.
\end{proposition}
\begin{proof}
The point is to prove the vanishing of
\[\begin{gathered}
{\rm Tr}_\Sigma\Big({\rm cl}(\nu)\pl_\lam\til{E}_0(\lam;x,y,y')
\til{E}_0^\sharp(-\lam;y',x,y)\Big)
\end{gathered}\]
so we use the explicit formula for $\til{E}_0(\lam)$ and
$\til{E}_0^\sharp(\lam)$ given in Lemma \ref{representations} and
\eqref{esharptil}. We have to deal with terms of the
form
\begin{equation}\label{2terms}
{\rm Tr}_\Sigma\Big(\iota_j^*{\rm cl}(X_j)\pl_\lam
(\chi^2_j\til{E}_{\hh^{d+1}}(\lam)\phi^1_j\eta_j^{-\lam-\ddemi})\gamma_{jk}^*
\phi^1_k\eta_k^{\lam-\ddemi}\til{E}_{\hh^{d+1}}^\sharp(-\lam)\chi^2_k{\iota_k}_*(m,m)\Big)
\end{equation}
where $\gamma_{jk}$ is the unique isometry of $\hh^{d+1}$
extending $\iota_k\circ\iota_j^{-1}:\iota_j(U_k\cap U_j)\to
\iota_k(U_k\cap U_j)$ (and which acts smoothly up to the
boundary), $X_j$ is the vector field $X_j:={\iota_j}_*(x\pl_x)$
and is the connection on $^0\Sigma(\bar{B})$. We  use the fact
that
$\gamma_{jk}^*R_{\hh^{d+1}}(\lam)=R_{\hh^{d+1}}(\lam)\gamma_{jk}^*$
since $\gamma_{jk}$ is an isometry so if
$\alpha_{jk}:=[\gamma_{jk}^*(x_0)/x_0]|_{x_0=0}\in
C^\infty(\rr^{d})$, then one deduces that
$\gamma_{jk}^*E^\sharp_{\hh^{d+1}}(-\lam)=\alpha_{jk}^{\lam-\ddemi}E^\sharp_{\hh^{d+1}}(-\lam)\gamma_{jk}^*$.
 Let us consider \eqref{2terms}, it can be written as
\[\begin{split}
A_{jk}(\lam;x_0,y_0,y_0'):=&\chi^2_j(x_0,y_0)\chi^2_k(\gamma_{jk}(x_0,y_0))\\
&\times{\rm Tr}_\Sigma\Big( {\rm cl}(X_j)(\pl_{\lam}[\til{E}_{\hh^{d+1}}(\lam;x_0,y_0,y_0')\eta_j^{-\lam}]\beta_{jk}(\lam;y_0')
\til{E}_{\hh^{d+1}}^\sharp(-\lam;y_0',x_0,y_0)\Big)
\end{split}\]
where $\iota_j(m)=(x_0,y_0)$ and
$\beta_{jk}(y',\lam):=\eta_j(y')^{\lam-d}\phi^1_j(y')\phi^1_k(\gamma_{jk}(y'))$.  We shall show that $A_{jk}$ vanishes for
algebraic reasons. First we recall from \eqref{resolvforHn} that
\[\til{E}_{\hh^{d+1}}(\lam;x_0,y_0,y_0')=f(\lam)x^{\lam}(x_0^2+|y_0-y'_0|^2)^{-\lam} U((x_0,y_0),y_0')\]
and a similar expression for $\til{E}_{\hh^{d+1}}^\sharp$, here $U$ is the parallel transport on $\hh^{d+1}$ extended to the
boundary (see Appendix \ref{paralleltr}) and $f(\lam)$ some explicit meromorphic function.
Thus using the fact that $U(m_0,y_0')U(y_0',m_0)={\rm Id}$, $A_{jk}$ can be written under the form
\[A_{jk}(\lam;m_0,y'_0)= b_{jk}(\lam;m_0,y_0'){\rm Tr}_\Sigma({\rm cl}(X_j)U(m_0,y_0') U(y_0',m_0))=
b_{jk}(\lam;m_0,y_0'){\rm Tr}_\Sigma({\rm cl}(X_j))
\]
for some $b_{jk}$ where $m_0=(x_0,y_0)$. But since the dimension $d>1$, the trace vanishes.
\end{proof}

We deduce from this formula
\begin{proposition}\label{secondred}
For $\lam\in\cc\setminus (\zz/2)$ not a pole of $S(\lam)$ and
$S(-\lam)$, the term  in Proposition \ref{firstred} has a limit as
$\eps\to 0$, given by
\begin{equation}
\begin{gathered}\label{sommedesS}
\lim_{\eps\to 0}\int_{x(m)>\eps}{\rm Tr}_\Sigma(D\Pi(\lam;m,m)){\rm dv}_{g}(m)
=-\frac{1}{2} {\rm Tr}\left(
{\rm cl}(\nu)[\pl_\lam \til{S}(\lam)\til{S}^\sharp_{\infty}(-\lam)+\pl_\lam
\til{S}_\infty(\lam)\til{S}^\sharp_0(-\lam)]\right).
\end{gathered}
\end{equation}
where ${\rm Tr}$ is the usual trace for trace class operators.
\end{proposition}
\begin{proof}
First when $d+1$ is odd, we know from the discussion before Lemma
\ref{formdlogzeta} that the term \eqref{sommedesE} has an
expansion as $\eps\to 0$ of the form
$A(\lam)+\sum_{j=0}^\infty\eps^{-2\lam}C^-_j(\lam)+\sum_{j=0}^\infty\eps^{2\lam+j}C^{+}_j(\lam)$
for some meromorphic functions $A(\lam),C^\pm_j(\lam)$.  But actually
the same result holds for the general AH manifolds where the
metric has constant curvature near $\infty$ and $d+1$ odd or even:
indeed, using the parametrix \eqref{RmoinsR0} and the fact that
the local trace ${\rm
Tr}_{\Sigma}(D(R_{\hh^{d+1}}(\lam)-R_{\hh^{d+1}}(-\lam)))=0$ as
explained in the Remark following Proposition \ref{p:CP}, it is
clear that ${\rm Tr}_{\Sigma}(D\Pi(\lam;m,m))$ is a function in
$x^{2\lam}C^{\infty}(\bar{X})+x^{-2\lam}C^{\infty}(\bar{X})$. Let
us then take the limit as $\eps\to 0$ in \eqref{sommedesE}. For
instance consider
\[ x^{-d}{\rm Tr}_\Sigma\Big({\rm cl}(\nu)\pl_\lam \til{E}(\lam;x,y,y'))
\til{E}_\infty^\sharp(-\lam;y',x,y)\Big),\]
we can use the arguments used in the proof of Theorem 3.10 of \cite{GuiAJM} 
(in the present case this is even simpler since they correspond only 
to the mixed terms there, which comes from the regularity 
$\til{E}^\sharp_\infty(-\lam)\in 
x^{-\lam+\ddemi}C^{\infty}(\pl\bar{X}\x\bar{X};\mc{E})$) and we obtain
\begin{multline}\label{termelog1}
x^{-d}{\rm Tr}_{\Sigma}\Big({\rm cl}(\nu)\pl_\lam\til{E}(\lam;x,y,y')
\til{E}_\infty^\sharp(-\lam;y',x,y)\Big) \\
=\log(x){\rm Tr}_\Sigma\Big({\rm cl}(\nu)\til{S}(\lam;y,y')
\til{S}^\sharp_\infty(-\lam;y',y)\Big)
+{\rm Tr}_\Sigma\Big({\rm cl}(\nu)\pl_\lam\til{S}(\lam;y,y')
\til{S}^\sharp_\infty(-\lam;y',y)\Big)\\+O(x\log(x)).
\end{multline}
In a similar way we have
\begin{multline}
\label{termelog2} x^{-d}{\rm Tr}_\Sigma\Big({\rm cl}(\nu)\pl_\lam
\til{E}_\infty(\lam;x,y,y') \til{E}_0^\sharp(-\lam;y',x,y)\Big)\\
=\log(x){\rm Tr}_\Sigma\Big({\rm
cl}(\nu)\til{S}_\infty(\lam;y,y')\til{S}^\sharp_0(-\lam;y',y)\Big)
+{\rm Tr}_\Sigma\Big({\rm
cl}(\nu)\pl_\lam\til{S}_\infty(\lam;y,y')\til{S}^\sharp_0(-\lam;y',y)\Big)\\+O(x\log(x)).
\end{multline}
Thus the sum of \eqref{termelog1}
and \eqref{termelog2} integrates in $y,y'$ to a function of $x$ of
the form $\alpha(\lam) \log(x)+\beta(\lam)+O(x\log(x))$ for some
meromorphic function $\alpha(\lam),\beta(\lam)$ which we can
express in terms of the scattering operators. But from the discussion before, we also know that
this trace has no $\log(x)$ coefficients and so $\alpha(\lam)=0$, which
ends the proof by letting $x\to 0$ and writing $\beta(\lam)$ in
terms of
$\til{S}(\lam),\til{S}_\infty(\lam),\til{S}^\sharp_\infty(-\lam)$
and $\til{S}^\sharp_0(-\lam)$ from \eqref{termelog1},
\eqref{termelog2}.
\end{proof}

Let us define the super trace of a trace class operator $A$ on
$L^2(\pl\bar{X},\Sigma)$ by
\[\stra(A):=\frac{1}{i}\tra({\rm cl}(\nu) A).\]

\begin{corollary}\label{lastformula}
Let $\lam\in \cc\setminus(\zz/2)$ be such that $S(z)$ and $S(-z)$
are analytic at $z=\lam$, then the super trace $\stra(\pl_\lam
\til{S}(\lam)\til{S}(-\mu))$ extends meromorphically in $\mu$ from
$\Re(\lam-\mu)<-d$  to
$\mu\in\cc$, it is analytic in $\mu=\lam$, and the following
identity holds
\[\lim_{\eps\to 0}\int_{x(m)>\eps}{\rm Tr}_\Sigma(D\Pi(\lam;m,m)){\rm dv}_{g}(m)=-\frac{i}{2}
\stra(\pl_\lam \til{S}(\lam)\til{S}(-\mu))|_{\mu=\lam}\]
\end{corollary}
\begin{proof} Since
$\pl_\lam\til{S}(\lam)\til{S}^\sharp_{\infty}(-\mu)$ and
$\pl_{\lam}\til{S}_\infty(\lam)\til{S}(-\mu)$ have smooth kernels,
it is clear that their super-trace extends meromorphically to
$\cc$ and is analytic at $\mu=\lam$ by assumption on $\lam$. Now
if we show
\begin{equation}\label{stra=0}
\stra(\pl_\lam \til{S}_0(\lam)\til{S}^\sharp_0(-\mu))=0
\end{equation}
then we have proved the corollary in view of Propositions
\ref{firstred} and \ref{secondred}. We have to study terms of the
form
\begin{equation}\label{gammajk}
\tra\Big(\iota_j^*{\rm cl}(\nu)\pl_\lam
[\eta_j^{-\lam-\ddemi}\phi^2_j\til{S}_{\hh^{d+1}}(\lam)\phi_j^1\eta_j^{-\lam-\ddemi}]
\gamma_{jk}^*\eta_k^{\mu-\ddemi}\phi^1_k\til{S}_{\hh^{d+1}}(-\mu)\phi^2_k\eta_k^{\mu-\ddemi}{\iota_k}_*
\Big)
\end{equation}
where $\gamma_{jk}$ is the unique isometry of $\hh^{d+1}$
extending $\iota_k\circ\iota_j^{-1}:\iota_j(U_k\cap U_j)\to
\iota_k(U_k\cap U_j)$, which acts also as a conformal
transformation on $\pl\bar{B}\subset \rr^{d}$.
As above we use the fact that
$\gamma_{jk}^*R_{\hh^{d+1}}(-\mu)=R_{\hh^{d+1}}(-\mu)\gamma_{jk}^*$
since $\gamma_{jk}$ is an isometry so if
$\alpha_{jk}:=[\gamma_{jk}^*(x_0)/x_0]|_{x=0}\in
C^\infty(\rr^{d})$, then one deduces that
$\gamma_{jk}^*S_{\hh^{d+1}}(-\mu)=\alpha_{jk}^{\mu-\ddemi}S_{\hh^{d+1}}(-\mu)\alpha_{jk}^{\mu-\ddemi}\gamma_{jk}^*$.
So the term \eqref{gammajk} is equal to
\begin{equation}\label{show0}
\begin{gathered}
\int_{\pl\bar{B}\x \pl\bar{B}}\tra_{\Sigma}\Big({\rm cl}(\nu)\pl_\lam [\eta_j(y)^{-\lam-\ddemi}\phi^2_j(y)\til{S}_{\hh^{d+1}}(\lam;y,y')\phi_j^1(y')\eta_j(y')^{-\lam-\ddemi}]\hspace{1.5cm}\\
\hspace{1.5cm}
\eta_{j}(y')^{\mu-\ddemi}\phi^1_k(\gamma_{jk}(y'))\til{S}_{\hh^{d+1}}(-\mu;y',y)
\phi^2_k(\gamma_{jk}(y))\eta_{j}(y)^{\mu-\ddemi}\Big)
(\eta_j(y)\eta_j(y'))^{d}dydy'.
\end{gathered}
\end{equation}
Using the explicit formula of $S_{\hh^{d+1}}(\lam;y,y')$, we see that the
local trace of the operator above, i.e. the integrand in
\eqref{show0}, can be written under the form
\[f_{jk}(\lam,\mu;y,y'){\rm Tr}_{\Sigma}({\rm cl}(\nu)U(y,y')U(y',y))\]
for some function $f_{jk}$ and where $U(y,y')$ is the parallel
transport map on spinors on $\hh^{d+1}$ studied in Section
\ref{paralleltr} and extended down to the boundary $\rr^{d}$. Thus
since $U(y,y')U(y',y)={\rm Id}$ and ${\rm Tr}_\Sigma(\rm
cl(\nu))=0$, we obtain that \eqref{show0} vanishes, which finishes
the proof.
\end{proof}

\subsection{Analysis of residues of $\stra(\pl_\lam\til{S}(\lam)\til{S}(-\lam))$}

Let us define $F(\lam)$ for the value at $\mu=\lam$ of the
meromorphic extension in $\mu$ of
$\stra(\pl_\lam\til{S}(\lam)\til{S}(-\mu))$
\begin{equation}\label{defF}
F(\lam):= \stra(\pl_\lam\til{S}(\lam)\til{S}(-\mu))|_{\mu=\lam}.
\end{equation}
It is clear from Corollary \ref{lastformula} that $F(\lam)$
is meromorphic in $\lam\in\cc$, but we want to prove that it has
only first order poles, the residue of which are integers.  Since
$\til{S}(\lam)$ is unitary on $\{\Re(\lam)=0,\lam\not=0\}$ it is
analytic at $\lam=0$, so one can define
\begin{equation}\label{defmcS}
\mc{S}_\pm(\lam):=\til{S}_\pm(\lam)\til{S}_\mp(0):
C^\infty(\pl\bar{X};{^0\Sigma}_\mp)\to
C^\infty(\pl\bar{X};{^0\Sigma}_\mp)
\end{equation}
which are the two diagonal components of
$\mc{S}(\lam):=\til{S}(\lam)\til{S}(0)$ in the splitting
${^0\Sigma}_+\oplus{^0\Sigma}_-$. These two operators are elliptic
pseudo-differential operators of complex order $2\lam$ by
Proposition \ref{spseudo}
\[\mc{S}_\pm(\lam)\in \Psi^{2\lam}(\pl\bar{X};{\rm End}({^0\Sigma}_\mp)),\]
and their principal symbol is $|\xi|^{2\lam}$. Let $D$ be the
Dirac operator on $(\pl\bar{X},h_0)$ and let $D_\mp=P_\pm
DI_\mp:C^{\infty}(\pl\bar{X};{^0\Sigma}_\mp)\to
C^{\infty}(\pl\bar{X};{^0\Sigma}_\pm)$ be the off diagonal
components of $D$. If $|D|_\mp:=(D_\pm D_\mp)^\demi$,  it is possible
to factorize $\mc{S}_\pm(\lam)$ by
$(1+|D|_\mp)^{-\lam}\mc{S}_\pm(\lam)(1+|D|_\mp)^{-\lam}$ and this
operator is of the form $1+K(\lam)$ for some meromorphic family of
compact operators on $L^2(\pl\bar{X};{^0\Sigma}_\mp)$, it is thus
Fredholm on this space. Then we can use the theory of
Gohberg-Sigal \cite{GSi} like in \cite{GZannals} or Section 2 of
\cite{GuiMRL} for these operators. In particular, one can define
the null multiplicities $N_{\lam_0}(\mc{S}_\pm(\lam))$ of
$\mc{S}_\pm(\lam)$ at a point $\lam_0$ as follows: by the theory
of \cite{GSi}, for a meromorphic family of operators
$L(\lam)=1+K(\lam)$ acting on a Hilbert space $\mc{H}$ with
$K(\lam)$ compact, with $L(\lam)$ invertible for some $\lam$, then
there exist holomorphically invertible operators
$U_1(\lam),U_2(\lam)$ near $\lam_0$, some $(k_l)_{l=0,\dots,m}\in
\zz\setminus\{0\}$ (with $m\in\nn$) and some orthogonal projectors
$P_l$ on $L^2(\pl\bar{X},{^0\Sigma}_\pm)$ such that ${\rm
rank}(P_l)=1$ if $l>0$, $P_iP_j=\delta_{ij}$ and
\begin{equation}\label{factor}
L(\lam)=U_1(\lam)\Big(
P_0+\sum_{l=1}^m(\lam-\lam_0)^{k_l}P_l\Big)U_2(\lam),
\end{equation}
then we define the null multiplicity at $\lam_0$ by
\begin{equation}\label{defNM}
N_{\lam_0}(L(\lam)):= \sum_{k_l>0}k_l.
\end{equation}
Note that, by \cite{GSi}, this is an integer depending only on
$\mc{S}_\pm(\lam)$
and not on the factorization \eqref{factor} and that
$N_{\lam_0}(L(\lam))=0$ if and only if $L(\lam)^{-1}$ is
holomorphic at $\lam=\lam_0$.

\begin{proposition}\label{residuesF}
The function $F(\lam)$ of \eqref{defF} is meromorphic in
$\lam\in\cc$, one has
\begin{equation}\label{firstident}
F(\lam)={\rm FP}_{\mu=\lam}{\rm
TR}(\pl_\lam\mc{S}_-(\lam)\mc{S}_-^{-1}(\mu))-{\rm
FP}_{\mu=\lam}{\rm TR}(\pl_\lam\mc{S}_+(\lam)\mc{S}_+^{-1}(\mu))
\end{equation}
where ${\rm TR}$ is the Kontsevich-Vishik trace of \cite{KV} and ${\rm
FP}_{\mu=\lam}$ means the finite part (or regular value) of the
meromorphic function of $\mu$ at $\mu=\lam$. The poles of
$F(\lam)$ are first order poles, with residue at a pole $\lam_0$
given by
\begin{equation*}\begin{split}
{\rm Res}_{\lam=\lam_0}F(\lam)=&\ {\rm Tr}({\rm Res}_{\lam=\lam_0}(\pl_\lam\mc{S}_-(\lam)\mc{S}_-(\lam)^{-1}))-{\rm Tr}({\rm Res}_{\lam=\lam_0}(\pl_\lam\mc{S}_+(\lam)\mc{S}_+(\lam)^{-1}))\\
=&\
(N_{\lam_0}(\mc{S}_-(\lam))-N_{\lam_0}(\mc{S}_-(\lam)^{-1}))-(N_{\lam_0}(\mc{S}_+(\lam))-N_{\lam_0}(\mc{S}_+(\lam)^{-1}))
\end{split}
\end{equation*}
where $N_{\lam_0}$ is the null multiplicity defined in
\eqref{defNM}.
\end{proposition}
\begin{proof}
The first statement is straightforward since
\[\stra(\pl_\lam\til{S}(\lam)\til{S}(-\mu))={\rm TR}(\pl_\lam\mc{S}_-(\lam)\mc{S}_-^{-1}(\mu))-{\rm TR}(\pl_\lam\mc{S}_+(\lam)\mc{S}_+^{-1}(\mu))\]
and we know from the work of Lesch \cite{L} that the
Kontsevich-Vishik trace of an analytic family of
log-polyhomogeneous operators $A(\mu)$ extend meromorphically to
$\mu\in\cc$, so it suffices to use the fact that
$\stra(\pl_\lam\til{S}(\lam)\til{S}(-\mu))$ analytically continues
to $\mu\in\cc$ and is analytic at $\mu=\lam$ to prove \eqref{firstident}.

As shown in Proposition \ref{secondred}, $F(\lam)$ can be written
as a trace of a meromorphic family of trace class operators,  more
precisely, using the fact that ${\rm cl}(\nu)$ anti-commutes
with $\til{S}(\lam)$ and $\til{S}_0(\lam)$ for all $\lam$ where
they are defined,
\[F(\lam)=\frac{1}{i}\tra\Big({\rm cl}(\nu)(\pl_\lam\til{S}(\lam)\til{S}(\lam)^{-1}-\pl_\lam\til{S}_0(\lam)\til{S}^\sharp_0(-\lam))\Big).\]
Consequently, the polar part of $F(\lam)$ at a pole $\lam_0$ is
given by the trace of the polar part (which is finite rank) of
${\rm cl}(\nu)(\pl_\lam\til{S}(\lam)\til{S}(\lam)^{-1}-\pl_\lam\til{S}_0(\lam)\til{S}^\sharp_0(-\lam))$.
But clearly from the explicit formula of $S_{\hh^{d+1}}(\lam)$, we see
that $\pl_\lam\til{S}_0(\lam)\til{S}^\sharp_0(-\lam)$ is
holomorphic in $\lam\in\cc$. Now use that $\til{S}(0)^2=\rm{Id}$
to write
\begin{equation}
\begin{split}
{\rm cl}(\nu) \pl_\lam\til{S}(\lam)\til{S}(\lam)^{-1}=& + iI_+\pl_\lam(\til{S}_-(\lam)\til{S}_+(0))(\til{S}_-(\lam)\til{S}_+(0))^{-1}P_+\\
&-iI_-\pl_\lam(\til{S}_+(\lam)\til{S}_-(0))(\til{S}_+(\lam)\til{S}_-(0))^{-1}P_-
\end{split}\end{equation}
and write a factorization of the form \eqref{factor} for
$\mc{S}_\pm(\lam)$, from which it is clear that
$\pl_\lam\mc{S}_\pm(\lam)\mc{S}_\pm(\lam)^{-1}$ has only first
order poles except possibly
\[\begin{gathered}
 U_1(\lam)(P_0+\sum_{l=1}^m(\lam-\lam_0)^{k_l}P_l)\pl_\lam U_2(\lam)U_2(\lam)^{-1}(\sum_{l=1}^m(\lam-\lam_0)^{-k_l}P_l)U_1(\lam)^{-1}\\
+ U_1(\lam)(\sum_{l=1}^m(\lam-\lam_0)^{k_l}P_l)\pl_\lam
U_2(\lam)U_2(\lam)^{-1}(P_0+\sum_{l=1}^m(\lam-\lam_0)^{-k_l}P_l)U_1(\lam)^{-1}
\end{gathered}.\]
But this is a finite rank operator and so by the cyclicity of the
trace, we deduce that the trace of this term is holomorphic in
$\lam$. To finish the proof, it suffices to apply the main result
of \cite{GSi}:
\[\tra({\rm Res}_{\lam=\lam_0}(\pl_\lam \mc{S}_\pm(\lam)\mc{S}_\pm(\lam)^{-1}))=
N_{\lam_0}(\mc{S}_\pm(\lam))-N_{\lam_0}(\mc{S}_\pm(\lam)^{-1}).\]
\end{proof}

Let us define the multiplicity of resonances as follows
\begin{equation}\label{defmultip}
m_\pm (\lam_0):={\rm rank}({\rm Res}_{\lam=\lam_0}R_\pm(\lam))
\end{equation}
We want to identify scattering poles and resonances.
\begin{proposition}\label{poleSpoleR}
Let $\lam_0\in\cc$, then the following identity holds
\begin{equation}\label{relationmult}
N_{\lam_0}(\mc{S}_\mp(-\lam))=m_\pm(\lam_0)+\indic_{-1/2-\nn_0}(\lam_0)\dim\ker
\mc{S}_\mp(-\lam_0).
\end{equation}
\end{proposition}
\begin{proof}
We just sketch the proof since it is very similar to that of
Theorem 1.1. of \cite{GuiMRL}, and we strongly encourage the
reader to look at \cite{GuiMRL}. The first thing to notice is that
$N_{\lam_0}(\mc{S}_\pm(-\lam))=N_{-\lam_0}(\mc{S}_\pm(\lam))$ and
that
$N_{\lam_0}(\mc{S}_\pm(\lam)^{-1})=N_{\lam_0}(\mc{S}_\mp(-\lam))$
since $\mc{S}_\pm(\lam)^{-1}=
\til{S}_{\pm}(0)\mc{S}_{\mp}(-\lam)\til{S}_\pm(0)^{-1}$. Remark
that $R_\pm(\lam)$ and $\mc{S}_\pm(\lam)$ are analytic in
$\{\Re(\lam)\geq 0\}$ and so the identity \eqref{relationmult}
is trivial (all terms are $0$) for $\Re(\lam_0)\geq 0$.

Now suppose that $\Re(\lam_0)<0$. First we prove that
\begin{equation}\label{sens1}
N_{\lam_0}(\mc{S}_\mp(-\lam))-\indic_{-1/2-\nn_0}(\lam_0)\dim\ker
\mc{S}_\mp(-\lam_0) \leq m_\pm(\lam_0).
\end{equation}
By \eqref{kernelEla} and \eqref{kernelSla}, $S_\pm(\lam)$ can be
represented for $\Re(\lam)<-\ddemi$ by
\[S_\pm(\lam;y,y')=\pm i[(xx')^{-\lam-\ddemi}R_\pm(\lam;x,y,x',y')]|_{x=x'=0}\]
and the expression can be
extended to $\lam\in\cc$ meromorphically as a singular integral
kernel using the blow-down maps like in \eqref{descriptionS}. Then
we can apply mutatis mutandis Lemma 3.2 of \cite{GuiMRL}, where
$S(\lam)$ there is replaced by $S_\pm(\lam)$ here, the function
$z(\lam)$ there is $\lam$ here, and we have to multiply the
factorization (3.11) of \cite{GuiMRL} by $\til{S}_\mp(0)$ on the
right, which is harmless since it does not depend on $\lam$. We
want to apply the factorization of $\mc{S}(\lam)$ obtained from
this Lemma 3.2 of \cite{GuiMRL} to prove \eqref{sens1}, in a way
similar to Corollary 3.3 of \cite{GuiMRL}. First the Corollary 3.3
in \cite{GuiMRL} can also be rewritten (using the notations of
\cite{GuiMRL}) under the form
\[N_{\lam_0}(\til{S}(n-\lam))-\indic_{-n/2-\nn_0}(\lam_0)\dim\ker\til{S}(n-\lam_0)\leq m_{\lam_0}(z'(\lam)R(\lam))\]
by using equation (3.19) in \cite{GuiMRL} if $\lam_0\in n/2-\nn$
and the fact that $c(n-\lam)$ is holomorphic at all $\lam_0\notin
n/2-\nn$. Then the proof of this Corollary 3.3 in \cite{GuiMRL}
can be copied word by word by replacing $\til{S}(n-\lam)$ and
$c(n-\lam)$ there by $\mc{S}_\mp(-\lam)$ and $C(-\lam)$ here, and
$m_{\lam_0}(z'(\lam)R(\lam))$ by $m_\pm(\lam_0)$. This finally
proves \eqref{sens1}.

Then we need to prove the converse inequality of \eqref{sens1}.
From Lemma \ref{functeq1}, Corollary
\ref{SvsE} and the fact that $\til{S}_\mp(0)\til{S}_\pm(0)={\rm
Id}|_{{^0\Sigma}_\mp}$, we deduce
\begin{equation}\label{eqfunct2}
R_\pm(\lam)-R_\mp(-\lam)=-E_{\mp}(-\lam)C(\lam)\mc{S}_\pm(\lam)\til{S}_{\pm}(0){\rm
cl}(\nu)E^{\sharp}_{\mp}(-\lam)
\end{equation}
which is the equivalent in our setting to the identity (3.15) of
\cite{GuiMRL}. Since the Lemma 3.4 of \cite{GuiMRL} is only based
on the identity (3.15) in \cite{GuiMRL}, the structure of the
resolvent kernel at the boundary and the unique continuation
principle of Mazzeo \cite{MaAJM}, the same proof applies and is
actually easier in our case since there is no pure point spectrum
thus no resonance in the physical sheet $\{\Re(\lam)>0\}$. This
implies
\[N_{\lam_0}(\mc{S}_\mp(-\lam))-\indic_{-1/2-\nn_0}(\lam_0)\dim\ker \mc{S}_\mp(-\lam_0)
\geq m_\pm(\lam_0).\] The idea of the proof of Lemma 3.4 of
\cite{GuiMRL} is to use \eqref{eqfunct2} to write the residue of
$R_\pm(\lam)$ at $\lam_0$ with $\Re(\lam_0)<0$ in terms of the
singular part of the Laurent expansion of $\mc{S}_\pm(\lam)$,
itself obtained from a factorization of the form \eqref{factor},
then use the fact that $R_\mp(-\lam)$, $E_\mp(-\lam)$ and
$E_\mp^\sharp(-\lam)$ are holomorphic in $\{\Re(\lam)<0\}$ and finally
count the rank of the residue in terms of the $k_l$ of the
factorization \eqref{factor}.
\end{proof}

\begin{theorem}
The function $F(\lam)$ is meromorphic with first order poles and
integer residues given by
\begin{equation*}\begin{split}
{\rm Res}_{\lam=\lam_0}F(\lam)=&\,m_+(\lam_0)-m_-(\lam_0)+\indic_{-1/2-\nn_0}(\lam_0){\rm Ind} (\mc{S}_-(-\lam_0))\\
=&\, m_+(\lam_0)-m_-(\lam_0).
\end{split}\end{equation*} for $\Re(\lam_0)\leq 0$, where
$m_\pm(\lam_0)$ is defined in \eqref{defmultip}.
\end{theorem}
\begin{proof}
Apply Proposition \ref{poleSpoleR} with Proposition
\ref{residuesF}. To see the index of $\mc{S}_{-}(-\lam_0)$
appearing, we also use that
$\mc{S}_+(-\lam_0)^*=\mc{S}_-(-\lam_0)$ for $\lam_0\in\rr$, which
come from the self-adjointness of $\til{S}(-\lam_0)$. The fact
that the index of $\mc{S}_{-}(-\lam_0)$ vanishes comes from the
invariance of the index by continuous deformation and the
invertibility of $\til{S}(\lam)$ except on a discrete set of
$\lam\in\cc$.
\end{proof}

We deduce directly our main theorem from this Theorem, Corollary \ref{lastformula} and the identity \eqref{zetavsresolvent}:
\begin{theorem}\label{mainth}
The odd Selberg zeta function $Z^o_{\Gamma,\Sigma}(\lam)$ on a
convex co-compact hyperbolic manifold $X_\Gamma$ of dimension
$2n+1$ has a meromorphic extension to $\cc$, is analytic in
$\{\Re(\lam)\geq 0\}$, and $\lam_0$ is a zero or pole if and only
if the meromorphic extension $R_+(\lam)$ or $R_-(\lam)$ of $(D\pm
i\lam)^{-1}$ from $\{\Re(\lam)>0\}$ to $\cc$ have a pole at
$\lam_0$, in which case the order of $\lam_0$ as a zero or pole of
$Z_{\Gamma,\Sigma}^o(\lam)$ (with the positive sign convention for
zeros) is given by
\[{\rm Rank }\, {\rm Res}_{\lam_0}R_-(\lam)-{\rm Rank }\,{\rm
Res}_{\lam_0}R_+(\lam).\]
\end{theorem}

\section{Eta invariant of odd signature operator and its structure
on Schottky space}

For  a $(4m-1)$-dimensional convex co-compact hyperbolic manifold
$X_\Gamma$, we consider the odd signature operator $A$ on odd
forms $\Lambda^{\mathrm{odd}}=\oplus_{p=0}^{2m}\Lambda^{2p-1}$
acting by $(-1)^{m+p}(\star d+d\star)$ over $\Lambda^{2p-1}$ as in
Millson's paper \cite{M}. Recall that $A^2=\Delta$ and $A$ is self
adjoint. We want to make a sense of
\begin{equation}\label{defeta}
\eta(A):=\frac{1}{{\sqrt{\pi}}}\int_0^\infty
t^{-\demi}\int_{X_\Gamma}{\rm
Tr}_{\Lambda}(Ae^{-t\Delta})(m)\,{\rm dv}(m)\, dt
\end{equation}
where ${\rm Tr}_{\Lambda}$ is the local trace on the bundle
$\Lambda^{\mathrm{odd}}$. First it is easy to see that ${\rm
Tr}_\Lambda(A e^{-t\Delta})$ is the same as ${\rm
Tr}_{\Lambda^{2m-1}}( \star d\, e^{-t\Delta})$ since the other
parts are off diagonal if we write $A e^{-t\Delta}$ as a matrix
with respect to the natural basis of $\Lambda^{\mathrm{odd}}$.
First we show that the local trace
$\mathrm{Tr}_{\Lambda^{2m-1}}(\star d\, e^{-t\Delta})$ is
integrable on $X_\Gamma$. As in the spinor bundle $\Sigma$, the
bundle of $(2m-1)$-forms can be understood as a homogeneous vector
bundle given by the representation $\Lambda^{2m-1}\phi$ with the
standard representation $\phi$ of $\mathrm{SO}(4m-1)$, which
decomposes into
\begin{equation}\label{e:decomp-lambda}
\Lambda^{2m-1}\phi\,
|_{M=\mathrm{SO}(4m-2)}=\Lambda^{2m-1}_+\bar{\phi}\oplus
\Lambda^{2m-1}_-\bar{\phi}\oplus \Lambda^{2m-2}\bar{\phi}
\end{equation}
where $\bar{\phi}$ denotes the standard representation of
$M=\mathrm{SO}(4m-2)$. As in the subsection \ref{ss-har}, there is
the $\Lambda^{2m-1}\phi$-radial function $P_t$ associated to
$\star d\, e^{-t\Delta}$, and we have the corresponding scalar
functions $p_t^\pm(r)$, $p_t^{2m-2}(r)$ of $P_t$ restricting to
the representation spaces on the right hand side of
\eqref{e:decomp-lambda}. Now, as in Proposition \ref{p:scalar}, we
have

\begin{proposition}\label{p:scalar2} The scalar components
$p_t^\pm(r)$, $p_t^{2m-2}(r)$ are given by
\begin{equation*}
p_t^\pm(r)=\pm\frac{(4m-1)\sinh(r)}{ i 2^{2m-1/2}
\pi^{2m+1/2}t^{3/2}} \left(-\frac{d}{d(\cosh r)}\right)^{2m-1} r
\sinh^{-1} (r)\, e^{-\frac{r^2}{4t}}, \qquad p_t^{2m-2}(r)\equiv
0.
\end{equation*}
\end{proposition}
\begin{proof} The equalities follow from Lemma 7.4 and Theorem 7.6 in
\cite{Pe} and Theorem 1.1 in \cite{M}.
\end{proof}

Using this proposition and repeating the same argument as in
Section \ref{geometricside}, one can easily show that
$\mathrm{Tr}_{\Lambda^{2m-1}}(\star d\, e^{-t\Delta})$ is
integrable over $X_\Gamma$. By the same argument as in Proposition
\ref{p:odd-heat} and Corollary \ref{c:asymp}, one can also obtain
the corresponding results, which  implies that the integral
$\int^\infty_0 \cdot \, dt$ in \eqref{defeta} converges. Hence the
eta invariant $\eta(A)$ given in \eqref{defeta} is well defined.

For $\Re(\lambda)>\delta_\Gamma$, we also have the Selberg
zeta function of odd type $Z^o_{\Gamma,\Lambda}(\lambda)$  just
putting  $\sigma_\pm=\Lambda^{2m-1}_\pm\bar{\phi}$ in
\eqref{e:def-Sel-halh} and \eqref{e:def-Sel-odd}, which coincides
with the one in \eqref{defoddzeta} introduced by Millson
\cite{M}. We first have a result similar to the case of spinor
bundle dealt with above:
\begin{theorem}
If $X_\Gamma:=\Gamma\backslash \hh^{4m-1}$ is a convex co-compact
hyperbolic manifold, then the local trace ${\rm
Tr}_{\Lambda^{2m-1}}(\star d\, e^{-t\Delta})$ is integrable on $X$
for all $t>0$, so that the integral \eqref{defeta} converges and
defines the eta invariant $\eta(A)$. Moreover, if the Poincar\'e
exponent $\delta_\Gamma<0$,
\begin{equation}\label{e:eta-odd-sign} e^{\pi
i\eta(A)}=Z^o_{\Gamma,\Lambda}(0).
\end{equation}
\end{theorem}
\begin{proof}We already showed the first claim. The equality
\eqref{e:eta-odd-sign} also easily follows from the corresponding
results to Proposition \ref{p:odd-heat} and the assumption
$\delta_\Gamma<0$, from which  we do not need to show the
meromorphic extension of $Z^o_{\Gamma,\Lambda}(\lambda)$.
\end{proof}

It turns out that this eta invariant $\eta(A)$ of the odd
signature operator $A$ has an intimate relationship with the
deformation space of the hyperbolic structures when $X_\Gamma$ is
$3$-dimensional. To explain this, first we review the work of
Zograf \cite{Zog}.

\subsection{Zograf factorization formula}
A \emph{marked Schottky group} is a discrete subgroup $\Gamma$ of
the linear fractional transformations
$\mathrm{PSL}(2,\mathbb{C})$, with distinguished free generators
$\gamma_1,\gamma_2,\ldots, \gamma_g$ satisfying the following
condition: there exist $2g$ smooth Jordan curves $C_r$,
$r=\pm1,\ldots,\pm g$, which form the oriented boundary of a
domain $\Omega_0$ in $\hat{\mathbb{C}}=\mathbb{C}\cup\{\infty\}$
such that $\gamma_r C_r = -C_{-r}$ for $r=1,\ldots,g$. If $\Omega$
is the union of images of $\Omega_0$ under $\Gamma$, then
$Y_\Gamma=\Gamma\backslash \Omega$ is a compact Riemann surface of
genus $g$. The action of $\Gamma$ on ${\mathbb{C}}$ naturally
extends to the action on $\mathbb{H}^3$ where
$\partial\mathbb{H}^3=\mathbb{C}$ and the quotient space
$X_\Gamma=\Gamma\backslash \mathbb{H}^3$ is a \emph{Schottky
hyperbolic $3$-manifold} whose boundary is the Riemann surface
$Y_\Gamma$. Here let us remark that $\delta_\Gamma+1$ is the
Hausdorff dimension of the limit set $\Lambda$ in $\partial
\hh^{3}$ of $\Gamma$ and $\delta_\Gamma$ is also the smallest
number such that $\prod_{\{\gamma\}}(1-q_\gamma^{s})$ absolutely
converges whenever $s>\delta_\Gamma+1$. The function
$\prod_{\{\gamma\}}(1-q_\gamma^{s})$ was briefly described in
\cite{Bow} where it was asserted without the proof that with the
values of $q_\gamma^s$ chosen appropriately, the infinite product
is defined for $\Re(s)>\delta_\Gamma+1$ and has an analytic
continuation to $\mathbb{C}$.

Each nontrivial element $\gamma\in\Gamma$ is loxodromic: there
exists a unique number $q_\gamma\in\mathbb{C}$ (the multiplier)
such that $0< |q_\gamma| < 1$ and $\gamma$ is conjugate in
$\mathrm{PSL}(2,\mathbb{C})$ to $z\mapsto q_\gamma z$, that is,
\[
\frac{\gamma z- a_\gamma}{\gamma z- b_\gamma} = q_\gamma
\frac{z-a_\gamma}{z-b_\gamma}
\]
for some $a_\gamma, b_\gamma\in\hat{\mathbb{C}}$ (the attracting
and repelling fixed points respectively). A marked Schottky group
with an ordered set of free generators $\gamma_1,\ldots,\gamma_g$
is \emph{normalized} if $a_{\gamma_1}=0$, $b_{\gamma_1}=\infty$,
$a_{\gamma_2}=1$. The \emph{Schottky space} $\mathfrak{S}_g$ is
the space of marked normalized Schottky groups with $g$
generators. It is a complex manifold of dimension $3g-3$, covering
the Riemann moduli space $\mathfrak{M}_g$ and with universal cover
the Teichm\"uller space $\mathfrak{T}_g$.

Like the Teichm\"uller space $\mathfrak{T}_g$, the Schottky space
$\mathfrak{S}_g$ has a natural K\"ahler metric, the Weil-Petersson
metric. In \cite{ZT}, Takhtajan-Zograf constructed a K\"ahler
potential $S$ called \emph{classical Liouville action} of the
Weil-Petersson metric on $\mathfrak{S}_g$, that is,
\begin{equation}\label{potential1}
\partial \overline{\partial} S= 2i\, \omega_{WP}
\end{equation}
where $\partial$ and $\overline{\partial}$ are the $(1,0)$ and
$(0,1)$ components of the de Rham differential $d$ on
$\mathfrak{S}_g$ respectively, and $\omega_{WP}$ is the symplectic
form of the Weil-Petersson metric. On the other hand, from the
local index theorem for families of
$\overline{\partial}$-operators in Takhtajan-Zograf \cite{ZT2},
the following equality also follows
\begin{equation}\label{potential2}
\partial\overline{\partial} \log
\frac{\mathrm{Det}\Delta}{\mathrm{det}\, \mathrm{Im}\, \tau} =
-\frac{i}{6\pi} \omega_{WP}
\end{equation}
where $\mathrm{Det}\Delta$ denotes the $\zeta$-regularized
determinant of the Laplacian $\Delta$  given by the hyperbolic
metric and $\tau$ denotes the period matrix. Comparing
\eqref{potential1} with \eqref{potential2}, one can expect a
nontrivial relationship between $S$ and $ \log
\frac{\mathrm{Det}\Delta}{\mathrm{det}\, \mathrm{Im}\, \tau}$.
Indeed, in \cite{Zog}, \cite{Zog2} Zograf proved

\begin{theorem} {\bf{[Zograf]}}
There exists a holomorphic function $F(\Gamma):\mathfrak{S}_g \to
\mathbb{C}$ such that
\begin{equation}\label{e:Zograf}
\frac{\mathrm{Det}\Delta}{\mathrm{det}\, \mathrm{Im}\, \tau} = c_g
\exp\left( -\frac{1}{12\pi} S\right) |F(\Gamma)|^2
\end{equation}
where $c_g$ is a constant depending only on $g$. For points in
$\mathfrak{S}_g$ corresponding to Schottky groups $\Gamma$ with
$\delta_\Gamma<0$, the function $F(\Gamma)$ is given by the
following absolutely convergent product:
\begin{equation}\label{e:Zog-F}
F(\Gamma)= \prod_{\{\gamma\}}\prod^\infty_{m=0} (1-
q_\gamma^{1+m})
\end{equation}
where $q_\gamma$ is the multiplier of $\gamma\in\Gamma$, and
$\{\gamma\}$ runs over all distinct primitive conjugacy classes in
$\Gamma$ excluding the identity.
\end{theorem}

Combining the equalities \eqref{e:Zograf} and \eqref{e:Zog-F},
these are called \emph{Zograf factorization formula}. This result
was extended by McIntyre-Takhtajan to the Schottky groups without
the condition for $\delta_\Gamma$ in \cite{Mc-Tak}. Here they used
the $\zeta$-regularized determinant of $\Delta_n$ acting on the
space of $n$-differentials so that the corresponding holomorphic
function is $F_n(\Gamma)= \prod_{\{\gamma\}}\prod^\infty_{m=0} (1-
q_\gamma^{n+m})$ which absolutely converges for any Schottky group
$\Gamma$ if $n>1$.

\subsection{Eta invariant as a functional over the Schottky space}
By the construction of $X_\Gamma$ and its boundary $Y_\Gamma$, the
eta invariant $\eta(A)$ can be understood as a functional over the
Schottky space $\mathfrak{S}_g$. Now a natural question is to
describe the eta invariant $\eta(A)$ as a functional over
$\mathfrak{S}_g$. For this, we have

\begin{theorem}\label{t:eta-F} Let $\mathfrak{S}_g^0$ to be a subset of $\mathfrak{S}_g$
consisting of normalized Schottky groups $\Gamma$'s  with the
property $\delta_\Gamma<0$. Then we have
\[
F(\Gamma)=|F(\Gamma)| \exp\left(-\frac{\pi i}{2}
\eta(A)\right)\qquad \text{over\ \ $\mathfrak{S}_g^0$},
\] in particular, $\eta(A)$ is a pluriharmonic function over
$\mathfrak{S}_g^0$.
\end{theorem}

\begin{proof} The proof is a simple application of the equality
\eqref{e:eta-odd-sign}. For this, as in Proposition
\ref{p:zeta-exp}, we rewrite $Z^o_{\Gamma,\Lambda}(\lambda)$ with
respect to the group $\mathrm{PSL}(2,\mathbb{C})$ as follows:
\begin{equation}\label{e:def-Sel-sign}
Z^o_{\Gamma,\Lambda}(\lambda)=\prod_{[\gamma]\in
\mathrm{P}\Gamma_{\mathrm{lox}}} \prod_{k,\ell=0}^\infty\ \frac{
\left(1- e^{i\theta_\gamma}\, (\mu_\gamma)^{-2k}\,
(\bar{\mu}_\gamma)^{-2\ell} \,
|\mu_\gamma|^{-2(\lambda+1)}\right)}{
\left(1-e^{-i\theta_\gamma}\, (\mu_\gamma)^{-2k}\,
(\bar{\mu}_\gamma)^{-2\ell}\,
|\mu_\gamma|^{-2(\lambda+1)}\right)}.
\end{equation}
Here $\gamma$ runs over the set of $\Gamma$-conjugacy classes of
the loxodromic elements in $\Gamma$ and a loxodromic element
$\gamma$ can be conjugated to a diagonal matrix with the diagonal
elements $\mu_\gamma=\exp(\frac12({l_\gamma+i\theta_\gamma}))$,
$\mu_\gamma^{-1}=\exp(-\frac12(l_\gamma+i\theta_\gamma))$ in
$\mathrm{PSL}(2,\mathbb{C})$ ($|\mu_\gamma|>1$). Let us remark
that the infinite product on the right hand side of
\eqref{e:def-Sel-sign} absolutely converges for
$\Re(\lambda)>\delta_\Gamma$, in particular, at $\lambda=0$ since
$\delta_\Gamma<0$.

Now comparing the definition of $q_\gamma$ and $\mu_\gamma$, one
can see that $q_\gamma=\mu_\gamma^{-2}$, that is,
$q_\gamma^{1/2}=\mu_\gamma^{-1}$. Hence the odd Selberg zeta
function $Z^o_{\Gamma,\Lambda}(0)$ has the following expression in
terms of $q_\gamma$,
\begin{align*}
Z^o_{\Gamma,\Lambda}(0)=&\prod_{[\gamma]\in
\mathrm{P}\Gamma_{\mathrm{lox}}} \prod_{k,\ell=0}^\infty\ \frac{
\left(1- (\bar{q}_\gamma q_\gamma^{-1})^{\frac12} \,
q_\gamma^{k}\, \bar{q}_\gamma^{\ell} \, (q_\gamma
\bar{q}_\gamma)^{\frac12}\right)}{ \left(1-({q}_\gamma
\bar{q}_\gamma^{\,-1})^{\frac12} \, q_\gamma^{k}\,
\bar{q}_\gamma^{\ell} \, (q_\gamma
\bar{q}_\gamma)^{\frac12}\right)} \\=& \prod_{[\gamma]\in
\mathrm{P}\Gamma_{\mathrm{lox}}} \prod_{k,\ell=0}^\infty\ \frac{
\left(1-  q_\gamma^{k}\, \bar{q}_\gamma^{\ell+1} \right)}{
\left(1- q_\gamma^{k+1}\, \bar{q}_\gamma^{\ell}\right) }
=\prod_{[\gamma]\in \mathrm{P}\Gamma_{\mathrm{lox}}}
\prod_{m=0}^\infty \frac{\left(1-
\bar{q}_\gamma^{1+m}\right)}{\left(1-q_\gamma^{1+m}\right)}.
\end{align*}
Combining this and \eqref{e:eta-odd-sign} completes the proof.

\end{proof}

\begin{remark*} In the proof of Theorem \ref{t:eta-F}, we assume the
condition $\delta_\Gamma<0$ which simplifies the proof in several
steps. But, one can expect that a similar result still holds over
whole Schottky space $\mathfrak{S}_g$. This extension to
$\mathfrak{S}_g$ is also related to the proof of the assertion of
Bowen in \cite{Bow} about the meromorphic extension of
$\prod_{\{\gamma\}}(1-q_\gamma^{s})$ over $\mathbb{C}$. These
problems will be discussed elsewhere.
\end{remark*}

\appendix
\section{Computation of parallel transport in the spinor bundle}\label{paralleltr}

Let $\tau_m^{m'}$ denote parallel transport in the tangent bundle of the upper
half-space model $\mathbb{H}^{d+1}$ of hyperbolic space, between points
$m=(x,y),m'=(x',y')$, along the unique geodesic linking them. We identify
$T_m \mathbb{H}^{d+1}$ with $\rz^{d+1}$ using the orthonormal basis
at $m$ given by
$\{x\pl_x,x\pl_{y_1},\ldots,x\pl_{y_d}\}$. We denote by $\tau(m,m')$
the matrix of the transformation $\tau_m^{m'}$ written in these bases.

\begin{proposition}\label{partr}
Let $r:=|y-y'|$, $\rho_{\rm ff}:=\sqrt{(x+x')^2+r^2}$.
The special orthogonal matrix $\tau(m,m')$ has the following coefficients:
\begin{align*}
\tau_{00}=&\, 1-2{r^2}/{\rho_{\rm ff}^2}\\
\tau_{0j}=&\, -2(x+x')(y_j-y_j')/\rho_{\rm ff}^2 && \text{for $j=1,\ldots,d$}\\
\tau_{j0}=&\, 2(x+x')(y_j-y_j')/\rho_{\rm ff}^2  && \text{for $j=1,\ldots,d$}\\
\tau_{jl}=&\, \delta_j^l-2(y_j-y_j')(y_l-y_l')/\rho_{\rm ff}^2&&
\text{for $j,l\in\{1,\ldots,d\}$}.
\end{align*}
which are smooth on the stretched product $\hh^{d+1}\x_0\hh^{d+1}$ defined in subsection \ref{stretched}.
\end{proposition}
\begin{proof}
Let $A$ be the translation by $(0,y')$ in $\rz^{d+1}$,
composed to the left by the homothety of factor $1/x'$. This isometry
of $\mathbb{H}^{d+1}$ maps $m'$ to $(1,0)$. In the above trivialization
of the tangent bundle, $A_*$ acts as the identity. Moreover, since it is
an isometry, $A$ transforms the geodesic from $m$ to $m'$ into the geodesic
from $A(m)$ to $(1,0)$ and preserves parallelism. Thus (as matrices)
\begin{align}\label{umrn}
\tau(m,m')=\tau(A(m),(1,0)),&& \text{where $A(m)=\left(\frac{x}{x'},
\frac{y-y'}{x'}\right)$.}
\end{align}
We now concentrate on $\tau(m,(1,0))$.

If $y=0$ it is clear that $\tau((x,0),(1,0))$ is just the identity matrix.
Suppose $y\neq 0$. Let $\pl_r:=\frac{1}{r}\sum_{j=1}^d y_j\pl_{y_j}$ denote
the radial vector field, defined outside the vertical line through the origin.
Define $R:=x\pl_r$, $X:=x\pl_x$. For each $j\geq 1$ set $e_j:=x\pl_{y_j}$,
and let $T_j:=e_j-\langle e_j,R\rangle R$ denote the component of $e_j$
which is tangent to the sphere $S^{d-1}$. The geodesic
from $m$ to $(1,0)$ lives in the totally geodesic plane $\Pi_m$ passing through
$(1,0)$ and $m$, which is a copy of the hyperbolic $2$-space. Along this plane
the vector fields $T_j$ extend smoothly at the vertical line through the origin.
It is clear that the vector fields $T_j$ are parallel along $\Pi_m$.

\begin{lemma} \label{ptp}
In the plane $\Pi_m$, parallel transport between $m$ and
$(1,0)$ is given by the complex number
\[\frac{-r+i(1+x)}{\ r+i(1+x)}\]
\end{lemma}
\begin{proof}
We use as (real) basis for $T\Pi_m$ the orthonormal vector fields $X$ and $R$.
The complex structure rotates $R$ to $X$. The formula is deduced from
the similar formula in $\mathbb{H}^2$.
\end{proof}
Equivalently, in the basis $\{X,R\}$, parallel transport is given
by the $2\times 2$ orthogonal matrix
\[ \rho^{-2}_{\rm ff}\,\begin{bmatrix} (x+1)^2-r^2 &
2r(x+1)\\-2r(x+1)& (x+1)^2-r^2
\end{bmatrix}. \]  We decompose a vector
$V=a_0X+\sum_{j=1}^d a_j e_j$ into its tangent, respectively
orthogonal components to $\Pi_m$ as follows:
\[V=a_0X+\left(\sum_{j=1}^d a_j \langle e_j,R\rangle\right) R
+\sum_{j=1}^d a_j T_j.\]
Since $T_j$ are parallel, the orthogonal component is constant during parallel
transport.
Using \eqref{umrn} and Lemma \ref{ptp}, we write
\begin{align*}\tau(m,m')(a_0X+\sum_{j=1}^d a_j e_j)=&a_0\tau(X)+\sum_{j=1}^d a_j
\langle e_j,R\rangle \tau(R)+\sum_{j=1}^d a_j T_j\\
=&a_0\frac{(x+x')^2-r^2}{\rho_{\rm ff}^2} X +2a_0 \frac{r(x+x')}{\rho_{\rm ff}^2} R +\sum_{j=1}^d a_j (e_j-\langle e_j,R\rangle R)\\
&+\sum_{j=1}^d a_j \langle e_j,R\rangle \left(-\frac{2r(x+x')}{\rho_{\rm ff}^2}X
+\frac{(x+x')^2-r^2}{\rho_{\rm ff}^2} R \right)
\end{align*}
from which the proposition follows, since
$(\rho_{\ff}, (y-y')/\rho_{\ff},
y',x/\rho_{\ff}, x'/\rho_{\ff})$
are smooth coordinates on the blow-up $\hh^{d+1}\x_0\hh^{d+1}$.
\end{proof}

The above oriented basis $\{X,e_1,\ldots,e_d\}$ of $T\mathbb{H}^{d+1}$
extends smoothly
to the boundary $\{0\}\times \rz^d$ of $\overline{\mathbb H}^{d+1}$ as
an orthonormal basis of the zero tangent bundle with respect to
the hyperbolic metric.
Therefore, by the Proposition above, $\tau(m,m')$ is a smooth section on ${\overline{\mathbb H}^{d+1}}
\times_0\overline{\mathbb H}^{d+1}$ in the pull-back vector bundle
\[\beta^* (\pi_1^* {}^0 T{\overline{\mathbb H}^{d+1}}
\boxtimes \pi_2^* {}^0 T^*{\overline{\mathbb H}^{d+1}}).\]

The orthonormal frame bundle $P_{\mathrm{SO}}$ of the $0$-tangent
bundle is trivialized over ${\overline{\mathbb H}^{d+1}}$ by the
frame $p=\{X,e_1,\ldots,e_d\}$, therefore the (unique) spin
structure $P_{\mathrm{Spin}}$ is identified with
${\overline{\mathbb H}^{d+1}}\times \mathrm{Spin}(d+1)$. Denote by
$\tilde{p}$ one of the lifts of $p$ to $P_{\mathrm{Spin}}$. By
definition of the lifted connection, parallel transport in
$P_{\mathrm{Spin}}$ of the section $\tilde{p}$ along the geodesic
from $m$ to $m'$ is $\tilde{p}\, U(m,m')$, where $U(m,m')$ the
unique lift of the $\mathrm{SO}(d+1)$-valued function $\tau(m,m')$
to the $\mathrm{Spin}(d+1)$ group, starting at the identity for
$m=m'$. Thus, parallel transport of a constant section $\sigma$
(with respect to the trivialization $\tilde{p}$) in the spinor
bundle is simply
\[\tau_m^{m'}[\tilde{p},\sigma]=[\tilde{p}\,U(m,m'),\sigma]=
[\tilde{p},U(m,m')\sigma] \]
where multiplication in the last term is the spinor representation.
By abuse of notation we write $U(m,m')$ for $\tau_m^{m'}$.

\begin{proposition} \label{ptspin}
Let $m=(x,y),m'=(x',y')\in \overline{\mathbb H}^{d+1}$. In the above
trivialization of the spinor bundle, parallel transport takes the form
\[U(m,m')=\frac{x+x'}{\rho}-\frac{r}{\rho}\cl(X)\cl(R).\]
\end{proposition}
\begin{proof}
We view the $\mathrm{Spin}(d+1)$ group inside the Clifford algebra
as the group generated by even Clifford products of unit vectors.
The projection $\pi:\mathrm{Spin}(d+1)\to \mathrm{SO}(d+1)$ is
given by the adjoint action in the Clifford algebra on vectors:
\[\pi(c)(V):=cVc^{-1},\]
the kernel being precisely $\{\pm 1\}$. We must therefore examine
the adjoint action of $A(m,m'):=\frac{x+x'}{\rho_{\rm
ff}}-\frac{r}{\rho_{\rm ff}}\cl(X)\cl(R)$ on ${}^0 T{\overline{\mathbb H}^{d+1}}$.
Note that any Clifford element of the form $\alpha+\beta \cl(X)\cl(R)$
with $\alpha^2+\beta^2=1$ belongs to the
Spin group. Next, $A^{-1}(m,m')=\frac{x+x'}{\rho_{\rm
ff}}+\frac{r}{\rho_{\rm ff}}\cl(X)\cl(R)$ so
\[\pi(A(m,m'))X=\left(\frac{(x+x')^2}{\rho_{\rm ff}^2}
-\frac{r^2}{\rho_{\rm ff}^2}\right) X
-2\frac{(x+x')r}{\rho_{\rm ff}^2}R\]
which coincides with the
action of $\tau(m,m')$ on $X$ from Proposition \ref{partr}.
Similarly, for the vector fields $T_j$ from the proof of Proposition
\ref{partr} we have $\pi(A(m,m'))T_j=T_j=\tau(m,m')T_j.$
Thus $\pi(A(m,m'))=\tau(m,m')$. The proof is finished by noting that
$A(m,m')$ was normalized so that $A(m,m)=1$.
\end{proof}

\begin{corollary}\label{upp'}
Let $m'=(1,0)$, $m=(0,r\omega)$. In the limit $r\to\infty$, the
parallel transport $U(m,m')$ tends to $-\cl(X)\cl(R)$.
\end{corollary}

\bibliographystyle{plain}
\def\cprime{$'$} \def\polhk#1{\setbox0=\hbox{#1}{\ooalign{\hidewidth
  \lower1.5ex\hbox{`}
\crcr\unhbox0}}} \providecommand{\bysame}{\leavevmode\hbox
to3em{\hrulefill}\thinspace}
\providecommand{\MR}{\relax\ifhmode\unskip\space\fi MR }
\providecommand{\MRhref}[2]{%
  \href{http://www.ams.org/mathscinet-getitem?mr=#1}{#2}
} \providecommand{\href}[2]{#2}

\end{document}